\documentclass[12pt,oneside,reqno]{amsart}
\usepackage{}
\usepackage{amssymb}
\usepackage[colorlinks=true, linkcolor=blue, citecolor=red]{hyperref}
\usepackage{bbm}
\usepackage{cases}
\usepackage{amsmath}
\usepackage{graphicx}
\usepackage{mathrsfs}
\usepackage{stmaryrd}
\usepackage{color}
\usepackage{soul}
\usepackage[dvipsnames]{xcolor}
\usepackage{amsfonts}
\usepackage{enumerate,amsmath,amssymb,amsthm}

\numberwithin{equation}{section}

\newcommand{\be}{\begin{eqnarray}}
\newcommand{\mE}{\end{eqnarray}}
\newcommand{\ce}{\begin{eqnarray*}}
\newcommand{\de}{\end{eqnarray*}}
\newtheorem{theorem}{Theorem}[section]
\newtheorem{lemma}[theorem]{Lemma}
\newtheorem{remark}[theorem]{Remark}
\newtheorem{definition}[theorem]{Definition}
\newtheorem{proposition}[theorem]{Proposition}
\newtheorem{example}[theorem]{Example}
\newtheorem{corollary}[theorem]{Corollary}

\def\eps{\varepsilon}
\def\up{\Upsilon}

\def\p{\partial}

\def\[{{\Big[}}
\def\]{{\Big]}}
\def\<{{\langle}}
\def\>{{\rangle}}
\def\({{\Big(}}
\def\){{\Big)}}

\def\bx{{\mathbf{x}}}

\def\dif{{\mathord{{\rm d}}}}

\def\no{\nonumber}
\def\={&\!\!=\!\!&}
\def\bt{\begin{theorem}}
\def\et{\end{theorem}}
\def\bl{\begin{lemma}}
\def\el{\end{lemma}}
\def\br{\begin{remark}}
\def\er{\end{remark}}

\def\bd{\begin{definition}}
\def\ed{\end{definition}}
\def\bp{\begin{proposition}}
\def\ep{\end{proposition}}
\def\bc{\begin{corollary}}
\def\ec{\end{corollary}}
\def\bx{\begin{example}}
\def\ex{\end{example}}

\def\cD{{\mathcal D}}

\def\cG{{\mathcal G}}

\def\cL{{\mathcal L}}

\def\cQ{{\mathcal Q}}
\def\cR{{\mathcal R}}

\def\cV{{\mathcal V}}

\def\mE{{\mathbb E}}

\def\mI{{\mathbb I}}

\def\mN{{\mathbb N}}

\def\mP{{\mathbb P}}

\def\mR{{\mathbb R}}

\def\sC{{\mathscr C}}

\def\sF{{\mathscr F}}
\def\sG{{\mathscr G}}

\def\sL{{\mathscr L}}

\def\sN{{\mathscr N}}

\def\sQ{{\mathscr Q}}

\def\sS{{\mathscr S}}

\def\sU{{\mathscr U}}
\def\sV{{\mathscr V}}

\def\geq{\geqslant}
\def\leq{\leqslant}

\allowdisplaybreaks

\begin{document}

\title{Averaging principle and normal deviations for multiscale stochastic  systems}


\date{}

\author{Michael R\"ockner \,\, and \,\, Longjie Xie}

\address{Michael R\"ockner:
	Fakult\"at f\"ur Mathematik, Universit\"at Bielefeld,
	33615, Bielefeld, Germany\\
	Email: roeckner@math.uni-bielefeld.de}

\address{Longjie Xie:
	School of Mathematics and Statistics, Jiangsu Normal University,
	Xuzhou, Jiangsu 221000, P. R. China\\
	Email: longjiexie@jsnu.edu.cn}

\thanks{ This work is supported by the DFG through CRC 1283, the Alexander-von-Humboldt foundation and NSFC (No. 12090011, 12071186, 11931004)}

\begin{abstract}
We study the asymptotic behavior for an  inhomogeneous multiscale  stochastic dynamical system  with  non-smooth coefficients. Depending on the averaging regime and the homogenization regime,  two   strong convergences in the averaging principle of functional law of large numbers type are established. Then we consider the small fluctuations of the  system around its average. Nine cases of functional central limit type theorems   are obtained. In particular,  even though the averaged equation for the original system is the same, the corresponding homogenization limit for the normal deviation can be quite different due to the difference in the   interactions between  the fast scales and the deviation scales. We provide  quite intuitive explanations for each case.
Furthermore, sharp rates  both for the strong convergences and the functional central limit theorems are obtained, and these convergences are shown to rely only on the regularity of the coefficients of the system with respect to the slow variable, and do not depend on their regularity with respect to the fast variable, which coincide with the intuition since in the limit equations the fast component has been totally averaged or homogenized out.

	\bigskip

	
	\vspace{2mm}
\end{abstract}

\maketitle

\tableofcontents

\section{Introduction}

Consider  the following fast-slow stochastic system in $\mR^{d_1+d_2}$:
\begin{equation} \label{sde00}
\left\{ \begin{aligned}
&\dif X^{\eps}_t =\eps^{-1}b(X^{\eps}_t,Y^{\eps}_t)\dif t+\eps^{-1/2}\sigma(X^{\eps}_t,Y_t^\eps)\dif W^1_t,\quad X^{\eps}_0=x\in\mR^{d_1},\\
&\dif Y^{\eps}_t =F(X^{\eps}_t, Y^{\eps}_t)\dif t+G(Y_t^\eps)\dif W_t^2,\quad\quad\quad\qquad\quad   Y^{\eps}_0=y\in\mR^{d_2},
\end{aligned} \right.
\end{equation}
where  $d_1,d_2\in\mN$, $b: \mR^{d_1}\times\mR^{d_2}\to\mR^{d_1}$, $F: \mR^{d_1}\times\mR^{d_2}\to\mR^{d_2}$,  $\sigma: \mR^{d_1}\times\mR^{d_2}\to\mR^{d_1}\otimes\mR^{d_1}$ and $G:  \mR^{d_2}\to\mR^{d_2}\otimes\mR^{d_2}$ are Borel measurable functions, $W^1_t$, $W^2_t$ are $d_1$, $d_2$-dimensional independent standard Brownian motions respectively, both defined on some probability space $(\Omega,\sF,\mP)$, and the small parameter $0<\eps\ll 1$ represents the separation  of  time scales  between   the fast motion $X_t^\eps$ (with time order $1/\eps$) and  the slow component $Y_t^\eps$.
Such multiscale models  appear frequently in many real world dynamical  systems. Typical examples   include climate weather interactions (see e.g. \cite{K,MTV}),  intracellular biochemical reactions (see e.g. \cite{BKRP,KK}), geophysical fluid flows (see e.g. \cite{GD}), stochastic volatility in finance (see e.g. \cite{FFK}),  etc. We refer the interested readers to the books \cite{Ku,PS} for a more  comprehensive overview.
In fact, as mentioned in \cite{SA}, almost all physical systems have a certain hierarchy in which not all components evolve at the same rate, and a mathematical description for such phenomena can be formulated by a singularly perturbed differential equation with a small parameter such as the one given in (\ref{sde00}).
Usually, the underlying system (\ref{sde00}) is difficult to deal with due to  the widely separated time scales and the cross interactions
between the  fast and slow   modes.
Thus a simplified equation
which governs the evolution of the system over a long time scale is highly desirable and is quite important for applications.

\vspace{2mm}
The intuitive idea for deriving such a simplified equation for system (\ref{sde00}) is based on the observation
that during the fast transients, the slow variable
remains ``constant" and by the time its changes become
noticeable, the fast variable has almost reached its
``quasi-steady state". Noting that after the natural time scaling $t\mapsto \eps t$, the process  $\tilde X_t^\eps:=X_{\eps t}^\eps$ satisfies
$$
\dif \tilde X_t^\eps=b(\tilde X^{\eps}_t,Y^{\eps}_{\eps t})\dif t+\sigma(\tilde X^{\eps}_t,Y_{\eps t}^\eps)\dif \tilde W^1_t,\quad \tilde X^{\eps}_0=x\in\mR^{d_1},
$$
where $\tilde W_t^1:=\eps^{-1/2}W^1_{\eps t}$ is a new Brownian motion.
Thus we need to consider the auxiliary process $X_t^y$ which is the solution of the following frozen stochastic differential equation (SDE for short): for fixed $y\in\mR^{d_2}$,
\begin{align}\label{sde2}
\dif X_t^y=b(X_t^y,y)\dif t+\sigma(X_t^y,y)\dif W_t^1,\quad X_0^y=x\in\mR^{d_1}.
\end{align}
The re-scaled  process $\tilde X_t^\eps$ will be asymptotically identical in distribution to $X_t^y$.
 Under certain recurrence conditions, the process $X_t^y$ admits a unique invariant measure $\mu^y(\dif x)$.
Then by averaging the coefficient with respect to parameters in the fast variable, the slow component $Y_t^\eps$ in system (\ref{sde00}) will converge  strongly (in the $L^2$-sense) as $\eps\to0$ to the solution of the following so-called averaged equation in $\mR^{d_2}$:
\begin{align}\label{sde11}
\dif \bar Y_t=\bar F(\bar Y_t)\dif t+G(\bar Y_t)\dif W_t^2,\quad\bar Y_0=y\in\mR^{d_2},
\end{align}
where the new averaged drift coefficient is defined by
$$
\bar F(y):=\!\int_{\mR^{d_1}}\!F(x,y)\mu^y(\dif x).
$$
The effective dynamics (\ref{sde11})  does not depend on the fast variable any more and thus is much simpler than
SDE (\ref{sde00}). This theory, known as the averaging principle, can be regarded as  classical  functional law of large numbers  and has been intensively studied in both the deterministic ($\sigma=G\equiv0$) and the stochastic context in the past decades,  see e.g. \cite{CFKM,D-L1,D-L2,GR,HLi,K1,KY2,Ki,Li}  and the references therein, see also \cite{Ce,CF,CL,WR} for similar results concerning stochastic partial differential equations (SPDEs for short).  Note that the diffusion coefficient $G$ in  SDE (\ref{sde00}) does not depend on the fast  variable $x$, otherwise, the strong convergence  may not be true (see e.g. \cite{Ki2}). Meanwhile, as a rule the averaging method requires certain regularities of the  coefficients of the original system (\ref{sde00})   to guarantee the above convergence, and we point out that  all the aforementioned papers assumed at least local Lipschitz continuity of all the coefficients.
For the averaging principle of SDE (\ref{sde00}) with non-smooth coefficients, we refer to \cite{RSX,V0}.

\vspace{2mm}
However, the effective equation (\ref{sde11}) is valid only in the limiting sense, and the time scale separation is never infinite in reality. For small but positive $\eps$, the slow process $Y_t^\eps$ will experience fluctuations around its averaged motion $\bar Y_t$.
To leading order, these fluctuations can be captured by characterizing the asymptotic behavior of the normalized difference
$$
Z_t^\eps:=\frac{Y_t^\eps-\bar Y_t}{\sqrt{\eps}}
$$
as $\eps$ tends to 0. Under extra regularity assumptions on the coefficients and when $G\equiv\mI_{d_2}$ (the $d_2\times d_2$ identity matrix), the deviation process $Z_t^\eps$ is known to converge weakly (i.e., in distribution) towards an Ornstein-Uhlenbeck type process $\bar Z_t$ with $\bar Z_t$ satisfying the following (linear) SDE in $\mR^{d_2}$:
\begin{align}\label{clt}
\dif \bar Z_t=\nabla_y\bar F(\bar Y_t)\bar Z_t\dif t+\zeta(\bar Y_t)\dif \tilde W_t,
\end{align}
where $\bar Y_t$ solves the averaged equation (\ref{sde11}), $\tilde W_t$ is another standard Brownian motion, and the new diffusion coefficient is given by
\begin{align}\label{cltd}
\zeta(y):=\sqrt{\int_0^\infty\!\!\int_{\mR^{d_1}}\mE\big[F(X_t^y(x),y)-\bar F(y)\big]\big[F(x,y)-\bar F(y)\big]^*\mu^y(\dif x)\dif t}.
\end{align}
Such result, also known as the Gaussian  approximation, is an analogue of the functional central limit theorem of Donsker. We refer the readers to the fundamental paper by Khasminskii \cite{K1}, see also \cite{HL,KV,KT,PTW,S1}  for further developments and \cite{Ce2,HP2}  for the corresponding central limit theorem type results for multiscale SPDEs. Besides having intrinsic interest, the central limit theorem  is also useful in applications. In particular, we can get the formal asymptotic expansion
$$
Y_t^\eps\stackrel{\cD}{\approx} \bar Y_t+\sqrt{\eps}\bar Z_t,
$$
where $\stackrel{\cD}{\approx}$ means approximate equality of probability distributions. Such  expansion has been introduced in
the context of stochastic climate models. In physics
this is also called the Van Kampen's  approximation (see e.g. \cite{Ar}),
which provides  better approximations for the original system (\ref{sde00}), see also \cite{BK,Ki2} and the references therein.
We  mention that other limit theorems for SDE (\ref{sde00}) have also been widely studied in the literature,
see e.g. \cite{BGTV,BDG,DS,Pu,S2,Ve} for the large deviations and \cite{Gu,Gu2,MS} for the moderate deviations.

\vspace{2mm}
In this paper, we study a  broader class of system, i.e.,
consider the following inhomogeneous multiscale SDE in $\mR^{d_1+d_2}$:
\begin{equation} \label{sde0}
\left\{ \begin{aligned}
&\dif X^{\eps}_t =\alpha_\eps^{-2}b(X^{\eps}_t,Y^{\eps}_t)\dif t+\beta_\eps^{-1}c(X^{\eps}_t,Y^{\eps}_t)\dif t+\alpha_\eps^{-1}\sigma(X^{\eps}_t,Y_t^\eps)\dif W^{1}_t,\\
&\dif Y^{\eps}_t =F(t,X^{\eps}_t, Y^{\eps}_t)\dif t+\gamma_\eps^{-1}H(t,X_t^\eps,Y_t^\eps)\dif t+G(t,Y_t^\eps)\dif W^{2}_t,\\
&X^{\eps}_0=x\in\mR^{d_1},\quad Y^{\eps}_0=y\in\mR^{d_2},
\end{aligned} \right.
\end{equation}
where the small parameters $\alpha_\eps,\beta_\eps,\gamma_\eps\to0$ as $\eps\to0$, and without loss of generality  we  assume $\alpha_\eps^2/\beta_\eps\to0$  as $\eps\to0$, and conventionally take $\beta_\eps\equiv1$ when $c\equiv0$ and $\gamma_\eps\equiv1$ when $H\equiv0$. The infinitesimal generator $\sL_\eps$ corresponding to  system (\ref{sde0}) has the form
$$
\sL_\eps:=\frac{1}{\alpha_\eps^{2}}\sL_0(x,y)+\frac{1}{\beta_\eps}\sL_3(x,y)+\frac{1}{\gamma_\eps}\sL_2(t,x,y)+\sL_1(t,x,y),
$$
where $\sL_0(x,y)$ is given by
\begin{align}
\sL_0:=\sL_0(x,y):=\sum_{i,j=1}^{d_1} a^{ij}(x,y)\frac{\p^2}{\p x_i\p x_j}+\sum_{i=1}^{d_1}b^i(x,y)\frac{\p}{\p x_i}\label{l0}
\end{align}
with $a(x,y):=\sigma\sigma^*(x,y)/2$ (where $\sigma^*$ denotes the transpose of $\sigma$), and
\begin{equation} \label{lll}
\begin{aligned}
&\sL_3:=\sL_3(x,y):=\sum_{i=1}^{d_1}c^i(x,y)\frac{\p}{\p x_i},\\
&\sL_2:=\sL_2(t,x,y):=\sum_{i=1}^{d_2}H^i(t,x,y)\frac{\p}{\p y_i},\\
&\sL_1:=\sL_1(t,x,y):=\sum_{i,j}^{d_2} \cG^{ij}(t,y)\frac{\p^2}{\p y_i\p y_j}+\sum_{i=1}^{d_2}F^i(t,x,y)\frac{\p}{\p y_i}\\
\end{aligned}
\end{equation}
with $\cG(t,y)=GG^*(t,y)/2$. Note that there exist two time scales in the fast motion $X_t^\eps$ and even the slow process $Y_t^\eps$ has a fast
varying component.
This is known to be important, in particular, for applications in the homogenization of second order parabolic and elliptic equations with singularly perturbed terms, which has its own interest in the theory of PDEs, see e.g. \cite{HP,KRS,Par} and \cite[Chapter IV]{Fr}.
The study of such generalized systems (\ref{sde0}) was first carried out by Papanicolaou, Stroock and Varadhan \cite{PSV} for a compact state space for the fast component and  time-independent coefficients when  $\alpha_\eps=\beta_\eps=\gamma_\eps$, see also \cite{Ba} for a similar result in terms of PDEs. Later, a  non-compact  homogeneous case with $c\equiv0$ and $\alpha_\eps=\gamma_\eps$ was studied in  \cite{P-V,P-V2} by using the method of the martingale problem and in \cite{KY} by the asymptotic expansion approach, see also  \cite{RX} for a more systematical study. However, all these papers concern only weak convergence of the slow process $Y_t^\eps$ in SDE (\ref{sde0}). In \cite{S1}, the convergence in probability of $Y_t^\eps$ and its small fluctuations around the averaged motion were studied in a particular homogeneous case where $\alpha_\eps=\delta/\sqrt{\eps}$, $\beta_\eps=\delta$,  $\gamma_\eps=\delta/\eps$ and with small noise perturbations, i.e., with $G$ replaced by $\sqrt{\eps}G$ in SDE (\ref{sde0}).  To the best of our knowledge, no strong convergence and functional central limit theorems type results as well as rates of convergence  in terms of $\eps\to0$ for SDE (\ref{sde0}) in the general case have been obtained so far.

\vspace{2mm}
We shall first study the strong convergence in the averaging principle for SDE (\ref{sde0}) with non-smooth coefficients. The main result   is given by {\bf Theorem \ref{main1}} below. It turns out that  depending on the orders how $\alpha_\eps,\beta_\eps,\gamma_\eps$ go to zero, we need to distinguish two different regimes of interactions, which lead to two different asymptotic behaviors for system (\ref{sde0}) as $\eps\to0$, i.e.,
\begin{equation}\label{regime}
\left\{\begin{aligned} &\lim_{\eps\to0}\frac{\alpha_\eps}{\gamma_\eps}=0\quad\text{and}\quad\lim_{\eps\to0}\frac{\alpha_\eps^2}{\beta_\eps\gamma_\eps}=0,\qquad \text{Regime 1};\\
&\lim_{\eps\to0}\frac{\alpha_\eps}{\gamma_\eps}=0\quad\text{and}\quad\alpha_\eps^2=\beta_\eps\gamma_\eps,\qquad\quad\,\, \text{Regime 2}.
\end{aligned}\right.
\end{equation}
If $\alpha_\eps$ and $\alpha_\eps^2$ go to zero faster than $\gamma_\eps$ and $\beta_\eps\gamma_\eps$ respectively (Regime 1), we show that the averaged equation for system (\ref{sde0}) coincides with the traditional case which corresponds to $c=H\equiv0$; whereas if $\alpha_\eps$ goes to zero faster than $\gamma_\eps$, while $\alpha_\eps^2$ and $\beta_\eps\gamma_\eps$ are of the same order (Regime 2) (which means that the term $\beta_\eps^{-1}c(X_t^\eps,Y_t^\eps)$ is varying fast enough),
then the averaging effect of term $c$ and the homogenization effect of term $H$ will occur in the effective dynamics.  Furthermore, unlike most previous publications (see e.g. \cite{KY,PSV,P-V2,S1}), we will mainly focus on the impact  of noises on the averaging principle for system (\ref{sde0}).
Namely,
we  prove that for non-degenerate noises, the averaging principle holds for system (\ref{sde0}) with only H\"older continuous drifts (the corresponding deterministic system  can even be ill-posed under such weak conditions on the coefficients), and the convergence in the averaging principle  relies only on the regularities of the coefficients with respect to the slow component ($y$-variable), and does not depend on their regularities with respect to the fast term ($x$-variable).  This coincides with the intuition, since in the limit equation the fast variable has been totally averaged out. See {\bf Remark \ref{ma1}} and {\bf Remark \ref{ma5}}  below for more detailed explanations and comparisons of our results with the  previous literature on the subject.

\vspace{2mm}
Our method to prove the strong convergence in Regime 1 and Regime 2 is unified and rather simple as we do not need  the classical time discretisation procedure, which is commonly used in the literature to prove the averaging principle (see e.g. \cite{Ce,K1,Li,V0,WR}). Two ingredients are crucial in our proof: Zvonkin's transformation  and the Poisson equation in the whole space. First of all, due to the low regularity of the coefficients, we shall use Zvonkin's argument as in \cite{RSX,V0} to transform the equations for $Y_t^\eps$ and its average into new ones with better coefficients.
Then we employ the result of Poisson equation established in \cite{RX} to prove the strong convergence for system (\ref{sde0}). In both regimes, rates of convergence are also obtained as easy by-products of our arguments. The convergence rates are known to be important for the analysis of numerical schemes for multiscale systems, see e.g. \cite{Br3,Br1,EL,KVa}. Moreover, it will play a crucial  role below for us to study the homogenization behavior for the fluctuations of $Y_t^\eps$ around its average in determining the deviation scales, which in turn implies that  the strong convergence rates obtained here are optimal.

\vspace{2mm}
After the strong convergence of the functional  law of large numbers type is established, we then proceed  to study the functional central limit theorem for system (\ref{sde0}). More precisely, we  will be interested in the asymptotic behavior for the normal deviations  of $Y_t^\eps$ from its averaged motion $\bar Y_t^k$ ($k=1,2$ which correspond to Regime 1 and Regime 2 in (\ref{regime})), i.e., to identify the limit of the normalized  difference
$$
Z_t^{k,\eps}:=\frac{Y_t^\eps-\bar Y_t^k}{\eta_\eps}
$$
with proper choice of deviation scale $\eta_\eps$ such that $\eta_\eps\to0$ as $\eps$ tends to 0. It turns out that the asymptotic limit for $Z_t^{k,\eps}$ is strongly linked to the interactions of the fast scales as well as the deviation scale $\eta_\eps$. Even though the law of large numbers type limit for $Y_t^\eps$ is the same, the homogenization behavior in the functional central limit theorems for the deviation process $Z_t^{k,\eps}$ can be quite different.  We need to distinguish three main cases:

\vspace{2mm}
\noindent{\bf Case 0:} $H\equiv 0$ in SDE (\ref{sde0}), i.e., there is no homogenization term in the slow equation. Note that even in this case, the system is still more general than the traditional ones due to the existence of the extra term $\beta_\eps^{-1}c(X_t^\eps,Y_t^{\eps})$ in the fast motion. We shall show that the limit equation  for the deviation process could be given in terms of the drift $c$ and the solution for an auxiliary Poisson equation involving the drift $F$.

\vspace{1mm}
\noindent{\bf Case 1:} $H\neq0$ with Regime 1 described in (\ref{regime}). In this case, we shall show that the averaging effect of the drift $c$  and the homogenization effect of the fast term $H$ may arise in the limit equation for $Z_t^{1,\eps}$, while the effects involving the drift $F$ as in Caes 0 will not  appear any more.

\vspace{1mm}
\noindent{\bf Case 2:} $H\neq0$ with Regime 2 described in (\ref{regime}). As mentioned before, in this case homogenization  has already occurred even  in the averaged equation for $Y_t^\eps$ (see (\ref{barf}) below). Thus we shall show that the second order homogenization  involving the drift $c$ and $H$ may arise in the limit equation for $Z_t^{2,\eps}$. Furthermore, the effects involving  the drift $F$ as in Caes 0 may occur again.

\vspace{2mm}
The main results are given by {\bf Theorem \ref{main3}}, {\bf Theorem \ref{main2}} and {\bf Theorem \ref{main4}}, respectively. Moreover, it is interesting to find  that in each case, depending on the choice of the deviation scale we still get a different   limit behavior  for $Z_t^{k,\eps}$. Several new terms will appear which seem  to be never observed before in the literature. In certain regimes the limit equation for $Z_t^{k,\eps}$ can even  be given by a linear random ordinary differential equation (ODE for short), while in certain situations an extra Gaussian part appears and the limit equation will be given by a linear SDE. We shall provide some quite intuitive explanations for each respective result, see  {\bf Remark \ref{ma2}}, {\bf Remark \ref{ma3}} and  {\bf Remark \ref{m4}} below. In particular, our results  lead to a deep understanding of the effects and the interactions between the extra averaging  term $\beta_\eps^{-1}c(X_t^\eps,Y_t^{\eps})$ and the homogenization term $\gamma_\eps^{-1}H(t,X_t^\eps,Y_t^{\eps})$ in system (\ref{sde0}).

\vspace{2mm}
We will  prove the functional central limit theorem for system (\ref{sde0}) in each regime in a very robust and unified way. Our method relies only on  the technique of Poisson equation, and neither involves an extra time discretisation procedure, nor martingale problem or tightness arguments (see e.g. \cite{Ce2,K1,KV,S1,WR}) and thus is  quite simple. Moreover, the conditions on the coefficients are  weaker  than in the known results in the literature even in the traditional case (i.e., $c=H\equiv0$), and rates of convergence are also obtained,  which we believe are rather sharp.  Furthermore, it will be pretty clear from our approach that which parts should be the leading terms in the fluctuations (whose effects arise in the homogenization procedure), which parts should be the lower order terms (whose fluctuations go to zero eventually) and what the deviation scales $\eta_\eps$ should be  in each regime  in order to observe non-trivial behavior  for the limits.

\vspace{2mm}
The rest of this paper proceeds as follows. In Section 2 we state our main results.  Section 3 is devoted to the preparation of the main tools that shall be used to prove the results. Then we prove Theorem \ref{main1} in Section 4,  Theorem \ref{main3} in Section 5, and Theorem \ref{main2} and Theorem \ref{main4} in Section 6, respectively.
Throughout our paper, we use the following convention: $C$ and $c$ with or without subscripts will denote positive constants, whose values may change in different places, and whose dependence on parameters can be traced from the calculations. Given a function space, the subscript $b$ will stand for boundness, while the subscript $p$ stands for polynomial growth in the $x$ variable. To be more precise, for a function $f(t,x,y,z)$ defined on $\mR_+\times\mR^{d_1}\times\mR^{d_2+d_2}$, by $f\in L^\infty_p:=L^{\infty}_p(\mR_+\times\mR^{d_1}\times\mR^{d_2+d_2})$ we mean there exist constants $C, m>0$ such that
 $$
 |f(t,x,y,z)|\leq C(1+|x|^m),\quad\forall t>0, x\in\mR^{d_1}, y, z\in\mR^{d_2},
 $$
 and $ C_p^{\gamma,\delta,\vartheta,\eta}:=C_{t,x,y,z}^{\gamma,\delta,\vartheta,\eta}(\mR_+\times\mR^{d_1}\times\mR^{d_2+d_2})$ with $0<\delta\leq 1$ denotes the space of all functions $f$ such that  for every fixed $x\in\mR^{d_1}$, $f(\cdot,x,\cdot,\cdot)\in C_b^{\gamma,\vartheta,\eta}(\mR_+\times\mR^{d_2+d_2})$ and for any $(t,y,z)\in\mR_+\times\mR^{d_2+d_2}$ and $x_1,x_2\in\mR^{d_1}$,
 \begin{align*}
 |f(t,x_1,y,z)-f(t,x_2,y,z)|\leq C|x_1-x_2|^\delta\big(1+|x_1|^m+|x_2|^m\big).
 \end{align*}

\section{Statement of main results}

\subsection{Strong convergence: functional law of large numbers}

Let us first introduce some basic assumptions. Throughout this paper, we shall always assume the following non-degeneracy conditions on the diffusion coefficients:

\vspace{2mm}
\noindent{\bf (A$_\sigma$):} the coefficient $a=\sigma\sigma^*/2$ is non-degenerate in $x$ uniformly with respect to $y$, i.e., \\
\indent\qquad\,\!\! there exists a $\lambda>1$ such that  for any $y\in\mR^{d_2}$,
$$
\lambda^{-1}|\xi|^2\leq |\sigma^*(x,y)\xi|^2\leq \lambda|\xi|^2,\ \ \forall\xi\in\mR^{d_1}.
$$

\vspace{1mm}
\noindent{\bf (A$_{G}$):} the coefficient $\cG=GG^*/2$ is non-degenerate in $y$ uniformly with respect to $t$, i.e., \\ \indent\qquad\,\!\! there exists a $\lambda>1$ such that for any $t>0$,
$$
\lambda^{-1}|\xi|^2\leq |G^*(t,y)\xi|^2\leq\lambda|\xi|^2,\ \ \forall\xi\in\mR^{d_2}.
$$
Recall that the frozen equation is given by SDE (\ref{sde2}). We make the following very weak recurrence assumption on the drift $b$  to ensure the existence of an  invariant measure $\mu^y(\dif x)$ for $X_t^y$ (cf. \cite{Ve5}):

\vspace{2mm}
\noindent{\bf (A$_b$):}\qquad\qquad\qquad\quad $\lim_{|x|\to\infty}\sup_y \<x,b(x,y)\>=-\infty$.

\vspace{2mm}
\noindent Note that the drift $c$ in SDE (\ref{sde0}) is not involved in the frozen equation. We need the following additional condition on $c$ to ensure the non-explosion of the solution $X_t^\eps$: for $\eps>0$ small enough, it holds that
\begin{align}\label{ac}
\lim_{|x|\to\infty}\sup_y \<x,b(x,y)+\eps c(x,y)\>=-\infty.
\end{align}
Concerning the fast term in the slow component of SDE (\ref{sde0}), it is natural to make the following assumption:

\vspace{2mm}
\noindent{\bf (A$_H$):} the drift $H$ is centered, i.e.,
\begin{align}\label{cen}
\int_{\mR^{d_1}}H(t,x,y)\mu^y(\dif x)=0,\quad\forall (t,y)\in\mR_+\times\mR^{d_2},
\end{align}
where $\mu^y(\dif x)$ is the invariant measure for SDE (\ref{sde2}).

\vspace{2mm}
Under (\ref{cen}) and according to Theorem \ref{popde} below, there exists a unique solution $\Phi(t,x,y)$ to the following Poisson equation in $\mR^{d_1}$:
\begin{align}
\sL_0(x,y)\Phi(t,x,y)=-H(t,x,y),\quad x\in\mR^{d_1},  \label{pde1}
\end{align}
where $\sL_0(x,y)$ is given by (\ref{l0}), and $(t,y)\in\mR_+\times\mR^{d_2}$ are regarded as parameters. We introduce the new averaged drift coefficients by
\begin{align}\label{barf}
\begin{split}
\bar F_1(t,y)&:=\int_{\mR^{d_1}}\!F(t,x,y)\mu^{y}(\dif x);\\
\bar F_2(t,y)&:=\int_{\mR^{d_1}}\!\big[F(t,x,y)+c(x,y)\cdot\nabla_x\Phi(t,x,y)\big]\mu^{y}(\dif x),
\end{split}
\end{align}
which correspond to Regime 1 and Regime 2 described in (\ref{regime}), respectively. Then the precise formulation of the averaged equation for SDE (\ref{sde0}) is as follows: for $k=1,2$,
\begin{align}\label{sde000}
\dif \bar Y^k_t=\bar F_k(t, \bar Y^k_t)\dif t+G(t,\bar Y^k_t)\dif W^2_t,\quad \bar Y^k_0=y\in\mR^{d_2}.
\end{align}

The following is the first main result of this paper.

\bt[Strong convergence]\label{main1}
Let {\bf (A$_\sigma$)}, {\bf (A$_b$)}, {\bf (A$_{G}$)}, {\bf (A$_H$)} and (\ref{ac}) hold, $\delta,\vartheta>0$ and  $\lim_{\eps\to0}\alpha_\eps^\vartheta/\gamma_\eps=0$. Then for any $T>0$ and every $q\geq 1$,

\vspace{1mm}
\noindent(i) (Regime 1) if  $b, \sigma\in C_b^{\delta,\vartheta}$, $F, H\in C_p^{\vartheta/2, \delta,\vartheta}$, $G\in C_b^{\vartheta/2, 1}$ and  $c\in L^\infty_p$, we have
\begin{align}\label{stro1}
\sup_{t\in[0,T]}\mE|Y_t^\eps-\bar Y^1_t|^q\leq C_T\Big(\frac{\alpha_\eps^{\vartheta\wedge1}}{\gamma_\eps}+\frac{\alpha_\eps^2}{\beta_\eps\gamma_\eps}\Big)^q;
\end{align}

\noindent(ii) (Regime 2) if  $b, \sigma\in C_b^{\delta,\vartheta}$, $F, H\in C_p^{\vartheta/2, \delta,\vartheta}$, $G\in C_b^{\vartheta/2, 1}$ and $c\in C_p^{\delta,\vartheta}$, we have
\begin{align}\label{stro2}
\sup_{t\in[0,T]}\mE|Y_t^\eps-\bar Y^2_t|^q\leq C_T\Big(\frac{\alpha_\eps^{\vartheta\wedge1}}{\gamma_\eps}+\frac{\alpha_\eps^2}{\beta_\eps}\Big)^q,
\end{align}
where for $k=1, 2$, $\bar Y_t^k$ is the unique strong solution  for SDE (\ref{sde000}), and $C_T>0$ is a constant independent of $\delta, \eps$.
\et

Let us list some important comments regarding the above result.

\br\label{ma1}

(i) [Non-smooth coefficients]. Note that homogenization occurs in Regime 2 and an additional drift part appears in the limit. In both regimes, we do not make any Lipschitz-type assumptions on the drift coefficients $b$, $c$, $F$ and $H$.  We mention that if $\sigma=0$ or $G=0$,  the corresponding deterministic system may even be ill-posed for only H\"older continuous coefficients. This reflects the regularization effects of the noises. In fact, under our assumptions we have for every $k=1,2$, $\bar F_k\in C_b^{\vartheta/2,\vartheta}$.  Thus, the strong well-posedness for system (\ref{sde0}) and SDE (\ref{sde000}) follows by \cite{Ver} or \cite[Theorem 1.3]{Zh1}.

We also point out that the above results still hold in the small noise perturbation case, i.e., with $G$ replaced by $\lambda_\eps G$, where $\lambda_\eps\to0$ as $\eps\to0$. Then we need to assume the coefficients $b$, $c$, $F$ and $H$ to be  Lipschitz continuous with respect to the $y$ variable in order to ensure the well-posedness of the averaged system. For the sake of simplicity, we do not
deal with this setting in the present article.

\vspace{1mm}
(ii)  [Dependence of convergence]. In both regimes, the convergence rates do not depend on the index $\delta$. This suggests that the convergence in the averaging principle relies only on the regularities of the coefficients in the original system with respect to the time variable and the $y$ (slow) variable, and does not depend on their regularities with respect to the x (fast) variable.

\vspace{1mm}
(iii) [Sharp rates]. The traditional result can be viewed as a particular case of Regime 1 by taking  $c=H\equiv0$ (i.e., $\beta_\eps=\gamma_\eps\equiv1$). In this case, our result implies  that the averaging principle holds for SDE (\ref{sde0}) with  a strong convergence rate $\alpha_\eps^\vartheta$ when the coefficients are $\vartheta$-H\"older continuous. This order is known to be optimal when $\vartheta=1$. In the general case and when $\vartheta=1$,
estimate (\ref{stro1}) means that in Regime 1, the slow process $Y_t^\eps$ will converge to $\bar Y_t^1$ strongly with rate $\tfrac{\alpha_\eps}{\gamma_\eps}+\frac{\alpha_\eps^2}{\beta_\eps\gamma_\eps}$, while estimate  (\ref{stro2}) suggests that in Regime 2, $Y_t^\eps$  converges to $\bar Y_t^2$ strongly with order $\frac{\alpha_\eps}{\gamma_\eps}+\frac{\alpha_\eps^2}{\beta_\eps}$. We shall show that these rates are also optimal by studying the respective functional central limit theorems.
\er

\subsection{Functional central limit theorem: without homogenization}

We first consider  SDE (\ref{sde0}) with $H\equiv0$, i.e., there is no fast term in the slow component. To avoid confusion of  notations, we shall denote by $Y_t^{0,\eps}$ the slow process. More precisely, consider
\begin{equation} \label{sde99}
\left\{ \begin{aligned}
&\dif X^{\eps}_t =\alpha_\eps^{-2}b(X^{\eps}_t,Y^{0,\eps}_t)\dif t+\beta_\eps^{-1}c(X^{\eps}_t,Y^{0,\eps}_t)\dif t+\alpha_\eps^{-1}\sigma(X^{\eps}_t,Y_t^{0,\eps})\dif W^{1}_t,\\
&\dif Y^{0,\eps}_t =F(t,X^{\eps}_t, Y^{0,\eps}_t)\dif t+G(t,Y_t^{0,\eps})\dif W^{2}_t,\\
&X^{\eps}_0=x\in\mR^{d_1},\quad Y^{0,\eps}_0=y\in\mR^{d_2}.
\end{aligned} \right.
\end{equation}
Note that even in this case, the above system  is still broader than the traditional ones due to the existence of the extra  term $\beta_\eps^{-1}c(X_t^\eps,Y_t^{0,\eps})$ in the fast equation.
According to Theorem \ref{main1} (i) with $H\equiv0$ (then $\gamma_\eps\equiv1$) and $\vartheta=1$, the slow process $Y_t^{0,\eps}$ will converge to $\bar Y_t^1$ strongly with a best  possible rate $\alpha_\eps+\alpha_\eps^2/\beta_\eps$. We intend to study the small fluctuations of $Y_t^{0,\eps}$ from its average $\bar Y_t^1$, i.e., to characterize the asymptotic behavior of the normalized difference
\begin{align}\label{z0}
Z_t^{0,\eps}:=\frac{Y_t^{0,\eps}-\bar Y^1_t}{\eta_\eps}
\end{align}
with  proper deviation scale $\eta_\eps$ such that $\eta_\eps\to0$ as $\eps\to0$.  It turns out that the limit behavior for $Z_t^{0,\eps}$ is strongly linked to the deviation scale $\eta_\eps$. Formally, if $\alpha_\eps$ goes to 0 faster than $\alpha_\eps^2/\beta_\eps$,
then the convergence rate of $Y_t^{0,\eps}$ to $\bar Y^1_t$ is dominated by $\alpha_\eps^2/\beta_\eps$. Thus one needs time of order $\alpha_\eps^2/\beta_\eps$ to observe non-trivial behavior for $Z_t^{0,\eps}$; while if $\alpha_\eps$ is of the same order or lower order than $\alpha_\eps^2/\beta_\eps$, then   $Y_t^{0,\eps}$ will converge  to $\bar Y^1_t$ with rate $\alpha_\eps$ and we shall need deviation scale $\alpha_\eps$ to observe non-trivial homogenization effects.
Consequently, the natural choice of the deviation scale $\eta_\eps$ should be divided into the following three regimes:
\begin{equation}\label{regime4}
\left\{\begin{aligned}
&\eta_\eps=\frac{\alpha_\eps^2}{\beta_\eps}\qquad\!\text{and}\qquad\lim_{\eps\to0}\frac{\beta_\eps}{\alpha_\eps}=0,\qquad\qquad\qquad\!\!\! \text{Regime 0.1};\\
&\eta_\eps=\alpha_\eps\qquad\text{and}\qquad\lim_{\eps\to0}\frac{\alpha_\eps}{\beta_\eps}=0,\qquad\qquad\qquad\!\!\! \text{Regime 0.2};\\
&\eta_\eps=\alpha_\eps=\beta_\eps,\qquad\qquad\qquad\qquad\qquad\qquad\quad\quad\!\!\text{Regime 0.3}.
\end{aligned}\right.
\end{equation}
Such choices of $\eta_\eps$ will also appear to be natural from our proof procedure. Let $\Upsilon(t,x,y)$ be the unique solution of the following  Poisson equation in $\mR^{d_1}$:
\begin{align}
\sL_0(x,y)\Upsilon(t,x,y)=-\big[F(t,x,y)-\bar F_1(t,y)\big]:=-\tilde F(t,x,y), \label{pde11}
\end{align}
where $\sL_0(x,y)$ is defined by (\ref{l0}), $\bar F_1$ is given by (\ref{barf}), and $(t,y)\in\mR_+\times\mR^{d_2}$ are regarded as parameters. Define
\begin{align}
\overline{c\cdot\nabla_x\up}(t,y)&:=\int_{\mR^{d_1}}c(x,y)\cdot\nabla_x\up(t,x,y)\mu^y(\dif x),\label{cup}\\
\overline{\tilde F\cdot\up^*}(t,y)&:=\int_{\mR^{d_1}}\tilde F(t,x,y)\cdot\up^*(t,x,y)\mu^y(\dif x).\no
\end{align}
Then the limit processes $\bar Z^0_{\ell,t}$ ($\ell=1,2,3$) for $Z_t^{0,\eps}$ corresponding to Regime 0.1-Regime 0.3 in (\ref{regime4}) turn out to satisfy the following linear equations:
\begin{align}\label{sdez0}
\begin{split}
\dif \bar Z_{1,t}^{0}&=\nabla_y\bar F_1(t,\bar Y_t^1)\bar Z_{1,t}^{0}\dif t+\nabla_yG(t,\bar Y_t^1)\bar Z_{1,t}^{0}\dif W_t^2+\overline{c\cdot\nabla_x\up}(t,\bar Y_t^1)\dif t;\\
\dif \bar Z_{2,t}^{0}&=\nabla_y\bar F_1(t,\bar Y_t^1)\bar Z_{2,t}^{0}\dif t+\nabla_yG(t,\bar Y_t^1)\bar Z_{2,t}^{0}\dif W_t^2+\sqrt{\overline{\tilde F\cdot\up^*}(t,\bar Y_t^1)}\dif \tilde W_t;\\
\dif \bar Z_{3,t}^{0}&=\nabla_y\bar F_1(t,\bar Y_t^1)\bar Z_{3,t}^{0}\dif t+\nabla_yG(t,\bar Y_t^1)\bar Z_{3,t}^{0}\dif  W_t^2\\
&\quad+\overline{c\cdot\nabla_x\up}(t,\bar Y_t^1)\dif t+\sqrt{\overline{\tilde F\cdot\up^*}(t,\bar Y_t^1)}\dif \tilde W_t,
\end{split}
\end{align}
with initial data $\bar Z_{\ell,0}^{0}=0$, where $\bar Y_t^1$ is the unique strong solution for SDE (\ref{sde000}) with $k=1$, and $\tilde W_t$ is another Brownian motion independent of $W_t^2$.

Our first functional  central limit theorem type result  is as follows.

\bt[Central limit theorem]\label{main3}
Let {\bf (A$_\sigma$)}, {\bf (A$_b$)}, {\bf (A$_{G}$)}  and (\ref{ac}) hold,  $0<\delta,\vartheta\leq 1$. Then for any $T>0$ and every $\varphi\in C_b^{4}(\mR^{d_2})$,

\vspace{1mm}
\noindent(i) (Regime 0.1) if $b, \sigma\in C_b^{\delta,1+\vartheta}$, $F\in C_p^{(1+\vartheta)/2, \delta,1+\vartheta}$, $G\in C_b^{1/2, 1+\vartheta}$ and  $c\in C_p^{\delta,\vartheta}$, we have
\begin{align*}
\sup_{t\in[0,T]}\Big|\mE[\varphi(Z_t^{0,\eps})]-\mE[\varphi(\bar Z^0_{1,t})]\Big|\leq C_T\Big(\frac{\beta_\eps^2}{\alpha_\eps^2}+\frac{\alpha_\eps^{2\vartheta}}{\beta_\eps^\vartheta}\Big);
\end{align*}

\noindent(ii) (Regime 0.2) if $b, \sigma\in C_b^{\delta,1+\vartheta}$, $F\in C_p^{(1+\vartheta)/2, \delta,1+\vartheta}$, $G\in C_b^{1/2, 1+\vartheta}$ and  $c\in L^\infty_p$, we have
\begin{align*}
\sup_{t\in[0,T]}\Big|\mE[\varphi(Z_t^{0,\eps})]-\mE[\varphi(\bar Z^0_{2,t})]\Big|\leq C_T\Big(\frac{\alpha_\eps}{\beta_\eps}+\alpha_\eps^{\vartheta}\Big);
\end{align*}
\noindent(iii) (Regime 0.3) if $b, \sigma\in C_b^{\delta,1+\vartheta}$, $F\in C_p^{(1+\vartheta)/2, \delta,1+\vartheta}$, $G\in C_b^{1/2, 1+\vartheta}$ and  $c\in C_p^{\delta,\vartheta}$, we have
\begin{align*}
\sup_{t\in[0,T]}\Big|\mE[\varphi(Z_t^{0,\eps})]-\mE[\varphi(\bar Z^0_{3,t})]\Big|\leq C_T\,\alpha_\eps^{\vartheta},
\end{align*}
where for $\ell=1,2,3$, $\bar Z^0_{\ell,t}$ satisfy the linear equation (\ref{sdez0}), and $C_T>0$ is a constant independent of $\delta, \eps$.
\et

\br\label{ma2}
(i) By Theorem \ref{popde} below, we have
$$
\overline{\tilde F\cdot\up^*}(t,y)=\int_0^\infty\!\!\int_{\mR^{d_1}}\mE\big[\tilde F(t,x,y)\tilde F^*(t,X_t^y(x),y)\big]\mu^y(\dif x)\dif t.
$$
Thus, in view of (\ref{clt}) and (\ref{cltd}) the classical result can be viewed as a particular case of Regime 0.2 by taking the coefficients to be time-independent, $c\equiv0$ (then $\beta_\eps\equiv1$) and $G\equiv \mI_{d_2}$. Even in this case, our result is still new in the sense that the conditions on the coefficients are weaker and the rate of convergence (i.e., $\alpha_\eps^\vartheta$) is obtained, which  again depends only on the regularity of the coefficients with respect to the slow variable.

\vspace{1mm}
(ii) The above result reveals the effect of the extra fast term $\beta_\eps^{-1}c(X_t^\eps,Y_t^{0,\eps})$ in  system (\ref{sde99}): even though it does not play any role in the functional law of large numbers   for $Y_t^{0,\eps}$, it does affect the  deviations of $Y_t^{0,\eps}$ from $\bar Y_t^1$ through the functional central limit theorem. Note that both in Regime 0.2 and Regime 0.3 there exists an additional Gaussian part involving  the drift $F$ in the limit equations. In particular, if $G\equiv\mI_{d_2}$, the limit is an Ornstein-Uhlenbeck type  process. While in Regime 0.1 (i.e., $\beta_\eps/\alpha_\eps\to0$, which implies the term $\beta_\eps^{-1}c(X_t^\eps,Y_t^{0,\eps})$ is varying fast enough), there exists only a new drift part involving the  term $c$  in the limit equation for $\bar Z^0_{1,t}$. In particular, if $G\equiv \mI_{d_2}$, then $\bar Z^0_{1,t}$ will satisfy a linear random ODE.

\vspace{1mm}
(iii) Let us give some intuitive explanations for the above result. If $\beta_\eps/\alpha_\eps\to0$ (Regime 0.1), then $Y_t^{0,\eps}$ will converge to $\bar Y_t^1$ with rate $\alpha_\eps^2/\beta_\eps$. This means that the fast term $\beta_\eps^{-1}c(X_t^\eps,Y_t^{0,\eps})$ is the dominant  term in the strong convergence. Thus its  averaging effect appears in the functional central limit theorem. While if $\alpha_\eps/\beta_\eps\to0$ (Regime 0.2),  then $Y_t^{0,\eps}$ converges to $\bar Y_t^1$ with order $\alpha_\eps$ (independent of $\beta_\eps$). This suggests that  the term $\beta_\eps^{-1}c(X_t^\eps,Y_t^{0,\eps})$ is of lower order now, whose effect will go to zero eventually in the homogenization procedure. Finally, when $\alpha_\eps=\beta_\eps$ (Regime 0.3), there is a balance between the fluctuations involving $\beta_\eps^{-1}c(X_t^\eps,Y_t^{0,\eps})$ and $F(t,X_t^\eps,Y_t^{0,\eps})$, and thus the effects of both terms can be observed in the limit equation.
\er

\subsection{Functional central limit theorem: homogenization case}
In this subsection, we consider SDE (\ref{sde0}) with $H\neq0$, i.e., there exists a fast varying term even in the slow component. According to Theorem \ref{main1}, the averaged equation for $Y_t^\eps$ can be divided into two cases: Regime 1 and Regime 2. We proceed to identify the asymptotic limit for the normalized difference in each regimes: for $k=1,2$,
$$
Z_t^{k,\eps}:=\frac{Y_t^\eps-\bar Y^k_t}{\eta_\eps}
$$
with suitable deviation scale $\eta_\eps$ such that $\eta_\eps\to0$ as $\eps\to0$.

Let us first consider  Regime 1 in (\ref{regime}). Note that in this case, SDE (\ref{sde0}) has the same averaged equation as system (\ref{sde99}). According to Theorem \ref{main1} (i) with $\vartheta=1$, the process $Y_t^\eps$ will converge  to $\bar Y_t^1$ strongly with the best  possible rate $\alpha_\eps/\gamma_\eps+\alpha_\eps^2/(\beta_\eps\gamma_\eps)$. Thus by the same formal discussions as before, we expect that the natural choice of the deviation scale $\eta_\eps$ in order to observe non-trivial homogenization behavior for $Z^{1,\eps}_t$ should be divided into the following three regimes:
\begin{equation}\label{regime2}
\left\{\begin{aligned}
&\eta_\eps=\frac{\alpha_\eps^2}{\beta_\eps\gamma_\eps}\quad\,\text{and}\qquad\lim_{\eps\to0}\frac{\beta_\eps}{\alpha_\eps}=0,\qquad\qquad\qquad\,\,\, \text{Regime 1.1};\\ &\eta_\eps=\frac{\alpha_\eps}{\gamma_\eps}\qquad\text{and}\qquad\lim_{\eps\to0}\frac{\alpha_\eps}{\beta_\eps}=0,\qquad\qquad\qquad\quad\!\!\! \text{Regime 1.2};\\
&\eta_\eps=\frac{\alpha_\eps}{\gamma_\eps}\qquad\text{and}\qquad\alpha_\eps=\beta_\eps,\qquad\qquad\qquad\qquad\text{Regime 1.3}.
\end{aligned}\right.
\end{equation}
Recall that $\Phi$ is the unique solution to the Poisson equation (\ref{pde1}), and define
\begin{align}
\overline{c\cdot\nabla_x\Phi}(t,y)&:=\int_{\mR^{d_1}}c(x,y)\cdot\nabla_x\Phi(t,x,y)\mu^{y}(\dif x);\label{cphi}\\
\overline{H\cdot\Phi^*}(t,y)&:=\int_{\mR^{d_1}}H(t,x,y)\cdot\Phi^*(t,x,y)\mu^{y}(\dif x).\label{hphi}
\end{align}
Then the limit processes $\bar Z^1_{\ell,t}$ ($\ell=1,2,3$)  for $Z_t^{1,\eps}$ corresponding to Regime 1.1-Regime 1.3 in (\ref{regime2}) shall be given by:
\begin{align}\label{sdez1}
\begin{split}
\dif \bar Z_{1,t}^{1}&=\nabla_y\bar F_1(t,\bar Y_t^1)\bar Z_{1,t}^{1}\dif t+\nabla_yG(t,\bar Y_t^1)\bar Z_{1,t}^{1}\dif W_t^2+\overline{c\cdot\nabla_x\Phi}(t,\bar Y_t^1)\dif t;\\
\dif \bar Z_{2,t}^{1}&=\nabla_y\bar F_1(t,\bar Y_t^1)\bar Z_{2,t}^{1}\dif t+\nabla_yG(t,\bar Y_t^1)\bar Z_{2,t}^{1}\dif W_t^2+\sqrt{\overline{H\cdot\Phi^*}(t,\bar Y_t^1)}\dif \tilde W_t;\\
\dif \bar Z_{3,t}^{1}&=\nabla_y\bar F_1(t,\bar Y_t^1)\bar Z_{3,t}^{1}\dif t+\nabla_yG(t,\bar Y_t^1)\bar Z_{3,t}^{1}\dif W_t^2\\
&\quad+\overline{c\cdot\nabla_x\Phi}(t,\bar Y_t^1)\dif t+\sqrt{\overline{H\cdot\Phi^*}(t,\bar Y_t^1)}\dif \tilde W_t,
\end{split}
\end{align}
with initial data $\bar Z_{\ell,0}^{1}=0$ ($\ell=1,2,3$),  where  $\tilde W_t$ is another Brownian motion independent of $W_t^2$.

The following is our second main result for the functional central limit theorems.

\bt[Central limit theorem: Regime 1]\label{main2}
Let {\bf (A$_\sigma$)}, {\bf (A$_b$)}, {\bf (A$_{G}$)}, {\bf (A$_H$)} and (\ref{ac}) hold,  $0<\delta,\vartheta\leq 1$.  Then for any $T>0$ and every $\varphi\in C_b^{4}(\mR^{d_2})$,

\vspace{1mm}
\noindent(i) (Regime 1.1) if $b, \sigma\in C_b^{\delta,1+\vartheta}$, $F,H\in C_p^{(1+\vartheta)/2, \delta,1+\vartheta}$, $G\in C_b^{1/2, 1+\vartheta}$ and  $c\in C_p^{\delta,\vartheta}$, we have
\begin{align*}
\sup_{t\in[0,T]}\Big|\mE[\varphi(Z_t^{1,\eps})]-\mE[\varphi(\bar Z^1_{1,t})]\Big|\leq C_T\Big(\gamma_\eps+\frac{\beta_\eps^2}{\alpha_\eps^2}+\frac{\alpha_\eps^{2\vartheta}}{\beta_\eps^\vartheta\gamma_\eps^\vartheta}\Big);
\end{align*}

\noindent(ii) (Regime 1.2) if $b, \sigma\in C_b^{\delta,1+\vartheta}$, $F,H\in C_p^{(1+\vartheta)/2, \delta,1+\vartheta}$, $G\in C_b^{1/2, 1+\vartheta}$ and  $c\in L^\infty_p$, we have
\begin{align*}
\sup_{t\in[0,T]}\Big|\mE[\varphi(Z_t^{1,\eps})]-\mE[\varphi(\bar Z^1_{2,t})]\Big|\leq C_T\Big(\gamma_\eps+\frac{\alpha_\eps}{\beta_\eps}+\frac{\alpha_\eps^\vartheta}{\gamma_\eps^\vartheta}\Big);
\end{align*}
\noindent(iii) (Regime 1.3) if $b, \sigma\in C_b^{\delta,1+\vartheta}$, $F,H\in C_p^{(1+\vartheta)/2, \delta,1+\vartheta}$, $G\in C_b^{1/2, 1+\vartheta}$ and  $c\in C_p^{\delta,\vartheta}$, we have
\begin{align*}
\sup_{t\in[0,T]}\Big|\mE[\varphi(Z_t^{1,\eps})]-\mE[\varphi(\bar Z^1_{3,t})]\Big|\leq C_T\Big(\gamma_\eps+\frac{\alpha_\eps^\vartheta}{\gamma_\eps^\vartheta}\Big),
\end{align*}
where for $\ell=1,2,3$, $\bar Z^1_{\ell,t}$ satisfy the linear equation (\ref{sdez1}), and $C_T>0$ is a constant independent of $\delta, \eps$.
\et

\br\label{ma3}
(i) Compared with SDE (\ref{sdez0}) and Theorem \ref{main3},   the homogenization effect of the drift $F$ never appears in SDE (\ref{sdez1}). In fact, all the corresponding terms involving $F$ (through the Poisson equation (\ref{pde11})) in SDE (\ref{sdez0}) are now replaced by the drift  $H$ (through the Poisson equation (\ref{pde1})). This is intuitively natural   since in this case the fluctuations will be dominated by the fast component $\gamma_\eps^{-1}H(t,X_t^\eps,Y_t^{\eps})$, and the term $F(t,X_t^\eps,Y_t^\eps)$ is only of lower order now, thus its  effect in the fluctuations will go to zero eventually in each regime.

\vspace{1mm}
(ii) In particular, if $G\equiv \mI_{d_2}$ and $F\equiv0$ then Theorem \ref{main1} asserts that $Y_t^\eps$ converges strongly to $\bar Y^1_t=y+W_t^2$. Thus the deviation process is given by
$$
Z_t^{1,\eps}=z+\frac{1}{\eta_\eps\gamma_\eps}\int_0^tH(s,X_s^\eps,Y_s^\eps)\dif s,
$$
which is an inhomogeneous integral functional of the Markov process $(X_t^\eps,Y_t^\eps)$. Theorem \ref{main2} provides the limit for $Z_t^{1,\eps}$ in each regime, which is of independent interest (see e.g. \cite{CCG,Gu}).

\vspace{1mm}
(iii) Let us give more intuitive explanations for Regime 1.1 and Regime 1.2. Under Regime 1.1, we have $\alpha_\eps/\gamma_\eps$ goes to 0 faster than $\alpha_\eps^2/(\beta_\eps\gamma_\eps)$, and the process $Y_t^{\eps}$ will converge to $\bar Y_t^1$ with rate $\alpha_\eps^2/(\beta_\eps\gamma_\eps)$. Thus the fast term $\beta_\eps^{-1}c(X_t^\eps,Y_t^{\eps})$ is the leading term in the convergence and its averaging effect appears in the limit equation of $\bar Z^1_{1,t}$. While in Regime 1.2, we have $\alpha_\eps^2/(\beta_\eps\gamma_\eps)$ goes to 0 faster than $\alpha_\eps/\gamma_\eps$, and  the process $Y_t^{\eps}$ converges to $\bar Y_t^1$ with order $\alpha_\eps/\gamma_\eps$ (independent of $\beta_\eps$). This implies that  the term  $\beta_\eps^{-1}c(X_t^\eps,Y_t^{\eps})$ is of lower order now, whereas $\gamma_\eps^{-1}H(t,X_t^\eps,Y_t^\eps)$ is the leading term and its  homogenization effect appears in the limit equation of $\bar Z^1_{2,t}$.
\er

Now we consider   Regime 2 in (\ref{regime}), where homogenization  already occurs  even in the functional law of large numbers. Recall that $\alpha_\eps^2=\beta_\eps\gamma_\eps$ in this case. In particular, we shall always have $\lim_{\eps\to0}\alpha_\eps/\beta_\eps=\infty$. According to Theorem \ref{main1} (ii) with $\vartheta=1$, the process $Y_t^\eps$ will converge to $\bar Y_t^2$ strongly with the best  possible rate $\alpha_\eps/\gamma_\eps+\alpha_\eps^2/\beta_\eps$. Thus the natural choice of the derivation scale $\eta_\eps$ in order to observe non-trivial homogenization behavior for $Z^{2,\eps}_t$ should be divided into the following three regimes:
\begin{equation}\label{regime3}
\left\{\begin{aligned}
&\eta_\eps=\frac{\alpha_\eps^2}{\beta_\eps}\qquad\!\text{and}\qquad\lim_{\eps\to0}\frac{\beta_\eps}{\alpha_\eps\gamma_\eps}=0,\qquad\quad\,\,\,\quad\qquad \text{Regime 2.1};\\ &\eta_\eps=\frac{\alpha_\eps}{\gamma_\eps}\qquad\!\text{and}\qquad\lim_{\eps\to0}\frac{\alpha_\eps\gamma_\eps}{\beta_\eps}=0,\qquad\qquad\qquad\quad\!\!\! \text{Regime 2.2};\\
&\eta_\eps=\frac{\alpha_\eps}{\gamma_\eps}=\frac{\alpha_\eps^2}{\beta_\eps},\qquad\qquad\qquad\qquad\qquad\qquad\qquad\,\quad\text{Regime 2.3}.
\end{aligned}\right.
\end{equation}
Recall that $\Phi$ is the unique solution to the Poisson equation (\ref{pde1}), and $\overline{c\cdot\nabla_x\Phi}$ is defined by (\ref{cphi}).
Let $\Psi$ solves the following Poisson equation:
\begin{align}\label{p2}
\sL_0(x,y)\Psi(t,x,y)=-\big[c(x,y)\cdot\nabla_x\Phi(t,x,y)-\overline{c\cdot\nabla_x\Phi}(t,y)\big],
\end{align}
where $\sL_0(x,y)$ is defined by (\ref{l0}), and $(t,y)\in\mR_+\times\mR^{d_2}$ are regarded as parameters. Define
\begin{align*}
\overline{c\cdot\nabla_x\Psi}(t,y)&:=\int_{\mR^{d_1}}c(x,y)\cdot\nabla_x\Psi(t,x,y)\mu^{y}(\dif x).
\end{align*}
Then the limit processes $\bar Z^2_{\ell,t}$ ($\ell=1,2,3$)  for $Z_t^{2,\eps}$ corresponding to Regime 2.1-Regime 2.3 in (\ref{regime3}) shall be given by:
\begin{align}\label{sdez2}
\begin{split}
\dif \bar Z_{1,t}^{2}&=\nabla_y\bar F_2(t,\bar Y_t^2)\bar Z_{1,t}^{2}\dif t+\nabla_yG(t,\bar Y_t^2)\bar Z_{1,t}^{2}\dif W_t^2\\
&\quad+\big(\overline{c\cdot\nabla_x\Upsilon}+\overline{c\cdot\nabla_x\Psi}\big)(t,\bar Y_t^2)\dif t;\\
\dif \bar Z_{2,t}^{2}&=\nabla_y\bar F_2(t,\bar Y_t^2)\bar Z_{2,t}^{2}\dif t+\nabla_yG(t,\bar Y_t^2)\bar Z_{2,t}^{2}\dif W_t^2+\sqrt{\overline{H\cdot\Phi^*}(t,\bar Y_t^2)}\dif\tilde W_t;\\
\dif \bar Z_{3,t}^{2}&=\nabla_y\bar F_2(t,\bar Y_t^2)\bar Z_{3,t}^{2}\dif t+\nabla_yG(t,\bar Y_t^2)\bar Z_{3,t}^{2}\dif W_t^2\\
&\quad+\big(\overline{c\cdot\nabla_x\Upsilon}+\overline{c\cdot\nabla_x\Psi}\big)(t,\bar Y_t^2)\dif t+\sqrt{\overline{H\cdot\Phi^*}(t,\bar Y_t^2)}\dif \tilde W_t
\end{split}
\end{align}
with initial data $\bar Z_{\ell,0}^{2}=0$, where $\overline{c\cdot\nabla_x\Upsilon}$ and $\overline{H\cdot\Phi^*}$ are defined by (\ref{cup}) and (\ref{hphi}), respectively, $\bar Y_t^2$ is the unique strong solution of SDE (\ref{sde000}) with $k=2$,  and $\tilde W_t$ is another Brownian motion independent of $W_t^2$.

Our  main result in this case is as follows.

\bt[Central limit theorem: Regime 2]\label{main4}
Let {\bf (A$_\sigma$)}, {\bf (A$_b$)}, {\bf (A$_{G}$)}, {\bf (A$_H$)} and (\ref{ac}) hold,   $0<\delta,\vartheta\leq 1$. Assume that $b, \sigma\in C_b^{\delta,1+\vartheta}$, $c\in C_p^{\delta,1+\vartheta}$, $F,H\in C_p^{(1+\vartheta)/2, \delta,1+\vartheta}$ and $G\in C_b^{1/2, 1+\vartheta}$. Then for any $T>0$ and every $\varphi\in C_b^{4}(\mR^{d_2})$,

\vspace{1mm}
\noindent(i) (Regime 2.1) we have
\begin{align*}
\sup_{t\in[0,T]}\Big|\mE[\varphi(Z_t^{2,\eps})]-\mE[\varphi(\bar Z^2_{1,t})]\Big|\leq C_T\Big(\frac{\beta_\eps^2}{\alpha_\eps^2\gamma_\eps^2}+\frac{\alpha_\eps^{2\vartheta}}{\beta_\eps^\vartheta}\Big);
\end{align*}

\noindent(ii) (Regime 2.2) we have
\begin{align*}
\sup_{t\in[0,T]}\Big|\mE[\varphi(Z_t^{2,\eps})]-\mE[\varphi(\bar Z^2_{2,t})]\Big|\leq C_T\Big(\frac{\alpha_\eps\gamma_\eps}{\beta_\eps}+\frac{\alpha_\eps^\vartheta}{\gamma_\eps^\vartheta}\Big);
\end{align*}
\noindent(iii) (Regime 2.3) we have
\begin{align*}
\sup_{t\in[0,T]}\Big|\mE[\varphi(Z_t^{2,\eps})]-\mE[\varphi(\bar Z^2_{3,t})]\Big|\leq C_T \frac{\alpha_\eps^\vartheta}{\gamma_\eps^\vartheta},
\end{align*}
where for $\ell=1,2,3$, $\bar Z^2_{\ell,t}$ satisfy the linear equation (\ref{sdez2}), and $C_T>0$ is a constant independent of $\delta, \eps$.
\et

\br\label{m4}
(i) Since homogenization  already occurs in the averaged equation of $\bar Y_t^2$, it is natural to expect that the second order homogenization will appear in the functional central limit theorem. This is exactly characterized through the function $\Psi$ in Regime 2.1 and Regime 2.3 which solves the Poisson equation (\ref{p2}). In  Regime 2.2, the slow process $Y_t^\eps$ will converge  to  $\bar Y_t^2$ with rate $\alpha_\eps/\gamma_\eps$ (independent of $\beta_\eps$), which means that  the term $\gamma_\eps^{-1}H(t,X_t^\eps,Y_t^\eps)$ is the only leading term, and thus the averaging effect involving $\beta_\eps^{-1}c(X_t^\eps,Y_t^\eps)$ does not arise.

\vspace{1mm}
(ii)  Note that the same homogenization behaviors involving the drift $F$ as in  SDE (\ref{sdez0}) appear again in Regime 2.1 and Regime 2.3. This implies that the component $\gamma_\eps^{-1}H(t,X_t^\eps,Y_t^\eps)$ is of the same order as the drift term $F(t,X_t^\eps,Y_t^\eps)$ in the  fluctuations, which should not be a contradiction since in Regime 2, the fast varying of $\gamma_\eps^{-1}H(t,X_t^\eps,Y_t^\eps)$ has already been homogenized out in the averaged equation.

\vspace{1mm}
(iii)   It is interesting to note that in order to observe all non-trivial behaviors of every term  simultaneously, we need to take $\gamma_\eps=\alpha_\eps^{1/2}$ and $\beta_\eps=\alpha_\eps^{3/2}$ (Regime 2.3) to balance the averaging effect of $\beta_\eps^{-1}c(X_t^\eps,Y_t^\eps)$ and the homogenization effect of $\gamma_\eps^{-1}H(t,X_t^\eps,Y_t^\eps)$.
\er

Finally, we mention that regime specific
analysis for system (\ref{sde0}) is also done in \cite{RX,S1}. We list the following comparisons with our results.

\br\label{ma5}
(i) In \cite{S1}, the convergence in probability and the central limit theorem for system (\ref{sde0}) was studied in a  homogeneous case where $\alpha_\eps=\delta/\sqrt{\eps}$, $\beta_\eps=\delta$,  $\gamma_\eps=\delta/\eps$  and with small noise perturbations, i.e., with $G$ replaced by $\sqrt{\eps}G$. Thus, the corresponding results in \cite{S1} can be seen as a particular case of Regime 2 and Regimes 2.1-2.3 of this paper. Moreover, we handle non-smooth coefficients, which is much more general than the case studied in \cite{S1}, and we obtain optimal rates of convergence. Note that all the above rates of convergence do not depend
on the regularity of the coefficients with respect to the fast variable. This reflects that the slow process is the main term in the limiting procedure of
the multi-scale system, which coincides with the intuition since the fast component has been
totally averaged or homogenized out in the limit equation.

The author in \cite{S1} also considered the cases where $\gamma_\eps\to\gamma\in(0,\infty)$ as $\eps\to0$. In this case, the part $\gamma_\eps^{-1}H(t,X_t^\eps,Y_t^\eps)$ is no longer a fast term and can be handled by using the same arguments as for the drift part $F$, which is sightly easier.

\vspace{1mm}
(ii)  In \cite{RX}, only convergence in distribution for system (\ref{sde0}) was studied. In Theorem \ref{main1}, we provide the strong convergence in the pathwise sense. Due to the low regularity conditions on the coefficients, Zvonkin's transformation (see Lemma \ref{zt} below) is needed to handle the non-smooth coefficients. Furthermore, the central limit results in this paper provide more delicate  characterisation   for the interactions between the  averaging term $\beta_\eps^{-1}c(X_t^\eps,Y_t^{\eps})$ and the homogenization term $\gamma_\eps^{-1}H(t,X_t^\eps,Y_t^{\eps})$ in system (\ref{sde0}).
\er

\vspace{1mm}
\noindent{\bf Notations:}
Since we shall prove the main results in a quite unified way, we introduce some notations here for brevity.
Let $\bar\sL_k$ be the infinitesimal generator for $\bar Y_t^k$, i.e., for $k=1,2$,
\begin{align}\label{LL1}
\bar\sL_k&:=\bar\sL_k(t,y):=\sum_{i,j=1}^{d_2} \cG^{ij}(t,y)\frac{\p^2}{\p y_i\p y_j}+\sum_{i=1}^{d_2}\bar F^i_k(t,y)\frac{\p}{\p y_i},
\end{align}
where $\bar F_k$ are defined by (\ref{barf}).
Note that the averaged process for $Y_t^{0,\eps}$ in SDE (\ref{sde99}) is also given by $\bar Y_t^1$, we let $\bar Y_t^0:=\bar Y_t^1$ and $\bar \sL_0:=\bar \sL_1$ for consistency.
We also introduce
\begin{align}
\bar \sL_3^0&:=\bar \sL_3^0(t,y,z):=\sum_{i=1}^{d_2}\big(\overline{c\cdot\nabla_x\up}\big)^i(t,y)\frac{\p}{\p z_i},\label{LL03}\\
\bar\sL^1_3&:=\bar\sL^1_3(t,y,z):=\sum_{i=1}^{d_2}\big(\overline{c\cdot\nabla_x\Phi}\big)^i(t,y)\frac{\p}{\p z_i},\label{LL13}\\
\bar\sL^2_3&:=\bar\sL^2_3(t,y,z):=\sum_{i=1}^{d_2}\big(\overline{c\cdot\nabla_x\Psi}\big)^i(t,y)\frac{\p}{\p z_i},\label{LL23}\\
\bar\sL^1_4&:=\bar\sL^1_4(t,y,z):=\sum_{i,j=1}^{d_2}\big(\overline{\tilde F\cdot\up^*}\big)^{ij}(t,y)\frac{\p^2}{\p z_i\p z_j}, \label{LL14}\\
\bar\sL^2_4&:=\bar\sL^2_4(t,y,z):=\sum_{i,j=1}^{d_2}\big(\overline{H\cdot\Phi^*}\big)^{ij}(t,y)\frac{\p^2}{\p z_i\p z_j},\label{LL24}
\end{align}
and for $k=1,2$,
\begin{align}\label{LL05}
\bar \sL_5^k:=\bar \sL_5^k(t,y,z)&:=\sum_{i=1}^{d_2}\big(\nabla_y\bar F_k(t,y) z\big)^{i}\frac{\p}{\p z_i}+\frac{1}{2}\sum_{i,j=1}^{d_2}\big(G(t,y)[\nabla_yG(t,y) z]^*\big)^{ij}\frac{\p^2}{\p y_i\p z_j}\no\\
&\quad+\frac{1}{2}\sum_{i,j=1}^{d_2}\big([\nabla_yG(t,y) z][\nabla_yG(t,y) z]^*\big)^{ij}\frac{\p^2}{\p z_i\p z_j}.
\end{align}
Then the  infinitesimal generator of $(\bar Y_t^k, \bar Z_{\ell,t}^k)$ can be written as $\bar\sL_k+\bar\cL_\ell^k$, where for $k=0,1,2$ and $\ell=1,2,3$, the operator $\bar\cL_\ell^k$ are defined as follows:\\
corresponding to $\bar Z^0_{\ell,t}$ ($\ell=1,2,3$) in SDE (\ref{sdez0}),
\begin{align}\label{rho0}
\bar\cL_1^0:=\bar\sL_5^1+\bar\sL_3^0,\quad\bar\cL_2^0:=\bar\sL_5^1+\bar\sL_4^1,\quad\bar\cL_3^0:=\bar\sL_5^1+\bar\sL_3^0+\bar\sL_4^1;
\end{align}
corresponding to  $\bar Z^1_{\ell,t}$ ($\ell=1,2,3$) in SDE (\ref{sdez1}),
\begin{align}\label{rho1}
\bar\cL_1^1:=\bar\sL_5^1+\bar\sL_3^1,\quad\bar\cL_2^1:=\bar\sL_5^1+\bar\sL_4^2,\quad\bar\cL_3^1:=\bar\sL_5^1+\bar\sL_3^1+\bar\sL_4^2;
\end{align}
and corresponding to   $\bar Z^2_{\ell,t}$ ($\ell=1,2,3$) in SDE (\ref{sdez2}),
\begin{align}\label{rho2}
\bar\cL_1^2:=\bar\sL_5^2+\bar\sL_3^0+\bar\sL_3^2,\quad\bar\cL_2^2:=\bar\sL_5^2+\bar\sL_4^2,\quad\bar\cL_3^2:=\bar\sL_5^2+\bar\sL_3^0+\bar\sL_3^2+\bar\sL_4^2.
\end{align}

\section{Poisson equation and Cauchy problem}

This section collects the main tools that we shall use to prove our results.
As  mentioned in the introduction, the Poisson equation will play an important role in the proof both for the strong convergence in the averaging principle and the central limit theorems.  Let us first recall some results in this direction.

Consider the following Poisson equation in $\mR^{d_1}$:
\begin{align}
\sL_0(x,y)u(x,y)=-f(x,y),  \label{pde0}
\end{align}
where $\sL_0(x,y)$ is defined by (\ref{l0}), and $y\in\mR^{d_2}$ is regarded as a parameter. Note that $\sL_0(x,y)$ is the infinitesimal generator of $X_t^y$ given by (\ref{sde2}). To ensure the well-posedness of equation (\ref{pde0}), we need to make the following centering condition on $f$:
\begin{align}\label{cen2}
\int_{\mR^{d_1}}f(x,y)\mu^y(\dif x)=0,\quad\forall y\in\mR^{d_2},
\end{align}
where $\mu^y(\dif x)$ is the invariant measure  for $X_t^y$.
The following result was proved in \cite[Theorem 2.1]{RX}, which will be used frequently below.

\bt\label{popde}
Let {\bf (A$_\sigma$)} and {\bf (A$_b$)} hold. Assume that $a, b\in C_b^{\delta,\vartheta}$ with $0<\delta\leq 1$ and $\vartheta\geq0$. Then for every function $f\in C_p^{\delta,\vartheta}$ satisfying (\ref{cen2}), there exists a unique solution $u\in C_p^{2+\delta,\vartheta}$ to equation (\ref{pde0}) satisfying (\ref{cen2}) which is given by
$$
u(x,y)=\int_0^\infty \mE f(X_t^y(x),y)\dif t.
$$
Moreover, there exists a constant $m>0$ such that:

\vspace{1mm}
\noindent
(i) for any $x\in\mR^{d_1}$ and $y\in\mR^{d_2}$,
\begin{align}\label{key1}
|u(x,y)|+|\nabla_xu(x,y)|+|\nabla_x^2u(x,y)|\leq C_0(1+|x|^m),
\end{align}
\indent\,\,\! where  $C_0>0$ depends only on $d_1,d_2$ and $\|a\|_{C_b^{\delta,0}}, \|b\|_{C_b^{\delta,0}}, [f]_{C_p^{\delta,0}}$;

\noindent(ii) for any $x\in\mR^{d_1}$,
\begin{align}\label{key3}
\|u(x,\cdot)\|_{C_b^\vartheta}\leq C_1(1+|x|^m),
\end{align}
\indent\,\,\! where  $C_1>0$ depends on $d_1,d_2$ and $\|a\|_{C_b^{\delta,\vartheta}}, \|b\|_{C_b^{\delta,\vartheta}}, [f]_{C_p^{\delta,\vartheta}}$.
\et

We will need to use  It\^o's formula for the solution of the Poisson equation with the fast and slow components in SDE (\ref{sde0}) plugged in for both variables, say $u(X_t^\eps,Y_t^\eps)$, which in turn requires at least two  derivatives
of $u$ with respect to the $x$ variable as well as the  $y$ variable.  In view of estimate (\ref{key1}), the derivatives with respect to  $x$  are not a problem since we can get them for free by virtue of the uniform  ellipticity property of the operator. However,
due to our low regularities of the coefficients with respect to the $y$ variable  (only H\"older continuous) and taking into account  (\ref{key3}), we cannot get the desired two derivatives for $u(x,\cdot)$ directly.
To overcome this problem, we use some mollification arguments.

Let $\rho_1:\mR\to[0,1]$ and $\rho_2:\mR^{d_2}\to[0,1]$ be two smooth radial convolution kernel functions
such that
$\int_\mR\rho_1(r)\dif r=\int_{\mR^{d_2}}\rho_2(y)\dif y=1$, and for any $k\geq 1$, there exist constants $C_k>0$ such that $|\nabla^k\rho_1(r)|\leq C_k\rho_1(r)$ and  $|\nabla^k\rho_2(y)|\leq C_k\rho_2(y)$. For every $n\in\mN^*$, set
$$
\rho_1^n(r):=n^2\rho_1(n^2r)\quad \text{and}\quad
\rho_2^n(y):=n^{d_2}\rho_2(ny).
$$
Given a function $f(t,x,y,z)$, we define the mollifying approximations of $f$ in $t$ and $y$ variables by
\begin{align}\label{fn}
f_n(t,x,y,z):=f*\rho_2^n*\rho_1^n:=\int_{\mR^{d_2+1}}f(t-s,x,y-y',z)\rho_2^{n}(y')\rho_1^n(s)\dif y'\dif s.
\end{align}
The following easy result can be proved similarly as in \cite[Lemma 4.1]{RX}, we omit the details.

\bl
Let $f\in C_p^{\vartheta/2,0,\vartheta,0}$ with $0<\vartheta\leq 2$ and define $f_n$ by (\ref{fn}). Then we  have
\begin{align}
\|f(\cdot,x,\cdot,\cdot)-f_n(\cdot,x,\cdot,\cdot)\|_\infty&\leq C_0n^{-\vartheta}(1+|x|^m),\label{n111}\\
\|\nabla_yf_n(\cdot,x,\cdot,\cdot)\|_{\infty}&\leq C_0n^{1-(\vartheta\wedge1)}(1+|x|^m),\label{n333}
\end{align}
and
\begin{align}\label{n222}
\|\p_tf_n(\cdot,x,\cdot,\cdot)\|_{\infty}+\|\nabla^2_yf_n(\cdot,x,\cdot,\cdot)\|_{\infty}\leq C_0n^{2-\vartheta}(1+|x|^m),
\end{align}
where $C_0>0$ is a constant independent of $n$.
\el

Given a function $h(t,x,y)$, we  denote its average with respect to the measure $\mu^y(\dif x)$ by $\bar h(t,y)$, i.e.,
\begin{align}\label{avh}
\bar h(t,y):=\int_{\mR^{d_1}}h(t,x,y)\mu^y(\dif x).
\end{align}
The following result specifies the regularity of an averaged function, which explains the assumptions we made on the coefficients in our main results.

\bl\label{coe}
Let {\bf (A$_\sigma$)} and {\bf (A$_b$)} hold. Assume that $a, b\in C_b^{\delta,\vartheta}$ with $0<\delta\leq 1$ and $\vartheta\geq0$. Then for every $h\in C_p^{\vartheta/2,\delta,\vartheta}$, we have $\bar h\in C_b^{\vartheta/2,\vartheta}$. In particular,

(i) under conditions in Theorem \ref{main1},  we have for every $k=1,2$,
$\bar F_k\in C_b^{\vartheta/2,\vartheta}$;

(ii) under conditions in Theorem \ref{main3}, Theorem \ref{main2} and Theorem \ref{main4}, we have for every $k=1,2$,
$$
\nabla_y\bar F_k, \overline{c\cdot\nabla_x\up}, \overline{\tilde F\cdot\up^*}, \overline{c\cdot\nabla_x\Phi}, \overline{H\cdot\Phi^*}, \overline{c\cdot\nabla_x\Psi}\in C_b^{\vartheta/2,\vartheta}.
$$
\el

\begin{proof}
The assertion that $\bar h\in C_b^{\vartheta/2,\vartheta}$ was proved in \cite[Lemma 3.2]{RX}. Then, under  the assumptions in Theorem \ref{main1} (Regime 1), the conclusion that $\bar F_1\in C_b^{\vartheta/2,\vartheta}$ follows directly. Recall that $\Phi$ solves (\ref{pde1}). By the assumptions in Theorem \ref{main1} (Regime 2) and Theorem \ref{popde}, we have $\Phi\in C_p^{\vartheta/2,2+\delta,\vartheta}$. This together with the condition that $c\in C_p^{\delta,\vartheta}$ implies that $c\cdot\nabla_x\Phi\in C_p^{\vartheta/2,2+\delta,\vartheta}$, which in turn yields $\bar F_2\in C_b^{\vartheta/2,\vartheta}$.  (ii) can be proved by the same arguments, so we omit the details.
\end{proof}

Another main tool we will use to prove the functional central limit theorems  is the Cauchy problem corresponding to the limit dynamics $(\bar Y_t^k, \bar Z^k_{\ell,t})$, $k=0,1,2$ and $\ell=1,2,3$. Note that the processes $\bar Y_t^k$ depend on the initial value $y$, while $\bar Z^k_{\ell,t}$ depend  on $y$ but with initial value $0$. Below, we shall write $\bar Y_t^k(y)$  when we want to stress the dependence on the
initial value, and use $\bar Z^k_{\ell,t}(y,z)$ to denote processes $\bar Z^k_{\ell,t}$ with initial point $z\in\mR^{d_2}$.
Fix a $T>0$ below, consider the following Cauchy problem on $[0,T]\times\mR^{d_2}\times\mR^{d_2}$:
\begin{equation}\left\{\begin{array}{l}\label{PDE}
\displaystyle
\p_t  u^k_\ell(t,y,z)-(\bar\sL_k+\bar\cL^k_\ell) u^k_\ell(t,y,z)=0,\quad t\in (0, T],\\
u^k_\ell(0,y,z)=\varphi(z),
\end{array}\right.
\end{equation}
where  $\bar \sL_k$ and $\bar\cL_\ell^k$ are defined by (\ref{LL1}), (\ref{rho0}), (\ref{rho1}) and (\ref{rho2}), respectively.
We have the following result.

\bt\label{pde}
For every  $k=0,1,2$, $\ell=1,2,3$ and $\varphi\in C_b^4(\mR^{d_2})$, there exists a unique solution $u^k_\ell\in C_b^{(2+\vartheta)/2,2+\vartheta,4}$ to equation (\ref{PDE}) which is given by
\begin{align}\label{ukl}
u^k_\ell(t,y,z)=\mE\varphi\big(\bar Z^k_{\ell,t}(y,z)\big).
\end{align}
\et

\begin{proof}
We only prove the assertion for $k=0$ and $\ell=1$.	Although this is not the most general one, we  choose this case since it carries the key difficulties. For simplicity, we shall write $u$ instead of $u^0_1$ for the solution. Without loss of generality, we may assume that the coefficients are smooth, and  focus on proving the a priori estimates for $u$. Since $\bar\sL_0+\bar\cL^0_1$ is the generator of the Markov process $(\bar Y_t^0, \bar Z^0_{1,t})$,
it is well known that the solution for (\ref{PDE}) will be given  by (\ref{ukl}). By the assumption $\varphi\in C_b^4(\mR^{d_2})$, the fact that $\bar Y_t^0$ does not depend on $z$, and since $\bar Z^0_{1,t}(y,z)$ satisfies the linear equation (\ref{sdez0}), it is easily checked that for every $(t,y)\in \mR_+\times\mR^{d_2}$,
we have $u(t,y,\cdot)\in  C_b^4(\mR^{d_2})$, and for $i=1,\cdots,4$,
$$
\|\nabla_z^iu(t,y,\cdot)\|_\infty\leq C_0\|\varphi\|_{C_b^4},
$$
where $C_0>0$ depends only on $\|\nabla_y\bar F_1\|_\infty$ and $\|\nabla_yG\|_\infty$.
It remains to prove that for every $z\in\mR^{d_2}$, $u(\cdot,\cdot,z) \in C_b^{(2+\vartheta)/2,2+\vartheta}$. To this end, we rewrite equation (\ref{PDE}) as
$$
\p_t  u(t,y,z)-\bar\sL_0 u(t,y,z)=\bar\cL^0_1u(t,y,z).
$$
By regarding $z$ as a parameter in the above equation, and recalling that $\cG$ is uniformly elliptic, it suffices to show that
$$
\nabla_zu(\cdot,\cdot,z), \nabla^2_zu(\cdot,\cdot,z)\in C_b^{ \vartheta/2,\vartheta}.
$$
Then the conclusion follows by classical PDE's result, see e.g. \cite[Chapter IV, Section 5]{La-So-Ur}.
For any $y_1,y_2\in\mR^{d_2}$, we  write
\begin{align*}
\big|\nabla_zu(t,y_1,z)&-\nabla_zu(t,y_2,z)\big|\leq \|\nabla_z\varphi\|_\infty\mE\big|\nabla_z\bar Z^0_{1,t}(y_1,z)-\nabla_z\bar Z^0_{1,t}(y_2,z)\big|\\
&\qquad\qquad+\mE\Big(\big|\nabla_z\varphi\big(\bar Z^0_{1,t}(y_1,z)\big)-\nabla_z\varphi\big(\bar Z^0_{1,t}(y_2,z)\big)\big|\cdot|\nabla_z\bar Z^0_{1,t}(y_2,z)|\Big)\\
&\leq C_1\mE\Big(\big|\nabla_z\bar Z^0_{1,t}(y_1,z)-\nabla_z\bar Z^0_{1,t}(y_2,z)\big|^2+\big|\bar Z^0_{1,t}(y_1,z)-\bar Z^0_{1,t}(y_2,z)\big|^2\Big)^{1/2}.
\end{align*}
Note that
\begin{align*}
\dif \nabla_z\bar Z^0_{1,t}=\nabla_y\bar F_1(t,\bar Y_t^1)\nabla_z\bar Z_{1,t}^{0}\dif t+\nabla_yG(t,\bar Y_t^1)\nabla_z\bar Z_{1,t}^{0}\dif W_t^2.
\end{align*}
Thus by the fact that $\nabla_y\bar F_k(t,\cdot), \nabla_yG(t,\cdot)\in C_b^{\vartheta}$, we deduce that
\begin{align}
\mE\big|\nabla_z\bar Z^0_{1,t}(y_1,z)&-\nabla_z\bar Z^0_{1,t}(y_2,z)\big|^2\leq C_2\mE\left(\int_0^t\big|\nabla_z\bar Z^0_{1,s}(y_1,z)-\nabla_z\bar Z^0_{1,s}(y_2,z)\big|^2\dif s\right)\no\\
&+C_2\mE\bigg(\int_0^t\big|\nabla_y\bar F_1(s,\bar Y_s^1(y_1))-\nabla_y\bar F_1(s,\bar Y_s^1(y_2))\big|^2\no\\
&\qquad\quad\qquad\quad+\big|\nabla_y\bar G(s,\bar Y_s^1(y_1))-\nabla_y\bar G(s,\bar Y_s^1(y_2))\big|^2\dif s\bigg)\no\\
&\leq C_2\mE\left(\int_0^t\big|\nabla_z\bar Z^0_{1,s}(y_1,z)-\nabla_z\bar Z^0_{1,s}(y_2,z)\big|^2\dif s\right)\no\\
&\quad+C_2\mE\bigg(\int_0^t\big|\bar Y_s^1(y_1)-\bar Y_s^1(y_2)\big|^{2\vartheta}\dif s\bigg).\label{y1y2}
\end{align}
It is well known that $y\mapsto \bar Y_t^1(y)$ is a $C^1$-diffeomorphism, i.e., for every $t\in[0,T]$,
$$
\mE\big|\bar Y_t^1(y_1)-\bar Y_t^1(y_2)\big|\leq C_T|y_1-y_2|.
$$
Taking this back into (\ref{y1y2}) and by Gronwall's inequality, we get
\begin{align*}
\mE\big|\nabla_z\bar Z^0_{1,t}(y_1,z)-\nabla_z\bar Z^0_{1,t}(y_2,z)\big|^2\leq C_3|y_1-y_2|^{2\vartheta}.
\end{align*}
By the same arguments and more easily, we also have
\begin{align*}
\mE\big|\bar Z^0_{1,t}(y_1,z)- \bar Z^0_{1,t}(y_2,z)\big|^2\leq C_4|y_1-y_2|^{2\vartheta},
\end{align*}
which in turn implies that
$$
\big|\nabla_zu(t,y_1,z)-\nabla_zu(t,y_2,z)\big|\leq C_5|y_1-y_2|^{\vartheta}.
$$
The corresponding regularity for $\nabla_zu$ with respect to $t$ variable and for $\nabla_z^2u$ can be proved similarly.
\end{proof}

\section{Strong convergence in the averaging principle}

Using the technique of Poisson equation, we shall first derive some fluctuation estimates in Subsection 4.1. Then we prove the strong convergence in the averaging principle of SDE (\ref{sde0}) in Subsection 4.2 by Zvonkin's transformation.

\subsection{Strong fluctuation estimates}

Given a function $h(t,x,y)$, recall that $\bar h(t,y)$ is defined by (\ref{avh}).
It is easy to see that $f(t,x,y):=h(t,x,y)-\bar h(t,y)$ satisfies the centering condition, i.e.,
\begin{align}\label{cen6}
\int_{\mR^{d_1}}f(t,x,y)\mu^y(\dif x)=0,\quad \forall (t,y)\in\mR_+\times\mR^{d_2}.
\end{align}
The following result gives an estimate for  the    fluctuations between $h(s,X_s^\eps,Y_s^\eps)$ and $\bar h(s, Y_s^\eps)$ over the time interval $[0,t]$.

\bl\label{key}
Let {\bf (A$_\sigma$)}, {\bf (A$_b$)} and (\ref{ac}) hold. Assume that $b,\sigma\in C_b^{\delta,\vartheta}$  with $0<\delta,\vartheta\leq 2$ and $c, F, H, G\in L^\infty_p$. Then for every $f\in C_p^{\vartheta/2,\delta,\vartheta}$  satisfying (\ref{cen6}) and any $q\geq 2$, we have
\begin{align}\label{fle1}
\mE\left|\int_0^tf(s,X_s^\eps,Y_s^\eps)\dif s\right|^q\leq C_t\Big(\alpha_\eps^{\vartheta\wedge1}+\frac{\alpha_\eps^2}{\beta_\eps}\Big)^q,
\end{align}
where $C_t>0$ is a constant independent of $\delta,\eps$.
\el

\br
We call (\ref{fle1})  a strong fluctuation estimate because we take the absolute value for the integral over $[0,t]$. Compared with Lemma \ref{key30} and Lemma \ref{key31} below, we shall see that the involved  martingale part will be one of the leading terms in the control of  error bounds in this case, and this is the main reason why the power $\vartheta\wedge1$ appears on the right hand side of (\ref{fle1}).
\er
\begin{proof}
By the assumptions that $f\in C_p^{\vartheta/2,\delta,\vartheta}$ satisfying  (\ref{cen6}), $b,\sigma\in C_b^{\delta,\vartheta}$ and  according to Theorem \ref{popde}, there exists a unique solution $\Phi^f(t,x,y)\in C^{\vartheta/2,2+\delta,\vartheta}_p$ to the following Poisson equation in $\mR^{d_1}$:
	\begin{align}\label{pn}
	\sL_0(x,y)\Phi^f(t,x,y)=-f(t,x,y),
	\end{align}
	where $(t,y)\in\mR_+\times\mR^{d_2}$ are regarded as parameters.
	Let $\Phi_n^f$ be the mollifyer of $\Phi^f$ defined as in (\ref{fn}) (which does not depend on the $z$-variable here).
Using It\^o's formula, we  have
	\begin{align*} \Phi_n^f(t,X_{t}^\eps,Y_{t}^\eps)&=\Phi_n^f(0,x,y)+\int_0^{t}\Big(\p_s+\sL_1+\gamma_\eps^{-1}\sL_2+\beta_\eps^{-1}\sL_3\Big)\Phi_n^f(s,X_s^\eps,Y_s^\eps)\dif s\\
	 &\quad+\frac{1}{{\alpha_\eps^2}}\int_0^{t}\sL_0\Phi_n^f(s,X_s^\eps,Y_s^\eps)\dif s+\frac{1}{\alpha_\eps}M^1_{n}(t)+M^2_{n}(t),
	\end{align*}
	where $\sL_1, \sL_2$ and $\sL_3$ are given  by (\ref{lll}), and for $i=1,2$, $M^i_{n}(t)$ are martingales defined by
	\begin{align*}
	 M^1_{n}(t)&:=\int_0^t\nabla_x\Phi_n^f(s,X_s^\eps,Y_s^\eps)\sigma(X_s^\eps, Y_s^\eps)\dif W_s^1,\\
	 M^2_{n}(t)&:=\int_0^t\nabla_y\Phi_n^f(s,X_s^\eps,Y_s^\eps)G(s,Y_s^\eps)\dif W_s^2.
	\end{align*}
	By (\ref{pn}) this in turn yields that
	\begin{align}\label{ito}
	\begin{split}
	\int^{t}_0\!f(s,X^{\eps}_s, Y^{\eps}_s)\dif s
	 &={\alpha_\eps^2}\Phi_n^f(0,x,y)-{\alpha_\eps^2}\Phi_n^f({t},X_{t}^\eps,Y_{t}^\eps)+{\alpha_\eps} M^1_{n}(t)+{\alpha_\eps^2}M^2_{n}(t)\\
	 &+{\alpha_\eps^2}\int_0^{t}\big(\p_s+\sL_1\big)\Phi_n^f(s,X_s^\eps,Y_s^\eps)\dif s\\
	 &+\frac{\alpha_\eps^2}{\gamma_\eps}\int_0^{t}\sL_2\Phi_n^f(s,X_s^\eps,Y_s^\eps)\dif s+\frac{\alpha_\eps^2}{\beta_\eps}\int_0^{t}\sL_3\Phi_n^f(s,X_s^\eps,Y_s^\eps)\dif s\\
	 &+\int_0^{t}\big(\sL_0\Phi^f_n-\sL_0\Phi^f\big)(s,X_s^\eps,Y_s^\eps)\dif s.
	\end{split}
	\end{align}
	As a result, we have for any $q\geq 2$,
	\begin{align*}
	\cQ(\eps)&:=\mE\left|\int_0^tf(s,X_s^\eps,Y_s^\eps)\dif s\right|^q\leq C_q\Big[{\alpha_\eps^{2q}}\Big(|\Phi_n^f(0,x,y)|^q+\mE|\Phi^f_n(t,X_t^\eps,Y_t^\eps)|^q\Big)\\
	&\quad+\alpha_\eps^{q}\,\mE|M_n^1(t)|^q+\alpha_\eps^{2q}\,\mE|M_n^2(t)|^q\Big] +C_q{\alpha_\eps^{2q}}\,\mE\left|\int_0^{t}\big(\p_s+\sL_1\big)\Phi_n^f(s,X_s^\eps,Y_s^\eps)\dif s\right|^q\\
	 &\quad+C_q\frac{\alpha_\eps^{2q}}{\gamma_\eps^q}\mE\left|\int_0^{t}\sL_2\Phi_n^f(s,X_s^\eps,Y_s^\eps)\dif s\right|^q+C_q\frac{\alpha_\eps^{2q}}{\beta_\eps^q}\mE\left|\int_0^{t}\sL_3\Phi_n^f(s,X_s^\eps,Y_s^\eps)\dif s\right|^q\\
	 &\quad+C_q\,\mE\left|\int_0^{t}\big(\sL_0\Phi_n^f-\sL_0\Phi^f\big)(s,X_s^\eps,Y_s^\eps)\dif s\right|^q=:\sum_{i=1}^5\cQ_{i}(\eps).
	\end{align*}
Under {\bf (A$_\sigma$)} and (\ref{ac}), it follows from \cite[Lemma 1]{Ve5} (see also \cite[Lemma 1]{P-V2} or \cite{RX}) that for any $m>0$, there exists a constant $C_0>0$ such that
	\begin{align}
\sup_{\eps\in(0,1/2)}\mE|X^{\eps}_t|^{m}\leq C_0(1+|x|^{m+2}).\label{ME}
	\end{align}
Thus by (\ref{key1}) we have  that there exist constants $m>0$ and $C_1>0$ independent of $n$ such that
	$$
	|\Phi_n^f(0,x,y)|^q+\mE|\Phi^f_n(t,X_t^\eps,Y_t^\eps)|^q\leq C_1\mE\big(1+|X_t^\eps|^{mq}\big)<\infty.
	$$	
At the same time, by H\"older's inequality,
	$$
	\mE|M_n^1(t)|^q\leq C_1\mE\left(\int_0^t\Big(1+|X_s^\eps|^{2m}\Big)\dif s\right)^{q/2}\leq C_1\mE\left(\int_0^t\Big(1+|X_s^\eps|^{mq}\Big)\dif s\right)<\infty,
	$$	
and in view of (\ref{n333}),
	$$
\mE|M_n^2(t)|^q\leq C_1n^{q(1-(\vartheta\wedge1))}\mE\left(\int_0^t\Big(1+|X_s^\eps|^{mq}\Big)\dif s\right)\leq C_1n^{q(1-(\vartheta\wedge1))}.
$$
Consequently, we get
$$
\cQ_1(\eps)\leq C_1\Big(\alpha_\eps^q+\alpha_\eps^{2q}n^{q(1-(\vartheta\wedge1))}\Big).
$$		
	To control the second term, by  (\ref{n222}) and the assumptions that $F, G\in L^\infty_p$, we deduce that
	\begin{align*}
	&\|\big(\p_s+\sL_1\big)\Phi_n^f(\cdot,x,\cdot)\|_\infty\leq C_2\big(1+|x|^m\big)\\
	&\qquad\times\Big(\|\p_s\Phi_n^f(\cdot,x,\cdot)\|_\infty +\sum_{\ell=1,2}\big\|\nabla_y^\ell\Phi_n^f(\cdot,x,\cdot)\big\|_\infty\Big)\leq C_2n^{2-\vartheta}(1+|x|^{2m}).
	\end{align*}
Taking into account  (\ref{ME}) yields
	$$
	\cQ_{2}(\eps)\leq C_2{\alpha_\eps^{2q}}n^{q(2-\vartheta)}\mE\left(\int_0^t\big(1+|X_s^\eps|^{2mq}\big)\dif s\right)\leq C_2{\alpha_\eps^{2q}}n^{q(2-\vartheta)}.
	$$
	Using the same argument as above, one can check that
	\begin{align*}
	\|\sL_2\Phi_n^f(\cdot,x,\cdot)\|_\infty\leq C_3(1+|x|^m)\|\nabla_y\Phi^f_n(\cdot,x,\cdot)\|_\infty\leq C_3 n^{1-(\vartheta\wedge1)}(1+|x|^{2m}).
	\end{align*}
	Thus we have
	\begin{align*}
	\cQ_{3}(\eps)
	&\leq C_3\frac{\alpha_\eps^{2q}}{\gamma_\eps^q} n^{q(1-(\vartheta\wedge1))}\mE\left(\int_0^t\big(1+|X_s^\eps|^{2mq}\big)\dif s\right)\leq C_3\frac{\alpha_\eps^{2q}}{\gamma_\eps^q} n^{q(1-(\vartheta\wedge1))}.
	\end{align*}
	Furthermore, by the assumption that $c\in L^\infty_p$, it follows directly that
	$$
	\cQ_{4}(\eps)\leq C_4\frac{\alpha_\eps^{2q}}{\beta_\eps^q}\mE\left(\int_0^t\big(1+|X_s^\eps|^{2mq}\big)\dif s\right)\leq C_4\frac{\alpha_\eps^{2q}}{\beta_\eps^q}.
	$$
	Finally, since $\nabla_x^2\Phi^f\in C_p^{\vartheta/2,\delta,\vartheta}$ and due to the fact that
	$$
	\nabla_x^2(\Phi^f_n)=(\nabla_x^2\Phi^f)*\rho_1^n*\rho_2^n,
	$$
	we  derive by (\ref{n111}) that
	\begin{align*}
	\cQ_{5}(\eps)
	&\leq C_5  \,\mE\left(\int_0^t\!\sum_{\ell=1,2}\big\|(\nabla_x^\ell\Phi_n^f-\nabla_x^\ell\Phi^f)(\cdot,X_s^\eps,\cdot)\big\|^q_\infty\dif s\right)\\
	&\leq C_5n^{-q\vartheta}\mE\left(\int_0^t\big(1+|X_s^\eps|^{mq}\big)\dif s\right)\leq C_5 n^{-q\vartheta}.
	\end{align*}
	Combining the above computations, we arrive at
	$$
	\cQ(\eps)\leq C_6\Big(\alpha_\eps+\alpha_\eps^{2}n^{1-(\vartheta\wedge1)}+{\alpha_\eps^2}n^{2-\vartheta}+\frac{\alpha_\eps^2}{\gamma_\eps} n^{1-(\vartheta\wedge1)}+\frac{\alpha_\eps^2}{\beta_\eps}+n^{-\vartheta}\Big)^q.
	$$	
	Taking $n={\alpha_\eps^{-1}}$, we thus get
	$$
	\cQ(\eps)\leq C_7\Big(\alpha_\eps+\alpha_\eps^{\vartheta}+\frac{\alpha_\eps^2}{\beta_\eps}\Big)^q\leq C_7\Big(\alpha_\eps^{\vartheta\wedge1}+\frac{\alpha_\eps^2}{\beta_\eps}\Big)^q,
	$$
	and the proof is finished.
\end{proof}

The above result can be regraded as a law of large numbers type fluctuation estimate. Now, we derive a central limit  type fluctuation estimate for the integral of $f(s,X_s^\eps,Y_s^\eps)$ over the time interval $[0,t]$. This will depend on the two regimes described in (\ref{regime}). Recall that $\Phi^f$ is the solution to the Poisson equation (\ref{pn}).
For simplify, we set
\begin{align*}
\overline{c\cdot\nabla_x\Phi^f}(t,y)&:=\int_{\mR^{d_1}}\!c(x,y)\cdot\nabla_x\Phi^f(t,x,y)\mu^{y}(\dif x).
\end{align*}
The following result will play an important role below.

\bl\label{key22}
Let {\bf (A$_\sigma$)}, {\bf (A$_b$)} and (\ref{ac})  hold. Assume that $b, \sigma\in C_b^{\delta,\vartheta}$ with $0<\delta,\vartheta\leq 2$,  $F, H, G\in L^\infty_p$ and $\lim_{\eps\to0}\alpha_\eps^\vartheta/\gamma_\eps=0$. Then for every $f\in C_p^{\vartheta/2,\delta,\vartheta}$  satisfying (\ref{cen6}) and any $q\geq 2$, the following hold:

\vspace{1mm}
\noindent(i) (Regime 1) if $c\in L^\infty_p$,  we have
\begin{align*}
\mE\left|\frac{1}{\gamma_\eps}\int_0^tf(s,X_s^\eps,Y_s^\eps)\dif s\right|^q\leq C_t\Big(\frac{\alpha_\eps^{\vartheta\wedge1}}{\gamma_\eps}+\frac{\alpha_\eps^2}{\beta_\eps\gamma_\eps}\Big)^q;
\end{align*}
(ii) (Regime 2)  if $c\in C_p^{\delta,\vartheta}$, we have
\begin{align*}
\mE\left|\frac{1}{\gamma_\eps}\int_0^tf(s,X_s^\eps,Y_s^\eps)\dif s-\int_0^t\overline{c\cdot\nabla_x\Phi^f}(s,Y_s^\eps)\dif s\right|^q\leq C_t\Big(\frac{\alpha_\eps^{\vartheta\wedge1}}{\gamma_\eps}+\frac{\alpha_\eps^2}{\beta_\eps}\Big)^q,
\end{align*}
where $C_t>0$ is a constant independent of $\delta, \eps$.
\el

\begin{proof}
We provide the proof for the two regimes separately.

\vspace{1mm}
\noindent
{\it (i) (Regime 1)} This case follows by Lemma \ref{key} directly since   no homogenization occurs. In fact, we have
\begin{align*}
\mE\left|\frac{1}{\gamma_\eps}\int_0^tf(s,X_s^\eps,Y_s^\eps)\dif s\right|^q= \frac{1}{\gamma_\eps^q}\mE\left|\int_0^tf(s,X_s^\eps,Y_s^\eps)\dif s\right|^q \leq C_1\frac{1}{\gamma_\eps^q}\Big(\alpha_\eps^{\vartheta\wedge1}+\frac{\alpha_\eps^2}{\beta_\eps}\Big)^q,
\end{align*}	
which in turn yields the desired result.

\vspace{1mm}
	\noindent{\it (ii) (Regime 2)}
	Recall that $\alpha_\eps^2=\beta_\eps\gamma_\eps$ in this case. As in the proof of Lemma \ref{key}, by (\ref{ito}) we have
	\begin{align*}
	 \hat\cQ(\eps)&:=\mE\left|\frac{1}{\gamma_\eps}\int_0^tf(s,X_s^\eps,Y_s^\eps)\dif s-\int_0^t\overline{c\cdot\nabla_x\Phi^f}(s,Y_s^\eps)\dif s\right|^q\\
&\leq \!C_2\bigg[\frac{\alpha_\eps^{2q}}{\gamma_\eps^q}\Big(|\Phi_n^f(0,x,y)|^q\!+\mE|\Phi^f_n(t,X_t^\eps,Y_t^\eps)|^q\Big)\!\!+ \!\frac{\alpha_\eps^{q}}{\gamma_\eps^q}\mE|M_n^1(t)|^q+\frac{\alpha_\eps^{2q}}{\gamma_\eps^q}\mE|M_n^2(t)|^q\bigg]\\
	 &+C_2\frac{\alpha_\eps^{2q}}{\gamma_\eps^q}\mE\left|\int_0^{t}\!\big(\p_s+\sL_1\big)\Phi_n^f(s,X_s^\eps,Y_s^\eps)\dif s\right|^q\!+C_2\frac{\alpha_\eps^{2q}}{\gamma_\eps^{2q}}\mE\left|\int_0^{t}\!\sL_2\Phi_n^f(s,X_s^\eps,Y_s^\eps)\dif s\right|^q\\
	 &+C_2\frac{1}{\gamma_\eps^{q}}\mE\left|\int_0^{t}\big(\sL_0\Phi_n^f-\sL_0\Phi^f\big)(s,X_s^\eps,Y_s^\eps)\dif s\right|^q\\
	&+C_2 \mE\left|\int_0^{t}\sL_3\Phi_n^f(s,X_s^\eps,Y_s^\eps)-\sL_3\Phi^f(s,X_s^\eps,Y_s^\eps)\dif s\right|^q\\
	 &+C_2\mE\left|\int_0^{t}\sL_3\Phi^f(s,X_s^\eps,Y_s^\eps)-\overline{c\cdot\nabla_x\Phi^f}(s,Y_s^\eps)\dif s\right|^q=:\sum_{i=1}^6\hat\cQ_i(\eps),
	\end{align*}
	where $\Phi_n^f$ is the mollifyer of $\Phi^f$ defined as in (\ref{fn}).
Following exactly the same arguments as in the proof of Lemma \ref{key}, one can check that
	$$
	\sum_{i=1}^5\hat\cQ_i(\eps)\leq C_3\Big(\frac{\alpha_\eps^q}{\gamma_\eps^q}+\frac{\alpha_\eps^{2q}}{\gamma_\eps^q}n^{q(2-\vartheta)}+\frac{\alpha_\eps^{2q}}{\gamma_\eps^{2q}}n^{q(1-(\vartheta\wedge1))}+\frac{1}{\gamma_\eps^q}n^{-q\vartheta}\Big).
	$$
	Taking $n=\alpha_\eps^{-1}$, we further get
		$$
	\sum_{i=1}^5\hat\cQ_i(\eps)\leq C_4\Big(\frac{\alpha_\eps^q}{\gamma_\eps^q}+\frac{\alpha_\eps^{q\vartheta}}{\gamma_\eps^q}\Big)\leq C_4\frac{\alpha_\eps^{q(\vartheta\wedge1)}}{\gamma^q_\eps}.
	$$	
	For the last term, note  that by definition we have
	$$
	 \sL_3\Phi^f(t,x,y)-\overline{c\cdot\nabla_x\Phi^f}(t,y)=c(x,y)\cdot\nabla_x\Phi^f(t,x,y)-\overline{c\cdot\nabla_x\Phi^f}(t,y),
	$$
	which satisfies the centering condition (\ref{cen6}). Furthermore, since $c\in C_p^{\delta,\vartheta}$, $\Phi^f\in C_p^{\vartheta/2,2+\delta,\vartheta}$ and by Lemma \ref{coe}, we  have $c\cdot\nabla_x\Phi^f-\overline{c\cdot\nabla_x\Phi^f}\in C_p^{\vartheta/2,\delta,\vartheta}$. As a direct result of Lemma \ref{key}, we  get
	\begin{align*}
	\hat\cQ_6(\eps)\leq C_5\Big(\alpha_\eps^{\vartheta\wedge1}+\frac{\alpha_\eps^2}{\beta_\eps}\Big)^q.
	\end{align*}
	Consequently, we arrive at
	$$
	\hat\cQ(\eps)\leq C_6\Big(\frac{\alpha_\eps^{\vartheta\wedge1}}{\gamma_\eps}+\alpha_\eps^{\vartheta\wedge1} +\frac{\alpha_\eps^2}{\beta_\eps}\Big)^q\leq C_6\Big(\frac{\alpha_\eps^{\vartheta\wedge1}}{\gamma_\eps} +\frac{\alpha_\eps^2}{\beta_\eps}\Big)^q,
	$$
	and the proof is finished.
\end{proof}

\subsection{Proof of Theorem \ref{main1}}

Throughout this subsection, we  fix  $T>0$ and always assume that the conditions in Theorem \ref{main1} hold.
Recall that $\bar Y_t^k$ ($k=1,2$) are given by (\ref{sde000}). Due to our low regularity assumptions on the drift coefficients of  system (\ref{sde0}), it seems to be not possible to prove the strong convergence of $Y_t^\eps$ to $\bar Y_t^k$ directly. For this reason, we  shall use Zvonkin's argument to transform the equations for $Y_t^\eps$ and $\bar Y_t$ into new ones.

Consider the following backward PDE on $[0,T]\times\mR^{d_2}$: for $k=1,2$,
\begin{equation}\left\{\begin{array}{l}\label{PDE2}
\displaystyle
\p_tv_k(t,y)+\bar\sL_kv_k(t,y)+\bar F_k(t,y)=0,\quad t\in [0, T),\\
v_k(T,y)=0,
\end{array}\right.
\end{equation}
where $\bar\sL_k$ are defined by (\ref{LL1}). By Lemma \ref{coe},  we have for every $k=1,2$, $\bar F_k\in C_b^{\vartheta/2,\vartheta}$ with $0<\vartheta\leq 1$.
Thus under {\bf (A$_{G}$)}, it is well known that there exits a unique solution $v_k\in L^\infty\big([0,T];C^{2+\vartheta}_{b}(\mR^{d_2})\big)\cap C_b^{1+\vartheta/2}\big([0,T];L^\infty(\mR^{d_2})\big)$ for equation \eqref{PDE2}, see e.g. \cite[Chapter IV, Section 5]{La-So-Ur}. Moreover, we can choose $T$ small enough so that for any $0<t\leq T$ (see e.g. \cite[Theorem 2]{FF} or \cite[Theorem 2]{FGP}),
$$
|\nabla_yv_k(t,y)|\leq 1/2,\quad\forall y\in\mR^{d_2}.
$$
Define the transformation function by
$$
\Gamma_k(t,y):=y+v_k(t,y).
$$
Then the map $y \mapsto\Gamma_k(t,y)$ is a $C^1$-diffeomorphism and for every $t\in(0,T]$ and $y\in\mR^{d_2}$,
\begin{align}\label{dd}
1/2\leq|\nabla_y\Gamma_k(t,y)|\leq 3/2.
\end{align}
Let us define two new processes by
\begin{align}\label{vv}
\bar V^k_t:=\Gamma_k(t,\bar{Y}^k_t)\quad\text{and}\quad V^{k,\eps}_t:=\Gamma_k(t,Y^{\eps}_t).
\end{align}
We have the following result.

\bl[Zvonkin's transformation]\label{zt}
Let $\bar V_t^k$ and $V_t^{k,\eps}$ be defined by (\ref{vv}). Then we have for $k=1,2$,
\begin{align}\label{z1}
\dif \bar V_t^k=G(t,\bar Y^k_t)\nabla_y \Gamma_k(t,\bar{Y}^k_t)\dif W^2_t,\qquad V^k_0=\Gamma_k(0, y)
\end{align}
and
\begin{align}\label{z2}
\dif V^{k,\eps}_t=\big[F(t,X^{\eps}_t, Y^{\eps}_t)&-\bar{F}_k(t,Y^{\eps}_t)\big]\nabla_y \Gamma_k(t,Y^{\eps}_t)\dif t+\gamma_\eps^{-1}H(t,X_t^\eps,Y_t^\eps)\nabla_y \Gamma_k(t,Y^{\eps}_t)\dif t\no\\
&+G(t,Y_t^\eps)\nabla_y \Gamma_k(t,Y^{\eps}_t)\dif W^2_t,\qquad V_0^\eps=\Gamma_k(0, y).
\end{align}
\el
\begin{proof}
	Using It\^o's formula, we have for $k=1,2$,
	\begin{align*}
	 v_k(t,Y_t^\eps)&=v_k(0,y)+\int_0^t\big(\p_s+\sL_1\big)v_k(s,Y_s^\eps)\dif s+\frac{1}{\gamma_\eps}\int_0^t\sL_2v_k(s,Y_s^\eps)\dif s\\
	&\quad+\int_0^tG(s,Y_s^\eps)\nabla_yv_k(s,Y_s^\eps)\dif W_s^2\\
	&=v_k(0,y)+\int_0^t\big(\p_s+\bar\sL_k\big)v_k(s,Y_s^\eps)\dif s\\
	&\quad+\int_0^t\big[F(s,X^{\eps}_s, Y^{\eps}_s)-\bar{F}_k(s,Y^{\eps}_s)\big]\nabla_y v_k(s,Y^{\eps}_s)\dif s\\
	&\quad+\frac{1}{\gamma_\eps}\int_0^t\sL_2v_k(s,Y_s^\eps)\dif s+\int_0^tG(s,Y_s^\eps)\nabla_yv_k(s,Y_s^\eps)\dif W_s^2,
	\end{align*}
	where $\sL_1$ and $\sL_2$ are given by (\ref{lll}). Taking into account (\ref{PDE2}), we further get
	\begin{align*}
	v_k(t,Y_t^\eps)
	&=v_k(0,y)-\!\int_0^t\bar F_k(s,Y_s^\eps)\dif s+\!\int_0^t\!\big[F(s,X^{\eps}_s, Y^{\eps}_s)-\bar{F}_k(s,Y^{\eps}_s)\big]\nabla_y v_k(s,Y^{\eps}_s)\dif s\\
	&\quad+\frac{1}{\gamma_\eps}\int_0^t\sL_2v_k(s,Y_s^\eps)\dif s+\int_0^tG(s,Y_s^\eps)\nabla_yv_k(s,Y_s^\eps)\dif W_s^2,
	\end{align*}
	This together with the equation for $Y_t^\eps$ and the fact that $\nabla_y\Gamma_k(t,y)=\mI_{d_2}+\nabla_y v_k(t,y)$ yields (\ref{z2}). The proof of (\ref{z1}) is easier and follows by the same argument.
\end{proof}

Now, we are in the position to give:

\begin{proof}[Proof of Theorem \ref{main1}]
	Let us first assume that $T>0$ is sufficiently small so that (\ref{dd}) holds.
	By the definition (\ref{vv}), we have for any $t\in [0, T]$, $k=1,2$ and $q\geq 2$ that
	\begin{align}\label{aaa}
	\mE\big|Y_{t}^\eps-\bar Y^k_{t}\big|^q\leq C_q\,\mE\big|V_{t}^{k,\eps}-\bar V^k_{t}\big|^q.
	\end{align}
	In view of  (\ref{z1}) and (\ref{z2}), we may write
	\begin{align*}
	V_t^{k,\eps}-\bar V^k_t&=\int^t_0\big[G(s,Y_s^\eps)\nabla_y \Gamma_k(s,Y^{\eps}_s)-G(s,\bar Y^k_s)\nabla_y \Gamma_k(s,\bar Y^k_s)\big]\dif W^2_s\\
	&\quad+\int_0^t\big[F(s,X^{\eps}_s, Y^{\eps}_s)-\bar{F}_k(s,Y^{\eps}_s)\big]\nabla_y \Gamma_k(s,Y^{\eps}_s)\dif s\\
	&\quad+\frac{1}{\gamma_\eps}\int_0^tH(s,X_s^\eps,Y_s^\eps)\nabla_y \Gamma_k(s,Y^{\eps}_s)\dif s.
	\end{align*}
	Taking expectation of both sides of the above equality, we get that there exists a constant $C_0>0$ such that
	\begin{align*}
	\mE|V_{t}^{k,\eps}-\bar V^k_{t}|^q&\leq C_0\mE\left(\int^{t}_0\big|G(s,Y_s^\eps)\nabla_y \Gamma_k(s,Y^{\eps}_s)-G(s,\bar Y^k_s)\nabla_y \Gamma_k(s,\bar Y^k_s)\big|^2\dif s\right)^{q/2}\\
	&\quad+C_0\mE\bigg|\int_0^t\big[F(s,X^{\eps}_s, Y^{\eps}_s)-\bar{F}_k(s,Y^{\eps}_s)\big]\nabla_y \Gamma_k(s,Y^{\eps}_s)\dif s\\
	 &\quad\quad\qquad+\frac{1}{\gamma_\eps}\int_0^tH(s,X_s^\eps,Y_s^\eps)\nabla_y \Gamma_k(s,Y^{\eps}_s)\dif s\bigg|^q=:\sQ_1(t,\eps)+\sQ_2(t,\eps).
	\end{align*}
	Note that for each $k=1,2$,
	$$
	y\mapsto G(t,\cdot)\nabla_y\Gamma_k(t,\cdot)\in  C_b^1(\mR^{d_2}).
	$$
	Thus we   have
	\begin{align}\label{ae}
	\sQ_1(t,\eps)\leq C_1\mE\left(\int^{t}_0|Y_s^\eps-\bar Y^k_s|^2\dif s\right)^{q/2}\leq C_t \,\mE\left(\int^{t}_0|Y_s^\eps-\bar Y^k_s|^q\dif s\right),
	\end{align}
	where $C_t>0$ is a constant independent of $\eps$.
		Below, we proceed to control the second term according to Regime 1 and Regime 2 in (\ref{regime}) separately.
		
	\vspace{1mm}
	\noindent
	{\bf Regime 1 ($k=1$).}
	In this case, note that  the function
	$[F(t,x,y)-\bar F_1(t,y)]\nabla_y\Gamma_1(t,y)$ satisfies the centering condition (\ref{cen6}) and belongs to $C_p^{\vartheta/2,\delta,\vartheta}$. Thus by Lemma \ref{key} we have
	\begin{align*}
	\mE\bigg|\int_0^t\big[F(s,X^{\eps}_s, Y^{\eps}_s)-\bar{F}_1(s,Y^{\eps}_s)\big]\nabla_y \Gamma_1(s,Y^{\eps}_s)\dif s\bigg|^q\leq C_2\Big(\alpha_\eps^{\vartheta\wedge1}+\frac{\alpha_\eps^2}{\beta_\eps}\Big)^q.
	\end{align*}
	At the same time, thanks to assumption {\bf (A$_H$)}, the function
	$H(t,x,y)\nabla_y\Gamma_1(t,y)$ also satisfies the centering condition (\ref{cen6}). It then follows by Lemma \ref{key22} (i) that
		\begin{align*}
	\mE\bigg|\frac{1}{\gamma_\eps}\int_0^tH(s,X_s^\eps,Y_s^\eps)\nabla_y \Gamma_1(s,Y^{\eps}_s)\dif s\bigg|^q\leq C_3\Big(\frac{\alpha_\eps^{\vartheta\wedge1}}{\gamma_\eps}+\frac{\alpha_\eps^2}{\beta_\eps\gamma_\eps}\Big)^q.
	\end{align*}
	Combining the above estimates, we get
	$$
	\sQ_2(t,\eps)\leq C_4\Big(\frac{\alpha_\eps^{\vartheta\wedge1}}{\gamma_\eps}+\frac{\alpha_\eps^2}{\beta_\eps\gamma_\eps}\Big)^q.
	$$
	
	Now, in view of (\ref{aaa}) and (\ref{ae}), we arrive at
	$$
	\mE\big|Y_{t}^\eps-\bar Y_{t}^1\big|^q\leq C_5\mE\left(\int^{t}_0|Y_s^\eps-\bar Y^1_s|^q\dif s\right)+C_6\Big(\frac{\alpha_\eps^{\vartheta\wedge1}}{\gamma_\eps}+\frac{\alpha_\eps^2}{\beta_\eps\gamma_\eps}\Big)^q,
	$$
	which in turn yields the desired result by Gronwall's inequality.
	For general $T>0$, the result can be proved by induction and analogous arguments.
	
	\vspace{1mm}
	\noindent
	{\bf Regime 2 ($k=2$).}
	In this case, recall that we have $\bar F_2(t,y)=\bar F_1(t,y)+\overline{c\cdot\nabla_x\Phi}(t,y)$, where $\Phi$ solves the Poisson equation (\ref{pde1}) and $\overline{c\cdot\nabla_x\Phi}$ is given by (\ref{cphi}). Thus, we deduce that
		\begin{align*}
	\sQ_2(t,\eps)\!&\leq C_2\mE\bigg|\int_0^t\big[F(s,X^{\eps}_s, Y^{\eps}_s)-\bar{F}_1(s,Y^{\eps}_s)\big]\nabla_y \Gamma_2(s,Y^{\eps}_s)\dif s\bigg|^q\\
	 &\quad+C_2\mE\bigg|\frac{1}{\gamma_\eps}\!\int_0^t\!H(s,X_s^\eps,Y_s^\eps)\nabla_y \Gamma_2(s,Y^{\eps}_s)\dif s\!-\!\!\int_0^t\!\overline{c\cdot\nabla_x\Phi}(s,Y_s^\eps)\nabla_y \Gamma_2(s,Y^{\eps}_s)\dif s \bigg|^q\\
	&=:\sQ_{21}(t,\eps)+\sQ_{22}(t,\eps).
	\end{align*}
Following exactly the same arguments as above, we  get
\begin{align*}
\sQ_{21}(t,\eps)\leq C_3\Big(\alpha_\eps^{\vartheta\wedge1}+\frac{\alpha_\eps^2}{\beta_\eps}\Big)^q.
\end{align*}
On the other hand, let $\tilde \Phi(t,x,y):=\Phi(t,x,y)\cdot\nabla_y\Gamma_2(t,y)$. Then one can check  that
$$
\sL_0(x,y)\tilde\Phi(t,x,y)=-H(t,x,y)\cdot\nabla_y\Gamma_2(t,y)
$$
and
$$
\overline{c\cdot\nabla_x\tilde\Phi}(t,y):=\int_{\mR^{d_1}}c(x,y)\cdot\nabla_x\tilde\Phi(t,x,y)\mu^y(\dif x)=\overline{c\cdot\nabla_x\Phi}(t,y)\cdot\nabla_y\Gamma_2(t,y).
$$
Consequently, by  Lemma \ref{key22} (ii)  we have
\begin{align*}
\sQ_{22}(t,\eps)\leq C_4\Big(\frac{\alpha_\eps^{\vartheta\wedge1}}{\gamma_\eps}+\frac{\alpha_\eps^2}{\beta_\eps}\Big)^q.
\end{align*}
Combining the above computations, we arrive at
$$
\mE\big|Y_{t}^\eps-\bar Y^2_{t}\big|^q\leq C_5\mE\left(\int^{t}_0|Y_s^\eps-\bar Y^2_s|^q\dif s\right)+C_6\Big(\frac{\alpha_\eps^{\vartheta\wedge1}}{\gamma_\eps}+\frac{\alpha_\eps^2}{\beta_\eps}\Big)^q,
$$
which in turn yields the desired result by Gronwall's inequality. Hence the whole proof is finished.
\end{proof}

\section{Central limit theorem without homogenization}

In this section, we study the central limit theorem for SDE (\ref{sde99}).  We shall first derive some weak  fluctuation estimates in Subsection 5.1. Then we prove Theorem \ref{main3} in Subsection 5.2.

\subsection{Weak fluctuation estimate (i)}

Recall that $Y^{0,\eps}_t$ converges strongly to $\bar Y^1_t$, and $Z_t^{0,\eps}$, $\bar Z^0_{\ell,t}$ ($\ell=1,2,3$) are defined by (\ref{z0}), (\ref{sdez0}), respectively. To prove the weak convergence of $Z_t^{0,\eps}$ to $\bar Z^0_{\ell,t}$, we shall view the process $(X^{\eps}_t, Y^{0,\eps}_t, \bar Y^1_t, Z_t^{0,\eps}, \bar Z^0_{\ell,t})$ as a whole system, i.e., we consider
\begin{equation*}
\left\{ \begin{aligned}
&\dif X^{\eps}_t =\alpha_\eps^{-2}b(X^{\eps}_t,Y^{0,\eps}_t)\dif t+\beta_\eps^{-1}c(X^{\eps}_t,Y^{0,\eps}_t)\dif t+\alpha_\eps^{-1}\sigma(X^{\eps}_t,Y_t^{0,\eps})\dif W^{1}_t,\qquad\,\,\,\,\,\! X_0=x,\\
&\dif Y^{0,\eps}_t =F(t,X^{\eps}_t, Y^{0,\eps}_t)\dif t+G(t,Y_t^{0,\eps})\dif W^{2}_t,\quad\qquad\qquad\qquad\qquad\qquad\quad\quad\,\,\, Y^{0,\eps}_0=y,\\
&\dif \bar Y^1_t=\bar F_1(t,\bar Y^1_t)\dif t+G(t,\bar Y^1_t)\dif W^2_t,\quad\qquad\quad\qquad\quad\qquad\quad\qquad\quad\quad\quad\quad\,\,\,\,\! \bar Y_0^1=y,\\
&\dif Z_t^{0,\eps}\!=\frac{1}{\eta_\eps}\Big[F(t,X_t^\eps,Y_t^{0,\eps})\!-\!\bar F_1(t,\bar Y^1_t)\Big]\dif t
+\!\frac{1}{\eta_\eps}\Big[G(t,Y_t^{0,\eps})\!-\!G(t,\bar Y^1_t)\Big]\dif W_t^2, \,\,Z_0^{0,\eps}=0,\\
&\dif \bar Z^0_{\ell,t}=\varrho^0_\ell(t,\bar Y_t^1,\bar Z^0_{\ell,t})\dif t+\chi_\ell^0(t,\bar Y_t^1,\bar Z^0_{\ell,t})\dif W_t^2, \quad\qquad\quad\qquad\quad\quad\!\!\!\!\!\quad\qquad\quad\quad\bar Z_{\ell, 0}^{1}=0,
\end{aligned} \right.
\end{equation*}
where for  $\ell=1,2,3$, $\eta_\eps$ is given in (\ref{regime4}), and $\varrho^0_\ell(t,y,z)$ and $\chi_\ell^0(t,y,z)$ denote  the drift and diffusion coefficients for $\bar Z^0_{\ell,t}$ in SDE (\ref{sdez0}), respectively. We write $\varrho^0_\ell$ and $\chi_\ell^0$ here just for simplicity, and we shall  not use them below.

Note that by definition, we can write
\begin{align*}
\dif Z_t^{0,\eps}
&=\frac{1}{\eta_\eps}\Big[F(t,X_t^\eps,Y_t^{0,\eps})-\bar F_1(t,Y_t^{0,\eps})\Big]\dif t\\
&\quad+\bigg(\frac{1}{\eta_\eps}\Big[\bar F_1(t,Y_t^{0,\eps})-\bar F_1(t,\bar Y^1_t)\Big]\dif t+\frac{1}{\eta_\eps}\Big[G(t,Y_t^{0,\eps})-G(t,\bar Y^1_t)\Big]\dif W_t^2\bigg).
\end{align*}
To shorten the notation, we define
\begin{align}
\sL^1_4&:=\sL^1_4(t,x,y,z):=\sum_{i=1}^{d_2}\big[F^i(t,x,y)-\bar F^i_1(t,y)\big]\frac{\p}{\p z_i},\label{LL104}
\end{align}
and for $k=1,2$,
\begin{align}\label{LL2}
\begin{split}
\sL_5^{k,\eps}&:=\sL_5^{k,\eps}(t,y,\bar y,z):=\eta_\eps^{-1}\sum_{i=1}^{d_2}\big[\bar F^i_k(t,y)-\bar F^i_k(t,\bar y)\big]\frac{\p}{\p z_i}\\
&\quad+\eta_\eps^{-2}/2\sum_{i,j=1}^{d_2}\Big(\big[G(t,y)-G(t,\bar y)\big]\big[G(t,y)-G(t,\bar y)\big]^*\Big)^{ij}\frac{\p^2}{\p z_i\p z_j}\\
&\quad+\eta_\eps^{-1}/2\sum_{i,j=1}^{d_2}\Big(G(t,y)\big[G(t,y)-G(t,\bar y)\big]^*\Big)^{ij}\frac{\p^2}{\p y_i\p z_j}.
\end{split}
\end{align}
Let $f(t,x,y,z)$ be a function satisfying the centering condition, i.e.,
\begin{align}\label{cen30}
\int_{\mR^{d_1}}f(t,x,y,z)\mu^y(\dif x)=0,\quad\forall (t,y,z)\in\mR_+\times\mR^{d_2}\times\mR^{d_2}.
\end{align}
Let $\Phi^f(t,x,y,z)$ denote the unique solution to the following Poisson equation:
\begin{align}\label{xyz0}
\sL_0(x,y)\Phi^f(t,x,y,z)=-f(t,x,y,z),
\end{align}
where $(t, y, z)\in\mR_+\times\mR^{d_2}\times\mR^{d_2}$ are regarded as parameters. We have the following fluctuation estimate for the process $(X_t^\eps,Y_t^{0,\eps},Z_t^{0,\eps})$.

\bl\label{key30}
Let {\bf (A$_\sigma$)}, {\bf (A$_b$)}, (\ref{ac})
hold and $0<\delta,\vartheta\leq 2$. Assume that $b,\sigma\in C_b^{\delta,1\vee\vartheta}$, $c\in L^\infty_p$, $F\in C_p^{0,\delta,1}$ and $G\in C_b^{0,1}$. Then for every $f\in C_p^{\vartheta/2,\delta,\vartheta,2}$  satisfying (\ref{cen30}), we have
\begin{align}
\mE\left(\int_0^tf(s,X_s^\eps,Y_s^{0,\eps},Z_s^{0,\eps})\dif s\right)&\leq C_t\bigg[ \alpha_\eps^{\vartheta}
+\frac{\alpha_\eps^2}{\beta_\eps}\mE\left(\int_0^t\sL_3\Phi^f(s,X_s^\eps,Y_s^{0,\eps},Z_s^{0,\eps})\dif s\right)\no\\
&\quad+\frac{\alpha_\eps^2}{\eta_\eps}\mE\left(\int_0^t\sL_4^1\Phi^f(s,X_s^\eps,Y_s^{0,\eps},Z_s^{0,\eps})\dif s\right)\bigg].\label{we1}
\end{align}
where $\sL_3$ and $\sL_4^1$ are given by (\ref{lll}) and (\ref{LL104}), respectively, and $C_t>0$ is a constant independent of $\delta, \eps$.
\el

\br
(i) Note that under the above assumptions and according
to Theorem \ref{popde}, we in fact have $\Phi^f\in C_p^{\vartheta/2,2+\delta,\vartheta,2}$. Thus we   get
\begin{align*}
&\mE\left(\int_0^t\sL_3\Phi^f(s,X_s^\eps,Y_s^{0,\eps},Z_s^{0,\eps})\dif s\right)+\mE\left(\int_0^t\sL_4^1\Phi^f(s,X_s^\eps,Y_s^{0,\eps},Z_s^{0,\eps})\dif s\right)\\
&\leq C_0\mE\left(\int_0^t\big(1+|X_s^\eps|^{2m}\big)\dif s\right)<\infty,
\end{align*}
which in turn yields that
\begin{align}\label{pp0}
\mE\left(\int_0^tf(s,X_s^\eps,Y_s^{0,\eps},Z_s^{0,\eps})\dif s\right)&\leq C_t\Big(\alpha_\eps^{\vartheta}+\frac{\alpha_\eps^2}{\beta_\eps}+\frac{\alpha_\eps^2}{\eta_\eps}\Big).
\end{align}
However,  the homogenization effects of the last two terms in (\ref{we1}) will appear when we study the central limit theorems, so we just keep them on the right hand side for later use.

\vspace{1mm}
(ii) Compared with Lemma \ref{key}, we do not take  absolute value for the integral over $[0,t]$ on the left hand side of (\ref{we1}). We shall see  that the  involved martingale part will not play any role in the control of the error bound in this case.
\er

\begin{proof}
	Let $\Phi^f_{n}$ be the mollifyer of $\Phi^f$ defined as in (\ref{fn}).
	Using It\^o's formula, we  have
	\begin{align*}
	&\Phi^f_{n}(t,X_t^\eps,Y_t^{0,\eps},Z_t^{0,\eps})=\Phi^f_{n}(0,x,y,0) +\frac{1}{\alpha_\eps}\tilde M_{n}^1(t)+\tilde M_n^2(t)+\frac{1}{\eta_\eps}\tilde M_{n}^3(t)\\
	 &\quad\quad+\int_0^t\big(\p_s+\sL_1+\sL_5^{1,\eps}\big)\Phi^f_{n}(s,X_s^\eps,Y_s^{0,\eps},Z_s^{0,\eps})\dif s\\
	 &\quad\quad+\frac{1}{\alpha_\eps^2}\int_0^t\sL_0\Phi^f_{n}(s,X_s^\eps,Y_s^{0,\eps},Z_s^{0,\eps})\dif s+\frac{1}{\beta_\eps}\int_0^t\sL_3\Phi^f_{n}(s,X_s^\eps,Y_s^{0,\eps},Z_s^{0,\eps})\dif s\\
	 &\quad\quad+\frac{1}{\eta_\eps}\int_0^t\sL_4^1\Phi^f_{n}(s,X_s^\eps,Y_s^{0,\eps},Z_s^{0,\eps})\dif s.
	\end{align*}
	where $\sL_1$ and $\sL_5^{1,\eps}$ are given by (\ref{lll}) and (\ref{LL2}), respectively, and $\tilde M_{n}^1(t)$, $\tilde M_{n}^2(t)$, $\tilde M_{n}^3(t)$ are martingales defined by
	\begin{align*}
	\tilde M_{n}^1(t)&:=\int_0^t\sigma(X_s^\eps,Y_s^{0,\eps})\nabla_x\Phi^f_n(s,X_s^\eps,Y_s^{0,\eps},Z_s^{0,\eps})\dif W^1_s,\\
	\tilde M_{n}^2(t)&:=\int_0^tG(s,Y_s^{0,\eps})\nabla_y\Phi^f_n(s,X_s^\eps,Y_s^{0,\eps},Z_s^{0,\eps})\dif W^2_s,\\
	\tilde M_{n}^3(t)&:=\int_0^t\big[G(s,Y_s^{0,\eps})-G(s,\bar Y_s^1)\big]\nabla_z\Phi^f_n(s,X_s^\eps,Y_s^{0,\eps},Z_s^{0,\eps})\dif W^2_s.
	\end{align*}
	Taking expectation and in view of (\ref{xyz0}), we have
	\begin{align*}
	 \cV(\eps)&:=\mE\!\left(\int_0^tf(s,X_s^\eps,Y_s^{0,\eps},Z_s^{0,\eps})\dif s\right)\!=\alpha_\eps^2\,\mE\big[\Phi^f_{n}(0,x,y,0)-\Phi^f_{n}(t,X_t^\eps,Y_t^{0,\eps},Z_t^{0,\eps})\big]\\
	 &\quad+\alpha_\eps^2\mE\left(\int_0^t\big(\p_s+\sL_1+\sL_5^{1,\eps}\big)\Phi^f_{n}(s,X_s^\eps,Y_s^{0,\eps},Z_s^{0,\eps})\dif s\right)\\
	 &\quad+\mE\left(\int_0^t\big[\sL_0\Phi^f_{n}-\sL_0\Phi^f\big](s,X_s^\eps,Y_s^{0,\eps},Z_s^{0,\eps})\dif s\right)\\
	 &\quad+\bigg[\frac{\alpha_\eps^2}{\beta_\eps}\mE\left(\int_0^t\big(\sL_3\Phi^f_{n}-\sL_3\Phi^f\big)(s,X_s^\eps,Y_s^{0,\eps},Z_s^{0,\eps})\dif s\right)\\
	 &\quad+\frac{\alpha_\eps^2}{\eta_\eps}\mE\left(\int_0^t\big(\sL_4^1\Phi^f_{n}-\sL_4^1\Phi^f\big)(s,X_s^\eps,Y_s^{0,\eps},Z_s^{0,\eps})\dif s\right)\bigg]\\
	 &\quad+\frac{\alpha_\eps^2}{\beta_\eps}\mE\left(\int_0^t\sL_3\Phi^f(s,X_s^\eps,Y_s^{0,\eps},Z_s^{0,\eps})\dif s\right)\\
	 &\quad+\frac{\alpha_\eps^2}{\eta_\eps}\mE\left(\int_0^t\sL_4^1\Phi^f(s,X_s^\eps,Y_s^{0,\eps},Z_s^{0,\eps})\dif s\!\right)=:\sum_{i=1}^6\cV_i(\eps).
	\end{align*}
	Let us handle the term involving $\sL_5^{1,\eps}$. Due to the assumptions that $b,\sigma\in C_b^{\delta,1\vee\vartheta}$, $F\in C_p^{0,\delta,1}$, and by  Lemma \ref{coe}, we have  $\bar F_1\in C_b^{0,1}$. This together with the condition  $G\in C_b^{0,1}$ and the mean value theorem yields  that for some $m>0$ and $C_1>0$,
	\begin{align}\label{52}
	 &\mE\left(\int_0^t\sL_5^{1,\eps}\Phi^f_{n}(s,X_s^\eps,Y_s^{0,\eps},Z_s^{0,\eps})\dif s\right)\leq C_1\|\nabla_z^2\Phi^f\|_{L^\infty_p}\mE\bigg(\int_0^t\big(1+|X_s^{\eps}|^m\big)\no\\
	&\qquad\qquad\qquad\qquad\times\bigg[\frac{|\bar F_1(s,Y_s^{0,\eps})-\bar F_1(s,\bar Y_s^1)|}{\eta_\eps}+\frac{|G(s,Y_s^{0,\eps})-G(s,\bar Y_s^1)|^2}{\eta_\eps^2}\bigg]\dif s\bigg)\no\\
	&\quad+C_1\|\nabla_y\nabla_z\Phi^f_n\|_{L^\infty_p}\mE\bigg(\int_0^t\big(1+|X_s^{\eps}|^m\big)\bigg[\frac{|G(s,Y_s^{0,\eps})-G(s,\bar Y_s^1)|}{\eta_\eps}\bigg]\dif s\bigg)\no\\
	&\leq C_1 n^{1-(\vartheta\wedge1)}\mE\left(\int_0^t\big(1+|X_s^{\eps}|^m\big)\big(1+|Z_s^{0,\eps}|^2\big)\dif s\right)\leq C_1 n^{1-(\vartheta\wedge1)}.
	\end{align}
	Then we may argue as in the proof of Lemma \ref{key} to get that
	$$
	\sum_{i=1}^4\cV_i(\eps)\leq C_2 \Big(\alpha_\eps^2+\alpha_\eps^2n^{1-(\vartheta\wedge1)}+\alpha_\eps^2n^{2-\vartheta}+n^{-\vartheta}+\frac{\alpha_\eps^2}{\beta_\eps}n^{-\vartheta} +\frac{\alpha_\eps^2}{\eta_\eps}n^{-\vartheta}\Big).
	$$
	Taking $n=\alpha_\eps^{-1}$, and noticing that $\alpha_\eps^2/\beta_\eps\to0$ and $\alpha_\eps^2/\eta_\eps\to0$ as $\eps\to0$, we get the desired result.
\end{proof}

\subsection{Proof of Theorem \ref{main3}}
We are now in the position to give:

\begin{proof}[Proof of Theorem \ref{main3}]
	For every $\varphi\in C_b^4(\mR^{d_2})$ and $\ell=1,2,3$, let $u^0_\ell$ be the solution to the Cauchy problem  (\ref{PDE}) on $[0,T]\times\mR^{d_2}\times\mR^{d_2}$ with $k=0$, and define
	$$
	\tilde u^0_\ell(t,y,z):=u^0_\ell(T-t,y,z),\quad t\in [0, T].
	$$
	Then for any $y\in\mR^{d_2}$, we have $\tilde u^0_\ell(T,y,z)\equiv\varphi(z)$ and $\tilde u^0_\ell(0,y,0)=\mE\varphi(\bar Z^0_{\ell,T})$. As a result, for $\ell=1,2,3$ we get that
	\begin{align*}
	\cR^0_\ell(\eps):=\mE\varphi(Z_{T}^{0,\eps})-\mE\varphi(\bar Z^0_{\ell,T})=\mE \tilde u^0_\ell(T,Y_T^{0,\eps},Z_T^{0,\eps})-\tilde u^0_\ell(0,y,0).
	\end{align*}
	According to Theorem \ref{pde} and by It\^o's formula,
	\begin{align*}
	\tilde u^0_\ell(T,Y_T^{0,\eps},Z_T^{0,\eps})&=\tilde u^0_\ell(0,y,0)+\int_0^T\!\Big(\p_s+\sL_1+\sL_5^{1,\eps}\Big)\tilde u^0_\ell(s,Y_s^{0,\eps},Z_s^{0,\eps})\dif s\\
	&\quad+\frac{1}{\eta_\eps}\int_0^T\!\!\sL_4^1\tilde u^0_\ell(s,Y_s^{0,\eps},Z_s^{0,\eps})\dif s+M^1_{\ell,T}+\frac{1}{\eta_\eps}M_{\ell,T}^2,
	\end{align*}
	where $M_{\ell,t}^1$ and $M_{\ell,t}^2$ are martingales given by
	\begin{align*}
	M_{\ell,t}^1&:=\int_0^tG(s,Y_s^{0,\eps})\nabla_y\tilde u^0_\ell(s,Y_s^{0,\eps},Z_s^{0,\eps})\dif W^2_s,\\
	M_{\ell,t}^2&:=\int_0^t\big[G(s,Y_s^{0,\eps})-G(s,\bar Y_s^1)\big]\nabla_z\tilde u^0_\ell(s,Y_s^{0,\eps},Z_s^{0,\eps})\dif W^2_s.
	\end{align*}
	Thus we further have
	\begin{align*}
	\cR^0_\ell(\eps)&\leq \mE\left(\int_0^T\!\Big[\p_s+\sL_1+\sL_5^{1,\eps}\Big]\tilde u^0_\ell(s,Y_s^{0,\eps},Z_s^{0,\eps})\dif s\right)\\
	&\quad+\frac{1}{\eta_\eps}\mE\left(\int_0^T\!\!\sL_4^1\tilde u^0_\ell(s,Y_s^{0,\eps},Z_s^{0,\eps})\dif s\right).
	\end{align*}
	Note that by definition,  the functions
		\begin{align*}
	\sL_4^1\tilde u^0_\ell(t,y,z)
	 =\big[F(t,x,y)-\bar F_1(t,y)\big]\nabla_z\tilde u^0_\ell(t,y,z)
	\end{align*}
 satisfy the centering condition (\ref{cen30}).
	Recall that  $\Upsilon(t,x,y)$ solves the   Poisson equation (\ref{pde11}).
	Define for $\ell=1,2,3$,
	$$
	\tilde\Upsilon^0_\ell(t,x,y,z):=\Upsilon(t,x,y)\nabla_z\tilde u^0_\ell(t,y,z).
	$$
	Then
	\begin{align*}
	\sL_0(x,y)\tilde\Upsilon^0_\ell(t,x,y,z)=-\big[F(t,x,y)-\bar F_1(t,y)\big]\nabla_z\tilde u^0_\ell(t,y,z).
	\end{align*}
	Furthermore, by Lemma \ref{coe}  we have $\bar F_1\in C_b^{(1+\vartheta)/2,1+\vartheta}$, which in turn implies that $[F-\bar F_1]\nabla_z\tilde u^0_\ell\in C_p^{(1+\vartheta)/2,\delta,1+\vartheta,2}$. Consequently, it follows  by  Lemma \ref{key30} that
	\begin{align*}
	&\frac{1}{\eta_\eps}\mE\left(\int_0^T\sL_4^1\tilde u^0_\ell(s,Y_s^{0,\eps},Z_s^{0,\eps})\dif s\right)\\
	&\leq C_T\bigg[\frac{\alpha_\eps^{1+\vartheta}}{\eta_\eps}
	 +\frac{\alpha_\eps^2}{\eta_\eps\beta_\eps}\mE\left(\int_0^T\sL_3\tilde\Upsilon^0_\ell(s,X_s^\eps,Y_s^{0,\eps},Z_s^{0,\eps})\dif s\right)\\
	 &\quad+\frac{\alpha_\eps^2}{\eta_\eps^2}\mE\left(\int_0^T\sL_4^1\tilde\Upsilon^0_\ell(s,X_s^\eps,Y_s^{0,\eps},Z_s^{0,\eps})\dif s\right)\bigg].
	\end{align*}
Thus we arrive at
	\begin{align*}
	\cR^0_\ell(\eps)&\leq C_T\bigg[\frac{\alpha_\eps^{1+\vartheta}}{\eta_\eps}+\mE\left(\int_0^T\!\Big[\p_s+\sL_1+\sL_5^{1,\eps}\Big]\tilde u^0_\ell(s,Y_s^{0,\eps},Z_s^{0,\eps})\dif s\right)\\
	 &\quad+\frac{\alpha_\eps^2}{\eta_\eps\beta_\eps}\mE\left(\int_0^T\sL_3\tilde\Upsilon^0_\ell(s,X_s^\eps,Y_s^{0,\eps},Z_s^{0,\eps})\dif s\right)\\
	 &\quad+\frac{\alpha_\eps^2}{\eta_\eps^2}\mE\left(\int_0^T\sL_4^1\tilde\Upsilon^0_\ell(s,X_s^\eps,Y_s^{0,\eps},Z_s^{0,\eps})\dif s\right)\bigg]=:C_T\bigg[\frac{\alpha_\eps^{1+\vartheta}}{\eta_\eps}+\sum_{i=1}^3\sC_i(\eps)\bigg].
	\end{align*}
	Below, we consider $\ell=1,2,3$ separately, which correspond to Regime 0.1-Regime 0.3 in the conclusion of Theorem \ref{main3}, respectively.
	
	\vspace{1mm}
	\noindent{\bf Case $\ell=1$.}  (Regime 0.1 in (\ref{regime4})). Note that we have
	$$
	 \frac{\alpha_\eps^2}{\eta_\eps\beta_\eps}=1\quad\text{and}\quad\frac{\alpha_\eps^2}{\eta_\eps^2} =\frac{\beta_\eps^2}{\alpha_\eps^2}\to0\quad\text{as}\quad\eps\to0.
	$$
	Thus we deduce
	$$
	\sC_3(\eps)\leq C_T\frac{\alpha_\eps^2}{\eta_\eps^2}\mE\left(\int_0^T\big(1+|X_s^\eps|^{2m}\big)\dif s\right)\leq C_T\frac{\beta_\eps^2}{\alpha_\eps^2}.
	$$
	Furthermore, recall that we have $\bar\cL_1^0=\bar \sL_3^0+\bar \sL^1_5$ in this case with $\bar \sL_3^0$ and $\bar \sL^1_5$ given by (\ref{LL03}) and (\ref{LL05}), respectively. As a result,
	$$
	\p_t\tilde u^0_1+\big(\bar\sL_1+\bar \sL_3^0+\bar \sL^1_5\big)\tilde u^0_1=0.
	$$
	This in turn yields that
	\begin{align*}
	\sC_1(\eps)+\sC_2(\eps)&= \mE\left(\int_0^T\!\big(\sL_1-\bar\sL_1\big)\tilde u^0_1(s,Y_s^{0,\eps},Z_s^{0,\eps})\dif s\right)\\
	 &\quad+\mE\left(\int_0^T\sL_3\tilde\Upsilon^0_1(s,X_s^\eps,Y_s^{0,\eps},Z_s^{0,\eps})-\bar\sL_3^0\tilde u^0_1(s,Y_s^{0,\eps},Z_s^{0,\eps})\dif s\right)\\
	&\quad+ \mE\left(\int_0^T\!\big(\sL_5^{1,\eps}-\bar \sL^1_5\big)\tilde u^0_1(s,Y_s^{0,\eps},Z_s^{0,\eps})\dif s\right)=:\sum_{i=1}^3\sU_i(\eps).
	\end{align*}
	Note  that the function
	$\big[\sL_1(t,x,y)-\bar\sL_1(t,y)\big]\tilde u^0_1(t,y,z)$ satisfies the centering condition (\ref{cen30}) and belongs to $C_p^{\vartheta/2,\delta,\vartheta,2}$. As a direct result of the estimate (\ref{pp0}), we have
	$$
	\sU_1(\eps)\leq C_T\Big(\alpha_\eps^{\vartheta}+\frac{\alpha_\eps^2}{\beta_\eps}+\frac{\alpha_\eps^2}{\eta_\eps}\Big).
	$$
	Similarly, by definition,
	\begin{align*}
	&\sL_3\tilde\Upsilon^0_1(t,x,y,z)-\bar\sL_3^0\tilde u^0_1(t,y,z)\\
	&=\big[c(x,y)\nabla_x\up(t,x,y)-\overline{c\cdot\nabla_x\up}(t,y)\big]\nabla_z\tilde u^0_1(t,y,z),
	\end{align*}
	which satisfies (\ref{cen30}) and belongs to $C_p^{\vartheta/2,\delta,\vartheta,2}$ due to the assumption that $c\in C_p^{\delta,\vartheta}$. Consequently, we also have
	$$
	\sU_2(\eps)\leq C_T\Big(\alpha_\eps^{\vartheta}+\frac{\alpha_\eps^2}{\beta_\eps}+\frac{\alpha_\eps^2}{\eta_\eps}\Big).
	$$
	Finally, by the fact that $\bar F(t,\cdot), G(t, \cdot)\in C_b^{1+\vartheta}$ and the mean value theorem, we deduce that
	\begin{align*}
	\sU_3(\eps)&\leq\mE\!\left(\int_0^T\!\!\bigg[\frac{\bar F_1(s,Y_s^{0,\eps})-\bar F_1(s,\bar Y_s^1)}{\eta_\eps}-\nabla_y\bar F_1(s,Y_s^{0,\eps})Z_s^{0,\eps}\bigg]\nabla_z\tilde u^0_1(s,Y_s^{0,\eps},Z_s^{0,\eps})\dif s\!\right)\\
	&\quad+\mE\bigg(\int_0^T\!\!\bigg[\frac{[G(s,Y_s^{0,\eps})-G(s,\bar Y_s^1)][G(s,Y_s^{0,\eps})-G(s,\bar Y_s^1)]^*}{\eta_\eps^2}\\
	 &\qquad\qquad\quad-\big[\nabla_yG(s,Y_s^{0,\eps})Z_s^{0,\eps}\big]\big[\nabla_yG(s,Y_s^{0,\eps})Z_s^{0,\eps}\big]^*\bigg]\nabla_z^2\tilde u^0_1(s,Y_s^{0,\eps},Z_s^{0,\eps})\dif s\bigg)\\
	 &\quad+\mE\bigg(\int_0^TG(s,Y_s^{0,\eps})\bigg[\frac{[G(s,Y_s^{0,\eps})-G(s,\bar Y_s^1)]}{\eta_\eps}\\
	 &\qquad\qquad\quad-\nabla_yG(s,Y_s^{0,\eps})Z_s^{0,\eps}\bigg]\nabla_y\nabla_z\tilde u^0_1(s,Y_s^{0,\eps},Z_s^{0,\eps})\dif s\bigg)\\
	&\leq C_T\mE\left(\int_0^T\big|Y_s^{0,\eps}-\bar Y_s^1\big|^\vartheta \big(1+|Z_s^{0,\eps}|^2\big)\dif s\right)\leq C_T\Big(\alpha_\eps+\frac{\alpha_\eps^{2}}{\beta_\eps}\Big)^\vartheta,
	\end{align*}
	where in the last inequality we also used H\"older's inequality and estimate (\ref{stro1}).
Based on the above estimates, we arrive at
	\begin{align*}
	\cR^0_1(\eps)\leq C_T\Big(\frac{\alpha_\eps^{1+\vartheta}}{\eta_\eps}+\frac{\beta_\eps^2}{\alpha_\eps^2}+\alpha_\eps^{\vartheta}+\frac{\alpha_\eps^2}{\beta_\eps}+\frac{\alpha_\eps^{2\vartheta}}{\beta_\eps^\vartheta}\Big)
	 \leq C_T\Big(\frac{\beta_\eps^2}{\alpha_\eps^2}+\frac{\alpha_\eps^{2\vartheta}}{\beta_\eps^\vartheta}\Big).
	\end{align*}

	\noindent{\bf Case $\ell=2$.} (Regime 0.2 in (\ref{regime4})). Note that we have
	$$
	 \frac{\alpha_\eps^2}{\eta_\eps^2}=1\quad\text{and}\quad\frac{\alpha_\eps^2}{\eta_\eps\beta_\eps}=\frac{\alpha_\eps}{\beta_\eps}\to0\quad\text{as}\quad\eps\to0.
	$$
	Thus by the assumption that $c\in L^\infty_p$,
	$$
	\sC_2(\eps)\leq C_T\frac{\alpha_\eps^2}{\eta_\eps\beta_\eps}\mE\left(\int_0^T\big(1+|X_s^\eps|^{2m}\big)\dif s\right)\leq C_T\frac{\alpha_\eps}{\beta_\eps}.
	$$
	Furthermore, recall that we have $\bar\cL^0_2=\bar \sL_4^1+\bar \sL^1_5$ in this case with $\bar\sL_4^1$ given by (\ref{LL14}). It follows that
	$$
	\p_t\tilde u^0_2+\big(\bar\sL_1+\bar\sL_4^1+\bar\sL^1_5\big)\tilde u^0_2=0.
	$$
	Consequently, we  get
	\begin{align*}
	\sC_1(\eps)+\sC_3(\eps)&\leq \mE\left(\int_0^T\!\big(\sL_1-\bar\sL_1\big)\tilde u^0_2(s,Y_s^{0,\eps},Z_s^{0,\eps})\dif s\right)\\
	 &+\mE\left(\int_0^T\sL_4^1\tilde\Upsilon^0_2(s,X_s^\eps,Y_s^{0,\eps},Z_s^{0,\eps})-\bar\sL_4^1\tilde u^0_2(s,Y_s^{0,\eps},Z_s^{0,\eps})\dif s\right)\\
	&+ \mE\left(\int_0^T\!\big(\sL_5^{1,\eps}-\bar \sL^1_5\big)\tilde u^0_2(s,Y_s^{0,\eps},Z_s^{0,\eps})\dif s\right).
	\end{align*}
	Note that by definition,
\begin{align*}
	&\sL_4^1\tilde\Upsilon^0_2(t,x,y,z)-\bar\sL_4^1\tilde u^0_2(t,y,z)\\
	&=\Big[\tilde F(t,x,y)\cdot\up^*(t,x,y)-\overline{\tilde F\cdot\up^*}(t,y)\Big]\nabla_z^2\tilde u^0_2(t,y,z).
	\end{align*}
	Following exactly the same arguments as above, we   get
	\begin{align*}
	\sC_1(\eps)+\sC_3(\eps)\leq C_T\Big(\alpha_\eps^\vartheta+\frac{\alpha_\eps^2}{\eta_\eps}+\frac{\alpha_\eps^{2^\vartheta}}{\beta_\eps^\vartheta}\Big)\leq C_T\Big(\alpha_\eps+\frac{\alpha_\eps^{2}}{\beta_\eps}\Big)^\vartheta.
	\end{align*}
	Based on the above estimates, we arrive at
	\begin{align*}
	\cR^0_2(\eps)\leq C_T\Big(\alpha_\eps^\vartheta+\frac{\alpha_\eps^{2\vartheta}}{\beta_\eps^\vartheta}+\frac{\alpha_\eps}{\beta_\eps}\Big)\leq C_T\Big(\alpha_\eps^\vartheta+\frac{\alpha_\eps}{\beta_\eps}\Big).
	\end{align*}

	\noindent{\bf Case $\ell=3$.}  (Regime 0.3 in (\ref{regime4})).  In this case,
	we have
	$$
	\p_t\tilde u^0_3+\big(\bar\sL_1+\bar\sL^1_5+\bar\sL_3^0+\bar\sL_4^1\big)\tilde u^0_3=0.
	$$
	Thus we have
	\begin{align*}
	\sum_{i=1}^3\sC_i(\eps)&\leq \mE\left(\int_0^T\!\big(\sL_1-\bar\sL_1\big)\tilde u^0_3(s,Y_s^{0,\eps},Z_s^{0,\eps})\dif s\right)\\
	 &\quad+\mE\left(\int_0^T\sL_3\tilde\Upsilon^0_3(s,X_s^\eps,Y_s^{0,\eps},Z_s^{0,\eps})-\bar\sL_3^0\tilde u^0_3(s,Y_s^\eps,Z_s^{0,\eps})\dif s\right)\\
	 &\quad+\mE\left(\int_0^T\sL_4^1\tilde\Upsilon^0_3(s,X_s^\eps,Y_s^{0,\eps},Z_s^{0,\eps})-\bar\sL_4^1\tilde u^0_3(s,Y_s^\eps,Z_s^{0,\eps})\dif s\right)\\
		&\quad+ \mE\left(\int_0^T\!\big(\sL_5^{1,\eps}-\bar \sL^1_5\big)\tilde u^0_3(s,Y_s^{0,\eps},Z_s^{0,\eps})\dif s\right).
	\end{align*}
	By combining the above two cases we  deduce that
	\begin{align*}
	\cR^0_3(\eps)\leq C_T\Big(\alpha_\eps^\vartheta+\frac{\alpha_\eps^{2\vartheta}}{\beta_\eps^\vartheta}+\frac{\alpha_\eps^2}{\eta_\eps}\Big)\leq C_T\alpha_\eps^\vartheta,
	\end{align*}
	and the whole proof is finished.
\end{proof}

\section{Central limit theorem with homogenization}

In this section, we study the functional central limit theorem for SDE (\ref{sde0}) by following the same procedure as in Section 5. We  first derive some weak type fluctuation estimates in Subsection 6.1. Then we prove Theorem \ref{main2} and Theorem \ref{main4} in Subsection 6.2 and Subsection 6.3, respectively. We shall mainly focus on the  differences in regard to the proof of Theorem \ref{main3}.

\subsection{Weak fluctuation estimate (ii)}

Recall that depending on Regime 1 and Regime 2 described in (\ref{regime}), the slow process $Y_t^\eps$ converges strongly to $\bar Y_t^1$ and $\bar Y_t^2$, respectively. As before, to prove the weak convergence of $Z_t^{k,\eps}$ to $\bar Z^k_{\ell,t}$ ($k=1,2$ and $\ell=1,2,3$), we will  view the process $(X^{\eps}_t, Y^{\eps}_t, \bar Y^k_t, Z_t^{k,\eps}, \bar Z^k_{\ell,t})$ as a whole system, i.e., we consider
\begin{equation}\label{sde7}
\left\{ \begin{aligned}
&\dif X^{\eps}_t =\alpha_\eps^{-2}b(X^{\eps}_t,Y^{\eps}_t)\dif t+\beta_\eps^{-1}c(X^{\eps}_t,Y^{\eps}_t)\dif t+\alpha_\eps^{-1}\sigma(X^{\eps}_t,Y_t^\eps)\dif W^{1}_t,\quad X_0=x,\\
&\dif Y^{\eps}_t =F(t,X^{\eps}_t, Y^{\eps}_t)\dif t+\gamma_\eps^{-1}H(t,X_t^\eps,Y_t^\eps)\dif t+G(t,Y_t^\eps)\dif W^{2}_t,\quad\qquad\! Y_0=y,\\
&\dif \bar Y^k_t=\bar F_k(t,\bar Y^k_t)\dif t+G(t,\bar Y^k_t)\dif W^2_t,\quad\qquad\quad\qquad\quad\qquad\quad\qquad\quad\, \bar Y_0^k=y,\\
&\dif Z_t^{k,\eps}=\frac{1}{\eta_\eps}\Big[F(t,X_t^\eps,Y_t^\eps)-\bar F_k(t,\bar Y^k_t)\Big]\dif t+\frac{1}{\eta_\eps\gamma_\eps}H(t,X_t^\eps,Y_t^\eps)\dif t\\
&\quad\quad\quad\,\,+\frac{1}{\eta_\eps}\Big[G(t,Y_t^\eps)-G(t,\bar Y^k_t)\Big]\dif W_t^2,\quad\qquad\quad\qquad\quad\quad\qquad\quad\!\!  Z_0^{k,\eps}=0,\\
&\dif \bar Z^k_{\ell,t}=\varrho^k_\ell(t,\bar Y_t^k,\bar Z^k_{\ell,t})\dif t+\chi_\ell^k(t,\bar Y_t^k,\bar Z^k_{\ell,t})\dif W_t^2 \quad\qquad\quad\qquad\quad\quad  \quad \bar Z_{\ell, 0}^{k}=0,
\end{aligned} \right.
\end{equation}
where for $k=1, 2$ and $\ell=1,2,3$, $\eta_\eps$ are given by (\ref{regime2}) and (\ref{regime3}),  $\varrho^k_\ell$ and $\chi_\ell^k$ denote  the  drift and diffusion coefficients for $\bar Z^1_{\ell,t}$ in SDE (\ref{sdez1}) and $\bar Z^2_{\ell,t}$ in SDE (\ref{sdez2}), respectively.

Note that by definition,
\begin{align*}
\dif Z_t^{1,\eps}
&=\frac{1}{\eta_\eps}\Big[F(t,X_t^\eps,Y_t^\eps)-\bar F_1(t,Y_t^\eps)\Big]\dif t+\frac{1}{\eta_\eps\gamma_\eps}H(t,X_t^\eps,Y_t^\eps)\dif t\\
&\quad+\bigg(\frac{1}{\eta_\eps}\Big[\bar F_1(t,Y_t^\eps)-\bar F_1(t,\bar Y^1_t)\Big]\dif t+\frac{1}{\eta_\eps}\Big[G(t,Y_t^\eps)-G(t,\bar Y^1_t)\Big]\dif W_t^2\bigg),
\end{align*}
and since $\bar F_2(t,y)=\bar F_1(t,y)+\overline{c\cdot\nabla_x\Phi}(t,y)$, we have
\begin{align*}
\dif Z_t^{2,\eps}
&=\frac{1}{\eta_\eps}\Big[F(t,X_t^\eps,Y_t^\eps)-\bar F_1(t,Y_t^\eps)\Big]\dif t\\
&\quad+\bigg(\frac{1}{\eta_\eps\gamma_\eps}H(t,X_t^\eps,Y_t^\eps)\dif t-\frac{1}{\eta_\eps}\overline{c\cdot\nabla_x\Phi}(t,Y_t^\eps)\dif t\bigg)\\
&\quad+\bigg(\frac{1}{\eta_\eps}\Big[\bar F_2(t,Y_t^\eps)-\bar F_2(t,\bar Y^2_t)\Big]\dif t+\frac{1}{\eta_\eps}\Big[G(t,Y_t^\eps)-G(t,\bar Y^2_t)\Big]\dif W_t^2\bigg).
\end{align*}
Recall that $\bar\sL_3^1$, $\sL_4^1$ and $\sL_5^{k,\eps}$ are defined by (\ref{LL13}), (\ref{LL104}) and (\ref{LL2}), respectively. We define
\begin{align}\label{LL4}
\sL^2_4:=\sL^2_4(t,x,y,z):=\sum_{i=1}^{d_2}H^i(t,x,y)\frac{\p}{\p z_i}.
\end{align}

Given a function $f(t,x,y,z)$ satisfying  the centering condition (\ref{cen30}),
 let  $\Phi^f(t,x,y,z)$ be
the solution to the  Poisson equation (\ref{xyz0}).
We have the following estimate for the fluctuations of $f(s,X_s^\eps,Y_s^\eps,Z_s^{k,\eps})$ over $[0,t]$.

\bl\label{key31}
Let {\bf (A$_\sigma$)}, {\bf (A$_b$)}, (\ref{ac})
 hold and $0<\delta,\vartheta\leq 2$. Assume that

\vspace{1mm}
\noindent(i) (Regime 1)  $b, \sigma\in C_b^{\delta,1\vee\vartheta}$,  $F\in C_p^{0,\delta,1}$, $G\in C_b^{0,1}$ and $c, H\in L^\infty_p$;\\
(ii) (Regime 2)  $b, \sigma\in C_b^{\delta,1\vee\vartheta}$, $c\in C_p^{0,1}$, $F, H\in C_p^{0,\delta,1}$ and $G\in C_b^{0,1}$.\\
Then for every $f\in C_p^{\vartheta/2,\delta,\vartheta,2}$  satisfying (\ref{cen30}), we have for $k=1,2$,
\begin{align}\label{we2}
\begin{split}
\mE\left(\int_0^tf(s,X_s^\eps,Y_s^\eps,Z_s^{k,\eps})\dif s\right)&\leq C_t\Big(\alpha_\eps^{\vartheta}+\frac{\alpha_\eps^2}{\eta_\eps}+\frac{\alpha_\eps^{1+(\vartheta\wedge1)}}{\gamma_\eps}\Big)\\
&\quad+\frac{\alpha_\eps^2}{\beta_\eps}\mE\left(\int_0^T\sL_3\Phi^f(s,X_s^\eps,Y_s^\eps,Z_s^{k,\eps})\dif s\right)\\
&\quad+\frac{\alpha_\eps^2}{\eta_\eps\gamma_\eps}\mE\left(\int_0^T\sL_4^2\Phi^f(s,X_s^\eps,Y_s^\eps,Z_s^{k,\eps})\dif s\right).
\end{split}
\end{align}
where $C_t>0$ is a constant independent of $\delta, \eps$.
\el

\br
Compared with Lemma \ref{key30}, the term involving  $\sL_4^1\Phi^f$ is replaced by $\sL_4^2\Phi^f$, since it is of lower order  now. As before, we  get
\begin{align*}
\mE\left(\int_0^T\sL_3\Phi^f(s,X_s^\eps,Y_s^\eps,Z_s^{k,\eps})\dif s\right)+\mE\left(\int_0^T\sL_4^2\Phi^f(s,X_s^\eps,Y_s^\eps,Z_s^{k,\eps})\dif s\right)<\infty.
\end{align*}
Thus we also have
\begin{align}\label{pp}
\mE\left(\int_0^Tf(s,X_s^\eps,Y_s^\eps,Z_s^{k,\eps})\dif s\right)&\leq C_T\Big(\alpha_\eps^{\vartheta}+\frac{\alpha_\eps^2}{\beta_\eps}+\frac{\alpha_\eps^{1+(\vartheta\wedge1)}}{\gamma_\eps}+\frac{\alpha_\eps^2}{\eta_\eps\gamma_\eps}\Big).
\end{align}
The homogenization effects of the last two terms on the right hand side of (\ref{we2}) will appear in the study of the functional central limit theorems.
\er

\begin{proof}
The proof follows by the same arguments as in the proof of Lemma \ref{key30}. We provide some details here in order to make  clear which parts should be the leading terms. We only prove (\ref{we2}) for $k=2$ (Regime 2), the case $k=1$ (Regime 1) can be proved similarly and is even easier since the operator $\bar\sL_3^1$ is not involved. Let $\Phi^f_{n}$ be the mollifyer of $\Phi^f$ defined as in (\ref{fn}).
Using It\^o's formula for SDE (\ref{sde7}), we have
\begin{align*}
&\Phi^f_{n}(t,X_t^\eps,Y_t^\eps,Z_t^{2,\eps})=\Phi^f_{n}(0,x,y,0) +\frac{1}{\alpha_\eps}\hat M_{n}^1(t)+\hat M_n^2(t)+\frac{1}{\eta_\eps}\hat M_{n}^3(t)\\
&\quad\quad+\int_0^t\big(\p_s+\sL_1+\sL_5^{2,\eps}\big)\Phi^f_{n}(s,X_s^\eps,Y_s^\eps,Z_s^{2,\eps})\dif s\\
&\quad\quad+\frac{1}{\alpha_\eps^2}\int_0^t\sL_0\Phi^f_{n}(s,X_s^\eps,Y_s^\eps,Z_s^{2,\eps})\dif s+\frac{1}{\beta_\eps}\int_0^t\sL_3\Phi^f_{n}(s,X_s^\eps,Y_s^\eps,Z_s^{2,\eps})\dif s\\
&\quad\quad+\frac{1}{\gamma_\eps}\int_0^t\sL_2\Phi^f_{n}(s,X_s^\eps,Y_s^\eps,Z_s^{2,\eps})\dif s+\frac{1}{\eta_\eps}\int_0^t\sL_4^1\Phi^f_{n}(s,X_s^\eps,Y_s^\eps,Z_s^{2,\eps})\dif s\\
&\quad\quad+\frac{1}{\eta_\eps\gamma_\eps}\int_0^t\sL_4^2\Phi^f_{n}(s,X_s^\eps,Y_s^\eps,Z_s^{2,\eps})\dif s-\frac{1}{\eta_\eps}\int_0^t\bar\sL_3^1\Phi^f_{n}(s,X_s^\eps,Y_s^\eps,Z_s^{2,\eps})\dif s,
\end{align*}
where  $\sL_4^2$ is given   by  (\ref{LL4}), $\hat M_{n}^1(t)$, $\hat M_{n}^2(t)$ and $\hat M_{n}^3(t)$ are martingales defined by
\begin{align*}
\hat M_{n}^1(t)&:=\int_0^t\sigma(X_s^\eps,Y_s^\eps)\nabla_x\Phi^f_n(s,X_s^\eps,Y_s^\eps,Z_s^{2,\eps})\dif W^1_s,\\
\hat M_{n}^2(t)&:=\int_0^tG(s,Y_s^\eps)\nabla_y\Phi^f_n(s,X_s^\eps,Y_s^\eps,Z_s^{2,\eps})\dif W^2_s,\\
\hat M_{n}^3(t)&:=\int_0^t\big[G(s,Y_s^\eps)-G(s,\bar Y_s^2)\big]\nabla_z\Phi^f_n(s,X_s^\eps,Y_s^\eps,Z_s^{2,\eps})\dif W^2_s.
\end{align*}
Taking expectation and in view of (\ref{xyz0}), we get
\begin{align*}
\hat\cV(\eps)&:=\mE\!\left(\int_0^tf(s,X_s^\eps,Y_s^\eps,Z_s^{2,\eps})\dif s\right)\!=\alpha_\eps^2\mE\big[\Phi^f_{n}(0,x,y,0)-\Phi^f_{n}(t,X_t^\eps,Y_t^\eps,Z_t^{2,\eps})\big]\\
&+\alpha_\eps^2\mE\left(\int_0^T\big(\p_s+\sL_1+\sL_5^{2,\eps}\big)\Phi^f_{n}(s,X_s^\eps,Y_s^\eps,Z_s^{2,\eps})\dif s\right)\\
&+\mE\left(\int_0^T\big[\sL_0\Phi^f_{n}-\sL_0\Phi^f\big](s,X_s^\eps,Y_s^\eps,Z_s^{2,\eps})\dif s\right)\\
&+\frac{\alpha_\eps^2}{\gamma_\eps}\mE\left(\int_0^T\sL_2\Phi^f_{n}(s,X_s^\eps,Y_s^\eps,Z_s^{2,\eps})\dif s\right)\\
&+\frac{\alpha_\eps^2}{\eta_\eps}\mE\left(\int_0^T\sL_4^1\Phi^f_{n}(s,X_s^\eps,Y_s^\eps,Z_s^{2,\eps})\dif s\!\right)\!-\!\frac{\alpha_\eps^2}{\eta_\eps}\mE\left(\int_0^T\bar\sL_3^1\Phi^f_{n}(s,X_s^\eps,Y_s^\eps,Z_s^{2,\eps})\dif s\right)\\
&+\bigg[\frac{\alpha_\eps^2}{\beta_\eps}\mE\left(\int_0^T\big(\sL_3\Phi^f_{n}-\sL_3\Phi^f\big)(s,X_s^\eps,Y_s^\eps,Z_s^{2,\eps})\dif s\right)\\
&\quad+\frac{\alpha_\eps^2}{\eta_\eps\gamma_\eps}\mE\left(\int_0^T\big(\sL_4^2\Phi^f_{n}-\sL_4^2\Phi^f\big)(s,X_s^\eps,Y_s^\eps,Z_s^{2,\eps})\dif s\right)\bigg]\\
&+\bigg[\frac{\alpha_\eps^2}{\beta_\eps}\mE\left(\int_0^T\sL_3\Phi^f(s,X_s^\eps,Y_s^\eps,Z_s^{2,\eps})\dif s\right)\\
&\quad+\frac{\alpha_\eps^2}{\eta_\eps\gamma_\eps}\mE\left(\int_0^T\sL_4^2\Phi^f(s,X_s^\eps,Y_s^\eps,Z_s^{2,\eps})\dif s\right)\bigg]=:\sum_{i=1}^8\hat\cV_i(\eps).
\end{align*}
Now, due to assumptions that $b,\sigma\in C_b^{\delta,1\vee\vartheta}$,  $c\in C_p^{0,1}$, $H\in C_p^{0,\delta,1}$ and  by Lemma \ref{coe},  we have  $\bar F_2\in C_b^{0,1}$.
Following the same argument as in (\ref{52}), we can  deduce that
$$
\mE\left(\int_0^T\sL_5^{2,\eps}\Phi^f_{n}(s,X_s^\eps,Y_s^\eps,Z_s^{2,\eps})\dif s\right)\leq C_T n^{1-(\vartheta\wedge1)}.
$$
As a result, we   further have
	$$
\sum_{i=1}^7\cV_i(\eps)\leq C_T\Big(\alpha_\eps^2n^{2-\vartheta}+n^{-\vartheta}+\frac{\alpha_\eps^2}{\gamma_\eps}n^{1-(\vartheta\wedge1)}+\frac{\alpha_\eps^2}{\eta_\eps}+\frac{\alpha_\eps^2}{\beta_\eps}n^{-\vartheta}+\frac{\alpha_\eps^2}{\eta_\eps\gamma_\eps}n^{-\vartheta}\Big).
$$
Taking $n=\alpha_\eps^{-1}$ and noticing that $\alpha_\eps^2/(\eta_\eps\gamma_\eps)\to0$ as $\eps\to0$, we get the desired result.
\end{proof}

\subsection{Proof of Theorem \ref{main2}}

We are now in the position to give:

\begin{proof}[Proof of Theorem \ref{main2}]
	We concentrate on the main differences in regard to the proof of Theorem \ref{main3}.
For every $\varphi\in C_b^4(\mR^{d_2})$ and $\ell=1,2,3$, let $u^1_\ell$ be the solution to the Cauchy problem  (\ref{PDE}) on $[0,T]\times\mR^{d_2}\times\mR^{d_2}$ with $k=1$, and define
$$
\tilde u^1_\ell(t,y,z):=u^1_\ell(T-t,y,z),\quad t\in [0, T].
$$
Then  we  write for $\ell=1,2,3$,
\begin{align*}
\cR^1_\ell(\eps):=\mE\varphi(Z_{T}^{1,\eps})-\mE\varphi(\bar Z^1_{\ell,T})=\mE \tilde u^1_\ell(T,Y_T^\eps,Z_T^{1,\eps})-\tilde u^1_\ell(0,y,0).
\end{align*}
By Theorem \ref{pde} and  using It\^o's formula for SDE (\ref{sde7}), we deduce that
\begin{align*}
\tilde u^1_\ell(T,Y_T^\eps,Z_T^{1,\eps})&=\tilde u^1_\ell(0,y,0)+\int_0^T\!\Big(\p_s+\sL_1+\sL_5^{1,\eps}\Big)\tilde u^1_\ell(s,Y_s^\eps,Z_s^{1,\eps})\dif s\\
&\quad+\frac{1}{\gamma_\eps}\int_0^T\!\!\sL_2\tilde u^1_\ell(s,Y_s^\eps,Z_s^{1,\eps})\dif s+\frac{1}{\eta_\eps}\int_0^T\!\!\sL_4^1\tilde u^1_\ell(s,Y_s^\eps,Z_s^{1,\eps})\dif s\\
&\quad+\frac{1}{\eta_\eps\gamma_\eps}\int_0^T\!\!\sL_4^2\tilde u^1_\ell(s,Y_s^\eps,Z_s^{1,\eps})\dif s+\hat M^1_{\ell,T}+\frac{1}{\eta_\eps}\hat M_{\ell,T}^2,
\end{align*}
where $\hat M_{\ell,t}^1$ and $\hat M_{\ell,t}^2$ are martingales defined by
\begin{align*}
\hat M_{\ell,t}^1&:=\int_0^tG(s,Y_s^\eps)\nabla_y\tilde u^1_\ell(s,Y_s^\eps,Z_s^{1,\eps})\dif W^2_s,\\
\hat M_{\ell,t}^2&:=\int_0^t\big[G(s,Y_s^\eps)-G(s,\bar Y_s^1)\big]\nabla_z\tilde u^1_\ell(s,Y_s^\eps,Z_s^{1,\eps})\dif W^2_s.
\end{align*}
As a result, we further have
\begin{align*}
\cR^1_\ell(\eps)&\leq \mE\left(\int_0^T\!\Big[\p_s+\sL_1+\sL_5^{1,\eps}\Big]\tilde u^1_\ell(s,Y_s^\eps,Z_s^{1,\eps})\dif s\right)\\
&\quad+\frac{1}{\gamma_\eps}\mE\left(\int_0^T\!\!\sL_2\tilde u^1_\ell(s,Y_s^\eps,Z_s^{1,\eps})\dif s\right)+\frac{1}{\eta_\eps}\mE\left(\int_0^T\!\!\sL_4^1\tilde u^1_\ell(s,Y_s^\eps,Z_s^{1,\eps})\dif s\right)\\
&\quad+\frac{1}{\eta_\eps\gamma_\eps}\mE\left(\int_0^T\!\!\sL_4^2\tilde u^1_\ell(s,Y_s^\eps,Z_s^{1,\eps})\dif s\right)=:\sum_{i=1}^4\sG_i(\eps).
\end{align*}
Note that by definition,
\begin{align*}
\sG_2(\eps)=\frac{1}{\gamma_\eps}\mE\left(\int_0^TH(s,X_s^\eps,Y_s^\eps)\nabla_y\tilde u^1_\ell(s,Y_s^\eps,Z_s^{1,\eps})\dif s\right).
\end{align*}
Thanks  to the assumption {\bf (A$_H$)}, it is easily checked  that for every $\ell=1,2,3$, the functions $H(t,x,y)\nabla_y\tilde u^1_\ell(t,y,z)$ satisfy the centering condition (\ref{cen30}) and belong to $C_p^{(1+\vartheta)/2,\delta,1+\vartheta,2}$.
As a result of estimate (\ref{pp}), we  get
$$
\sG_2(\eps)\leq C_1\Big(\frac{\alpha_\eps^{1+\vartheta}}{\gamma_\eps}+\frac{\alpha_\eps^2}{\gamma_\eps\beta_\eps}+\frac{\alpha_\eps^2}{\eta_\eps\gamma_\eps^2}\Big).
$$
Similarly,  we also  have
$$
\sG_3(\eps)\leq C_2\Big(\frac{\alpha_\eps^{1+\vartheta}}{\eta_\eps}+\frac{\alpha_\eps^2}{\eta_\eps\beta_\eps}+\frac{\alpha_\eps^2}{\eta_\eps^2\gamma_\eps}\Big).
$$
Concerning the last term, note that
\begin{align*}
\sG_4(\eps)=\frac{1}{\eta_\eps\gamma_\eps}\mE\left(\int_0^TH(s,X_s^\eps,Y_s^\eps)\nabla_z\tilde u^1_\ell(s,Y_s^\eps,Z_s^{1,\eps})\dif s\right).
\end{align*}
For every $\ell=1,2,3$, the functions $H(t,x,y)\nabla_z\tilde u^1_\ell(t,y,z)$ satisfy the centering condition (\ref{cen30}). Recall that $\Phi(t,x,y)$ solves the Poisson equation (\ref{pde1}),
and define
$$
\hat\Phi^1_\ell(t,x,y,z):=\Phi(t,x,y)\nabla_z\tilde u^1_\ell(t,y,z).
$$
Then $\hat\Phi^1_\ell(t,x,y,z)$ satisfies
\begin{align*}
\sL_0(x,y)\hat\Phi^1_\ell(t,x,y,z)=-H(t,x,y)\nabla_z\tilde u^1_\ell(t,y,z).
\end{align*}
Consequently, we  use Lemma \ref{key31} to deduce that
\begin{align*}
\sG_4(\eps)&\leq C_3\Big(\frac{\alpha_\eps^{1+\vartheta}}{\eta_\eps\gamma_\eps}+\frac{\alpha_\eps^2}{\eta_\eps^2\gamma_\eps}+\frac{\alpha_\eps^2}{\eta_\eps\gamma_\eps^2}\Big)\\
&\quad+\frac{\alpha_\eps^2}{\eta_\eps\gamma_\eps\beta_\eps}\mE\left(\int_0^T\sL_3\hat\Phi^1_{\ell}(s,X_s^\eps,Y_s^\eps,Z_s^{1,\eps})\dif s\right)\\
&\quad+\frac{\alpha_\eps^2}{\eta_\eps^2\gamma_\eps^2}\mE\left(\int_0^T\sL_4^2\hat\Phi^1_{\ell}(s,X_s^\eps,Y_s^\eps,Z_s^{1,\eps})\dif s\right).
\end{align*}
Combining the above  estimates, we arrive at
\begin{align}
\cR^1_\ell(\eps)&\leq C_4\Big(\frac{\alpha_\eps^{1+\vartheta}}{\eta_\eps\gamma_\eps}+\frac{\alpha_\eps^2}{\gamma_\eps\beta_\eps}+\frac{\alpha_\eps^2}{\eta_\eps\beta_\eps}+\frac{\alpha_\eps^2}{\eta_\eps^2\gamma_\eps}+\frac{\alpha_\eps^2}{\eta_\eps\gamma_\eps^2}\Big)\no\\
&\quad+\mE\left(\int_0^T\!\Big[\p_s+\sL_1+\sL_5^{1,\eps}\Big]\tilde u^1_\ell(s,Y_s^\eps,Z_s^{1,\eps})\dif s\right)\no\\
&\quad+\frac{\alpha_\eps^2}{\eta_\eps\gamma_\eps\beta_\eps}\mE\left(\int_0^T\sL_3\hat\Phi^1_\ell(s,X_s^\eps,Y_s^\eps,Z_s^{1,\eps})\dif s\right)\no\\
&\quad+\frac{\alpha_\eps^2}{\eta_\eps^2\gamma_\eps^2}\mE\left(\int_0^T\sL_4^2\hat\Phi^1_\ell(s,X_s^\eps,Y_s^\eps,Z_s^{1,\eps})\dif s\right)\no\\
&=:C_4\Big(\frac{\alpha_\eps^{1+\vartheta}}{\eta_\eps\gamma_\eps}+\frac{\alpha_\eps^2}{\gamma_\eps\beta_\eps}+\frac{\alpha_\eps^2}{\eta_\eps\beta_\eps}+\frac{\alpha_\eps^2}{\eta_\eps^2\gamma_\eps}+\frac{\alpha_\eps^2}{\eta_\eps\gamma_\eps^2}\Big)+\sum_{i=1}^3\sV_i(\eps).\label{r1}
\end{align}
Below we consider $\ell=1,2,3$ separately, which correspond to Regime 1.1-Regime 1.3 in the conclusion of Theorem \ref{main2}, respectively.

\vspace{1mm}
\noindent{\bf Case $\ell=1$.}  (Regime 1.1 in (\ref{regime2})). Note that we have
$$
\frac{\alpha_\eps^2}{\eta_\eps\gamma_\eps\beta_\eps}=1\quad\text{and}\quad\frac{\alpha_\eps^2}{\eta_\eps^2\gamma_\eps^2}=\frac{\beta_\eps^2}{\alpha_\eps^2}\to0\quad\text{as}\quad\eps\to0.
$$
Thus we deduce that
$$
\sV_3(\eps)\leq C_T\frac{\alpha_\eps^2}{\eta_\eps^2\gamma_\eps^2}\mE\left(\int_0^T\big(1+|X_s^\eps|^{2m}\big)\dif s\right)\leq C_T\frac{\alpha_\eps^2}{\eta_\eps^2\gamma_\eps^2}
$$
Furthermore, recall that we have $\bar\cL_1^1=\bar \sL^1_5+\bar \sL_3^1$ in this case with $\bar\sL_3^1$ given by (\ref{LL13}). As a result,
$$
\p_s\tilde u^1_1+\big(\bar\sL_1+\bar \sL^1_5+\bar \sL_3^1)\tilde u^1_1=0.
$$
This in turn yields that
\begin{align*}
\sV_1(\eps)+\sV_2(\eps)&\leq \mE\left(\int_0^T\!\big(\sL_1-\bar\sL_1\big)\tilde u^1_1(s,Y_s^\eps,Z_s^{1,\eps})\dif s\right)\\
&\quad+\mE\left(\int_0^T\sL_3\hat\Phi^1_1(s,X_s^\eps,Y_s^\eps,Z_s^{1,\eps})-\bar\sL_3^1\tilde u^1_1(s,Y_s^\eps,Z_s^{1,\eps})\dif s\right)\\
&\quad+ \mE\left(\int_0^T\!\big(\sL_5^{1,\eps}-\bar \sL^1_5\big)\tilde u^1_1(s,Y_s^\eps,Z_s^{1,\eps})\dif s\right).
\end{align*}
By definition,
\begin{align*}
&\sL_3\hat\Phi^1_1(t,x,y,z)-\bar\sL_3^1\tilde u^1_1(t,y,z)\\
&=\big[c(x,y)\cdot\nabla_x\Phi(t,x,y)-\overline{c\cdot\nabla_x\Phi}(t,y)\big]\cdot\nabla_z\tilde u^1_1(t,y,z),
\end{align*}
which satisfies the centering condition (\ref{cen30}).
Using the assumption that $c\in C_p^{\delta,\vartheta}$ and exactly the same arguments as in the proof of Theorem \ref{main3} (Case $\ell=1$), we  get
\begin{align*}
\sV_1(\eps)+\sV_2(\eps)&\leq C_T\Big(\alpha_\eps^\vartheta+\frac{\alpha_\eps^2}{\beta_\eps}+\frac{\alpha_\eps^{1+\vartheta}}{\gamma_\eps}+\frac{\alpha_\eps^2}{\eta_\eps\gamma_\eps}\Big)+C_T\Big(\frac{\alpha_\eps}{\gamma_\eps}+\frac{\alpha_\eps^{2}}{\beta_\eps\gamma_\eps}\Big)^\vartheta\\
&\leq C_T\Big(\frac{\alpha_\eps^2}{\beta_\eps}+\frac{\alpha_\eps^2}{\eta_\eps\gamma_\eps}+\frac{\alpha_\eps^{2\vartheta}}{\beta_\eps^\vartheta\gamma_\eps^\vartheta}\Big).
\end{align*}
Combining the above computations with (\ref{r1}), we arrive at
\begin{align*}
\cR^1_1(\eps)&\leq C_T\Big(\frac{\alpha_\eps^2}{\gamma_\eps\beta_\eps}+\frac{\alpha_\eps^2}{\eta_\eps\beta_\eps}+\frac{\alpha_\eps^{1+\vartheta}}{\eta_\eps\gamma_\eps}+\frac{\alpha_\eps^2}{\eta_\eps^2\gamma_\eps^2}+\frac{\alpha_\eps^{2\vartheta}}{\beta_\eps^\vartheta\gamma_\eps^\vartheta}\Big)\\
&\leq C_T\Big(\gamma_\eps+\frac{\beta_\eps^2}{\alpha_\eps^2}+\frac{\alpha_\eps^{2\vartheta}}{\beta_\eps^\vartheta\gamma_\eps^\vartheta}\Big).
\end{align*}

\noindent{\bf Case $\ell=2$.}  (Regime 1.2 in (\ref{regime2})). Note that we have
$$
\frac{\alpha_\eps^2}{\eta_\eps^2\gamma_\eps^2}=1\quad\text{and}\quad\frac{\alpha_\eps^2}{\eta_\eps\gamma_\eps\beta_\eps}=\frac{\alpha_\eps}{\beta_\eps}\to0\quad\text{as}\quad\eps\to0.
$$
Thus we deduce that
$$
\sV_2(\eps)\leq C_T\frac{\alpha_\eps^2}{\eta_\eps\gamma_\eps\beta_\eps}\mE\left(\int_0^T\big(1+|X_s^\eps|^{2m}\big)\dif s\right)\leq C_T\frac{\alpha_\eps}{\beta_\eps}.
$$
Furthermore, recall that we have $\bar\cL_2^1=\bar \sL^1_5+\bar \sL_4^2$ in this case with $\bar\sL_4^2$ given by (\ref{LL24}). It then follows that
$$
\p_s\tilde u^1_2+\big(\bar\sL_1+\bar\sL^1_5+\bar\sL_4^2\big)\tilde u^1_2=0.
$$
We thus have
\begin{align*}
\sV_1(\eps)+\sV_3(\eps)&\leq \mE\left(\int_0^T\!\big(\sL_1-\bar\sL_1\big)\tilde u^1_2(s,Y_s^\eps,Z_s^{1,\eps})\dif s\right)\\
&\quad+\mE\left(\int_0^T\sL_4^2\hat\Phi^1_2(s,X_s^\eps,Y_s^\eps,Z_s^{1,\eps})-\bar\sL_4^2\tilde u^1_2(s,Y_s^\eps,Z_s^{1,\eps})\dif s\right)\\
&\quad+ \mE\left(\int_0^T\!\big(\sL_5^{1,\eps}-\bar \sL^1_5\big)\tilde u^1_2(s,Y_s^\eps,Z_s^{1,\eps})\dif s\right).
\end{align*}
By definition,
\begin{align*}
&\sL_4^2\hat\Phi^1_2(t,x,y,z)-\bar\sL_4^2\tilde u^1_2(t,y,z)\\
&=\big[H(t,x,y)\cdot\Phi^*(t,x,y)-\overline{H\cdot\Phi^*}(t,y)\big]\cdot\nabla_z^2\tilde u^1_2(t,y,z),
\end{align*}
which satisfies the centering condition (\ref{cen30}).
We can employ the  same argument used before to deduce that
\begin{align*}
\sV_1(\eps)+\sV_3(\eps)&\leq C_T\Big(\alpha_\eps^\vartheta+\frac{\alpha_\eps^2}{\beta_\eps}+\frac{\alpha_\eps^{1+\vartheta}}{\gamma_\eps}+\frac{\alpha_\eps^2}{\eta_\eps\gamma_\eps}\Big)+C_T\Big(\frac{\alpha_\eps}{\gamma_\eps}+\frac{\alpha_\eps^{2}}{\beta_\eps\gamma_\eps}\Big)^\vartheta\\
&\leq C_T\Big(\frac{\alpha_\eps^2}{\beta_\eps}+\frac{\alpha_\eps^2}{\eta_\eps\gamma_\eps}+\frac{\alpha_\eps^{\vartheta}}{\gamma_\eps^\vartheta}\Big).
\end{align*}
Combining the above computations with (\ref{r1}), we arrive at
\begin{align*}
\cR^1_2(\eps)&\leq C_T\Big(\frac{\alpha_\eps^{1+\vartheta}}{\eta_\eps\gamma_\eps}+\frac{\alpha_\eps}{\beta_\eps}+\frac{\alpha_\eps^2}{\eta_\eps^2\gamma_\eps}+\frac{\alpha_\eps^2}{\eta_\eps\gamma_\eps^2}+\frac{\alpha_\eps^{\vartheta}}{\gamma_\eps^\vartheta}\Big)\\
&\leq C_T\Big(\gamma_\eps+\frac{\alpha_\eps}{\beta_\eps}+\frac{\alpha_\eps^\vartheta}{\gamma_\eps^\vartheta}\Big).
\end{align*}

\noindent{\bf Case $\ell=3$.}  (Regime 1.3 in (\ref{regime2})).  In this case,
we have
$$
\p_s\tilde u^1_2+\big(\bar\sL_1+\bar\sL^1_5+\bar\sL_3^1+\bar\sL_4^2\big)\tilde u^1_3=0.
$$
Thus we can write
\begin{align*}
\sum_{i=1}^3\sV_i(\eps)&\leq \mE\left(\int_0^T\!\big(\sL_1-\bar\sL_1\big)\tilde u^1_3(s,Y_s^\eps,Z_s^{1,\eps})\dif s\right)\\
&\quad+\mE\left(\int_0^T\sL_3\hat\Phi^1_3(s,X_s^\eps,Y_s^\eps,Z_s^{1,\eps})-\bar\sL_3^1\tilde u^1_3(s,Y_s^\eps,Z_s^{1,\eps})\dif s\right)\\
&\quad+\mE\left(\int_0^T\sL_4^2\hat\Phi^1_3(s,X_s^\eps,Y_s^\eps,Z_s^{1,\eps})-\bar\sL_4^2\tilde u^1_3(s,Y_s^\eps,Z_s^{1,\eps})\dif s\right)\\
&\quad+ \mE\left(\int_0^T\!\big(\sL_5^{1,\eps}-\bar \sL^1_5\big)\tilde u^1_3(s,Y_s^\eps,Z_s^{1,\eps})\dif s\right).
\end{align*}
By combining the above two cases we have
\begin{align*}
\cR^1_3(\eps)\leq C_T\Big(\frac{\alpha_\eps^{1+\vartheta}}{\eta_\eps\gamma_\eps}+\frac{\alpha_\eps}{\eta_\eps}+\frac{\alpha_\eps^2}{\eta_\eps^2\gamma_\eps}+\frac{\alpha_\eps^2}{\eta_\eps\gamma_\eps^2}+\frac{\alpha_\eps^\vartheta}{\gamma_\eps^\vartheta}\Big)\leq C_T\Big(\gamma_\eps+\frac{\alpha_\eps^\vartheta}{\gamma_\eps^\vartheta}\Big).
\end{align*}
Hence the whole proof is finished.
\end{proof}

\subsection{Proof of Theorem \ref{main4}}

The central limit theorems in Regime 2 of (\ref{regime}) will be the most complicated cases since homogenization already appears in the  law of large numbers.
Now, we proceed  to give:

\begin{proof}[Proof of Theorem \ref{main4}]
For every $\varphi\in C_b^4(\mR^{d_2})$ and $\ell=1,2,3$, let $u^2_\ell$ be the solution to the Cauchy problem  (\ref{PDE}) on $[0,T]\times\mR^{d_2}\times\mR^{d_2}$ with $k=2$, and define
	$$
	\tilde u^2_\ell(t,y,z):=u^2_\ell(T-t,y,z),\quad t\in [0, T].
	$$
	Then following the same arguments as in the proof of Theorem \ref{main2}, we  have for $\ell=1,2,3$ that
	\begin{align*}
	\cR^2_\ell(\eps):=\mE\varphi(Z_T^{2,\eps})-\mE\varphi(\bar Z^2_{\ell,T})=\mE \tilde u^2_\ell(T,Y_T^\eps,Z_T^{2,\eps})-\tilde u^2_\ell(0,y,0).
	\end{align*}
Using It\^o's formula and taking expectation, we further get
	\begin{align*}
	\cR^2_\ell(\eps)&\leq \mE\left(\int_0^T\!\Big[\p_s+\sL_1+\sL_5^{2,\eps}\Big]\tilde u^2_\ell(s,Y_s^\eps,Z_s^{2,\eps})\dif s\right)\\
	&\quad+\frac{1}{\gamma_\eps}\mE\left(\int_0^T\!\!\sL_2\tilde u^2_\ell(s,Y_s^\eps,Z_s^{2,\eps})\dif s\right)+\frac{1}{\eta_\eps}\mE\left(\int_0^T\!\!\sL_4^1\tilde u^2_\ell(s,Y_s^\eps,Z_s^{2,\eps})\dif s\right)\\
	 &\quad+\left[\frac{1}{\eta_\eps\gamma_\eps}\mE\left(\int_0^T\!\!\sL_4^2\tilde u^2_\ell(s,Y_s^\eps,Z_s^{2,\eps})\dif s\right)-\frac{1}{\eta_\eps}\mE\left(\int_0^T\!\!\bar\sL_3^1\tilde u^2_\ell(s,Y_s^\eps,Z_s^{2,\eps})\dif s\right)\right]\\
&=:\sum_{i=1}^4\sN_i(\eps),
	\end{align*}
where $\bar \sL_3^1$, $\sL_4^1$  and $\sL_4^2$ are given by (\ref{LL13}), (\ref{LL104}) and (\ref{LL4}), respectively. Recall that $\Phi(t,x,y)$ is the solution of the Poisson equation (\ref{pde1}).
Define
$$
\tilde \Phi^2_\ell(t,x,y,z):=\Phi(t,x,y)\nabla_y\tilde u^2_\ell(t,y,z).
$$
Then  $\tilde \Phi^2_\ell(t,x,y,z)$ satisfies
\begin{align*}
\sL_0(x,y)\tilde \Phi^2_\ell(t,x,y,z)=-H(t,x,y)\nabla_y\tilde u^2_\ell(t,y,z).
\end{align*}
Note that for every $\ell=1,2,3$, we have $H\cdot\nabla_y\tilde u^2_\ell\in C_p^{(1+\vartheta)/2,\delta,1+\vartheta,2}$. As a result of Lemma \ref{key31} we obtain
\begin{align*}
\sN_2(\eps)&=\frac{1}{\gamma_\eps}\mE\left(\int_0^TH(s,X_s^\eps,Y_s^\eps)\nabla_y\tilde u^2_\ell(s,Y_s^\eps,Z_s^{2,\eps})\dif s\right)\\
&\leq C_T\Big(\frac{\alpha_\eps^{1+\vartheta}}{\gamma_\eps}+\frac{\alpha_\eps^2}{\gamma_\eps\eta_\eps}+\frac{\alpha_\eps^2}{\gamma_\eps^2}\Big)\\
&\quad+\frac{\alpha_\eps^2}{\gamma_\eps\beta_\eps}\mE\left(\int_0^T\sL_3\tilde \Phi^2_\ell(s,X_s^\eps,Y_s^\eps,Z_s^{2,\eps})\dif s\right)\\
&\quad+\frac{\alpha_\eps^2}{\eta_\eps\gamma_\eps^2}\mE\left(\int_0^T\sL_4^2\tilde \Phi^2_\ell(s,X_s^\eps,Y_s^\eps,Z_s^{2,\eps})\dif s\right).
\end{align*}
Since $\alpha_\eps^2=\beta_\eps\gamma_\eps$ in Regime 2, we further get
	\begin{align*}
\sN_2(\eps)&\leq C_T\Big(\frac{\alpha_\eps^{1+\vartheta}}{\gamma_\eps}+\frac{\alpha_\eps^2}{\eta_\eps\gamma_\eps^2}\Big)+\mE\left(\int_0^T\sL_3\tilde\Phi^2_\ell(s,X_s^\eps,Y_s^\eps,Z_s^{2,\eps})\dif s\right).
\end{align*}
Similarly, recall that $\Upsilon(t,x,y)$ solves the  Poisson equation (\ref{pde11}),
and define
$$
\tilde\Upsilon^2_\ell(t,x,y,z):=\Upsilon(t,x,y)\nabla_z\tilde u^2_\ell(t,y,z).
$$
Then we have for $\ell=1,2,3$,
\begin{align*}
\sL_0(x,y)\tilde\Upsilon^2_\ell(t,x,y,z)=-\big[F(t,x,y)-\bar F_1(t,y)\big]\nabla_z\tilde u^2_\ell(t,y,z).
\end{align*}
Consequently, we  use Lemma \ref{key31} again to deduce that
	\begin{align*}
\sN_3(\eps)&=\frac{1}{\eta_\eps}\mE\left(\int_0^T\big[F(s,X_s^\eps,Y_s^\eps)-\bar F_1(s,Y_s^\eps)\big]\nabla_z\tilde u^2_\ell(s,Y_s^\eps,Z_s^{2,\eps})\dif s\right)\\
&\leq C_T\Big(\frac{\alpha_\eps^{1+\vartheta}}{\eta_\eps}+\frac{\alpha_\eps^2}{\eta_\eps^2}+\frac{\alpha_\eps^2}{\eta_\eps\gamma_\eps}\Big)+\frac{\alpha_\eps^2}{\eta_\eps\beta_\eps}\mE\left(\int_0^T\sL_3\tilde\Upsilon^2_\ell(s,X_s^\eps,Y_s^\eps,Z_s^{2,\eps})\dif s\right)\\
&\quad+\frac{\alpha_\eps^2}{\eta_\eps^2\gamma_\eps}\mE\left(\int_0^T\sL_4^2\tilde \Phi^2_\ell(s,X_s^\eps,Y_s^\eps,Z_s^{2,\eps})\dif s\right)\\
&\leq C_T\Big(\frac{\alpha_\eps^{1+\vartheta}}{\eta_\eps}+\frac{\alpha_\eps^2}{\eta_\eps^2\gamma_\eps}\Big)+\frac{\alpha_\eps^2}{\eta_\eps\beta_\eps}\mE\left(\int_0^T\sL_3\tilde\Upsilon^2_\ell(s,X_s^\eps,Y_s^\eps,Z_s^{2,\eps})\dif s\right).
\end{align*}	
	Finally, define
	$$
	\hat\Phi^2_\ell(t,x,y,z):=\Phi(t,x,y)\nabla_z\tilde u^2_\ell(t,y,z).
	$$
	Then we have
	\begin{align*}
	\sL_0(x,y)\hat\Phi^2_\ell(t,x,y,z)=-H(t,x,y)\nabla_z\tilde u^2_\ell(t,y,z),
	\end{align*}
	which in turn yields that
	\begin{align*}
	 &\mE\left(\frac{1}{\eta_\eps\gamma_\eps}\int_0^TH(s,X_s^\eps,Y_s^\eps)\nabla_z\tilde u^2_\ell(s,Y_s^\eps,Z_s^{2,\eps})\dif s\right)\\
	&\leq C_T\Big(\frac{\alpha_\eps^{1+\vartheta}}{\eta_\eps\gamma_\eps}+\frac{\alpha_\eps^{2}}{\eta_\eps^2\gamma_\eps}+\frac{\alpha_\eps^{2}}{\eta_\eps\gamma^2_\eps}\Big)+\frac{\alpha_\eps^2}{\eta_\eps\gamma_\eps\beta_\eps}\mE\left(\int_0^T\sL_3\hat\Phi^2_\ell(s,X_s^\eps,Y_s^\eps,Z_s^{2,\eps})\dif s\right)\\
	 &\quad+\frac{\alpha_\eps^2}{\eta_\eps^2\gamma_\eps^2}\mE\left(\int_0^T\sL_4^2\hat\Phi^2_\ell(s,X_s^\eps,Y_s^\eps,Z_s^{2,\eps})\dif s\right).
	\end{align*}
Taken into account  the definition of $\sN_4(\eps)$ and the fact that $\alpha_\eps^2=\beta_\eps\gamma_\eps$, we obtain
\begin{align*}
\sN_4(\eps)&\leq C_T\Big(\frac{\alpha_\eps^{1+\vartheta}}{\eta_\eps\gamma_\eps}+\frac{\alpha_\eps^2}{\eta_\eps^2\gamma_\eps}+\frac{\alpha_\eps^2}{\eta_\eps\gamma_\eps^2}\Big)+\frac{\alpha_\eps^2}{\eta_\eps^2\gamma_\eps^2}\mE\left(\int_0^T\sL_4^2\hat\Phi^2_\ell(s,X_s^\eps,Y_s^\eps,Z_s^{2,\eps})\dif s\right)\\
&\quad+\frac{1}{\eta_\eps}\mE\left(\int_0^T\sL_3\hat\Phi^2_\ell(s,X_s^\eps,Y_s^\eps,Z_s^{2,\eps})-\bar\sL_3^1\tilde u^2_\ell(s,Y_s^\eps,Z_s^{2,\eps})\dif s\right).
\end{align*}	
Note that by definition,
\begin{align*}
&\sL_3\hat\Phi^2_\ell(t,x,y,z)-\bar\sL_3^1\tilde u^2_\ell(t,y,z)\\
&=\big[c(x,y)\cdot\nabla_x\Phi(t,x,y)-\overline{c\cdot\nabla_x\Phi}(t,y)\big]\cdot\nabla_z\tilde u^2_\ell(t,y,z),
\end{align*}	
and recall that $\Psi$ solves the Poisson equation (\ref{p2}).
Define
	$$
	\tilde\Psi^2_\ell(t,x,y,z):=\Psi(t,x,y)\nabla_z\tilde u^2_\ell(t,y,z).
	$$	
Then we have
$$
\sL_0(x,y)\tilde\Psi^2_\ell(t,x,y,z)=-\big[c(x,y)\cdot\nabla_x\Phi(t,x,y)-\overline{c\cdot\nabla_x\Phi}(t,y)\big]\cdot\nabla_z\tilde u^2_\ell(t,y,z).
$$
By the assumption that $c\in C_p^{\delta,1+\vartheta}$ and Lemma \ref{key31}, we further have
	\begin{align*}
&\frac{1}{\eta_\eps}\mE\left(\int_0^T\sL_3\hat\Phi^2_\ell(s,X_s^\eps,Y_s^\eps,Z_s^{2,\eps})-\bar\sL_3^1\tilde u^2_\ell(s,Y_s^\eps,Z_s^{2,\eps})\dif s\right)\\
&\leq C_T\Big(\frac{\alpha_\eps^{1+\vartheta}}{\eta_\eps}+\frac{\alpha_\eps^2}{\eta_\eps^2\gamma_\eps}+\frac{\alpha_\eps^2}{\eta_\eps\gamma_\eps}\Big)+\frac{\alpha_\eps^2}{\eta_\eps\beta_\eps}\mE\left(\int_0^T\sL_3\tilde\Psi^2_\ell(s,X_s^\eps,Y_s^\eps,Z_s^{2,\eps})\dif s\right).
\end{align*}	
Combining the above computations, we arrive at
	\begin{align*}
	&\cR^2_\ell(\eps)\leq C_T\Big(\frac{\alpha_\eps^{1+\vartheta}}{\eta_\eps\gamma_\eps}+\frac{\alpha_\eps^2}{\eta_\eps^2\gamma_\eps}+\frac{\alpha_\eps^2}{\eta_\eps\gamma_\eps^2}\Big)+\mE\left(\int_0^T\!\Big[\p_s+\sL_1+\sL_5^{2,\eps}\Big]\tilde u^2_\ell(s,Y_s^\eps,Z_s^{2,\eps})\dif s\right)\\
	 &\quad+\mE\left(\int_0^T\sL_3\tilde\Phi^2_\ell(s,X_s^\eps,Y_s^\eps,Z_s^{2,\eps})\dif s\right)+\frac{\alpha_\eps^2}{\eta_\eps^2\gamma_\eps^2}\mE\left(\int_0^T\sL_4^2\hat\Phi^2_\ell(s,X_s^\eps,Y_s^\eps,Z_s^{2,\eps})\dif s\right)\\
	 &\quad+\frac{\alpha_\eps^2}{\eta_\eps\beta_\eps}\mE\left(\int_0^T\sL_3\tilde\Upsilon^2_\ell(s,X_s^\eps,Y_s^\eps,Z_s^{2,\eps})\dif s\right)\\
	 &\quad+\frac{\alpha_\eps^2}{\eta_\eps\beta_\eps}\mE\left(\int_0^T\sL_3\tilde\Psi^2_\ell(s,X_s^\eps,Y_s^\eps,Z_s^{2,\eps})\dif s\right)\\
	 &\quad=:C_T\Big(\frac{\alpha_\eps^{1+\vartheta}}{\eta_\eps\gamma_\eps}+\frac{\alpha_\eps^2}{\eta_\eps^2\gamma_\eps}+\frac{\alpha_\eps^2}{\eta_\eps\gamma_\eps^2}\Big)+\sum_{i=1}^5\sS_i(\eps).
	\end{align*}
Below, we consider $\ell=1,2,3$ separately, which correspond to Regime 2.1-Regime 2.3 in the conclusion of Theorem \ref{main4}, respectively.
	
	\vspace{1mm}
	\noindent{\bf Case $\ell=1$.}  (Regime 2.1 in (\ref{regime3})). Note that we have
	$$
	 \frac{\alpha_\eps^2}{\eta_\eps\beta_\eps}=1\quad\text{and}\quad\frac{\alpha_\eps^2}{\eta_\eps^2\gamma_\eps^2}=\frac{\beta_\eps^2}{\alpha_\eps^2\gamma_\eps^2}\to0\quad\text{as}\quad\eps\to0.
	$$
	Thus we deduce that
	$$
	\sS_3(\eps)\leq C_T\frac{\alpha_\eps^2}{\eta_\eps^2\gamma_\eps^2}\mE\left(\int_0^T\big(1+|X_s^\eps|^{2m}\big) \dif s\right)\leq C_T\frac{\beta_\eps^2}{\alpha_\eps^2\gamma_\eps^2}.
	$$
	Furthermore, recall that we  have $\bar \sL_2=\bar\sL_1+\overline{c\cdot\nabla_x\Phi}(t,y)\cdot\nabla_y$ and $\bar\cL^2_1=\bar \sL^2_5+\bar\sL_3^0+\bar\sL_3^2$ in this case with $\bar \sL_3^2$ given by (\ref{LL23}). As a result,
	$$
	\p_s\tilde u^2_1+(\bar\sL_2+\bar \sL^2_5+\bar\sL_3^0+\bar\sL_3^2)\tilde u^2_1=0.
	$$
	This in turn yields that
	\begin{align*}
	&\sS_1(\eps)+\sS_2(\eps)+\sS_4(\eps)+\sS_5(\eps)\\
	&\leq \mE\left(\int_0^T\!\big(\sL_1-\bar\sL_1\big)\tilde u^2_1(s,Y_s^\eps,Z_s^{2,\eps})\dif s\right)\\
	&\quad+ \mE\left(\int_0^T\!\big(\sL_5^{2,\eps}-\bar \sL^2_5\big)\tilde u^2_1(s,Y_s^\eps,Z_s^{2,\eps})\dif s\right)\\
	&\quad+\mE\left(\int_0^T\sL_3\tilde\Phi^2_1(s,X_s^\eps,Y_s^\eps,Z_s^{2,\eps})-\overline{c\cdot\nabla_x\Phi}(s,Y_s^\eps)\cdot\nabla_y\tilde u^2_1(s,Y_s^\eps,Z_s^{2,\eps})\dif s\right)\\
	&\quad+\mE\left(\int_0^T\sL_3\tilde\Upsilon^2_1(s,X_s^\eps,Y_s^\eps,Z_s^{2,\eps})-\bar\sL_3^0\tilde u^2_1(s,Y_s^\eps,Z_s^{2,\eps})\dif s\right)\\
	&\quad+\mE\left(\int_0^T\sL_3\tilde\Psi^2_1(s,X_s^\eps,Y_s^\eps,Z_s^{2,\eps})-\bar\sL_3^2\tilde u^2_1(s,Y_s^\eps,Z_s^{2,\eps})\dif s\right).
	\end{align*}
	By definition, we have
	\begin{align*}
	&\sL_3\tilde\Phi^2_1(t,x,y,z)-\overline{c\cdot\nabla_x\Phi}(t,y)\cdot\nabla_y\tilde u^2_1(t,y,z)\\
	&=\big[c(x,y)\cdot\nabla_x\Phi(t,x,y)-\overline{c\cdot\nabla_x\Phi}(t,y)\big]\cdot\nabla_y\tilde u^2_1(t,y,z),\\
	&\sL_3\tilde\Upsilon^2_1(t,x,y,z)-\bar\sL_3^0\tilde u^2_1(t,y,z)\\
	&=\big[c(x,y)\cdot\nabla_x\up(t,x,y)-\overline{c\cdot\nabla_x\up}(t,y)\big]\cdot\nabla_z\tilde u^2_1(t,y,z),\\
	&\sL_3\tilde\Psi^2_1(t,x,y,z)-\bar\sL_3^2\tilde u^2_1(t,y,z)\\
	&=\big[c(x,y)\cdot\nabla_x\Psi(t,x,y)-\overline{c\cdot\nabla_x\Psi}(t,y)\big]\cdot\nabla_z\tilde u^2_1(t,y,z).
	\end{align*}
	Using (\ref{stro2}) and the same arguments as before, we  deduce that
	\begin{align*}
\sS_1(\eps)+\sS_2(\eps)+\sS_4(\eps)+\sS_5(\eps)&\leq C_T\Big(\alpha_\eps^{\vartheta}+\frac{\alpha_\eps^2}{\beta_\eps} +\frac{\alpha_\eps^2}{\eta_\eps\gamma_\eps}+\frac{\alpha_\eps^{\vartheta}}{\gamma_\eps^\vartheta}+\frac{\alpha_\eps^{2\vartheta}}{\beta_\eps^\vartheta}\Big)\\
 &\leq C_T\Big(\alpha_\eps^\vartheta+\frac{\alpha_\eps^2}{\eta_\eps\gamma_\eps}+\frac{\alpha_\eps^{2\vartheta}}{\beta_\eps^\vartheta}\Big).
	\end{align*}
	Combining the above computations, we arrive at
	\begin{align*}
	\cR^2_1(\eps)\leq C_T\Big(\alpha_\eps^\vartheta+\frac{\alpha_\eps^{2\vartheta}}{\beta_\eps^\vartheta}+\frac{\beta_\eps^2}{\alpha_\eps^2\gamma_\eps^2}\Big)\leq C_T\Big(\frac{\alpha_\eps^{2\vartheta}}{\beta_\eps^\vartheta}+\frac{\beta_\eps^2}{\alpha_\eps^2\gamma_\eps^2}\Big).
	\end{align*}
	
\noindent{\bf Case $\ell=2$.}  (Regime 2.1 in (\ref{regime3})).  Note that in this case, we have
	$$
	 \frac{\alpha_\eps^2}{\eta_\eps^2\gamma_\eps^2}=1\quad\text{and}\quad\frac{\alpha_\eps^2}{\eta_\eps\beta_\eps}=\frac{\alpha_\eps\gamma_\eps}{\beta_\eps}\to0\quad\text{as}\quad\eps\to0..
	$$
	Thus we deduce that
	$$
	\sS_4(\eps)+\sS_5(\eps)\leq C_T\frac{\alpha_\eps^2}{\eta_\eps\beta_\eps}\mE\left(\int_0^T\big(1+|X_s^\eps|^{2m}\big) \dif s\right)\leq C_T\frac{\alpha_\eps\gamma_\eps}{\beta_\eps}.
	$$
	Furthermore, recall that we  have $\bar\cL^2_2=\bar \sL^2_5+\bar\sL_4^2$ in this case. As a result,
	$$
	\p_s\tilde u^2_2+(\bar\sL_2+\bar \sL^2_5+\bar\sL_4^2)\tilde u^2_2=0.
	$$
	This in turn yields that
	\begin{align*}
	&\sS_1(\eps)+\sS_2(\eps)+\sS_3(\eps)\\
	&\leq \mE\left(\int_0^T\!\big(\sL_1-\bar\sL_1\big)\tilde u^2_2(s,Y_s^\eps,Z_s^{2,\eps})\dif s\right)\\
	&\quad+ \mE\left(\int_0^T\!\big(\sL_5^{2,\eps}-\bar \sL^2_5\big)\tilde u^2_2(s,Y_s^\eps,Z_s^{2,\eps})\dif s\right)\\
	&\quad+\mE\left(\int_0^T\sL_3\tilde\Phi^2_2(s,X_s^\eps,Y_s^\eps,Z_s^{2,\eps})-\overline{c\cdot\nabla_x\Phi}(s,Y_s^\eps)\cdot\nabla_y\tilde u^2_2(s,Y_s^\eps,Z_s^{2,\eps})\dif s\right)\\
	&\quad+\mE\left(\int_0^T\sL_4^2\hat\Phi^2_2(s,X_s^\eps,Y_s^\eps,Z_s^{2,\eps})-\bar\sL_4^2\tilde u^2_2(s,Y_s^\eps,Z_s^{2,\eps})\dif s\right).
	\end{align*}
	By definition, we have
	\begin{align*}
	&\sL_4^2\hat\Phi^2_2(t,x,y,z)-\bar\sL_4^2\tilde u^2_2(t,y,z)\\
	&=\big[H(t,x,y)\cdot\Phi^*(t,x,y)-\overline{H\cdot\Phi^*}(t,y)\big]\cdot\nabla^2_z\tilde u^2_2(t,y,z).
	\end{align*}
	Then we can get
	\begin{align*}
	\sS_1(\eps)+\sS_2(\eps)+\sS_3(\eps)\leq C_T\Big(\alpha_\eps^{\vartheta}+\frac{\alpha_\eps^2}{\beta_\eps} +\frac{\alpha_\eps^2}{\eta_\eps\gamma_\eps}+\frac{\alpha_\eps^{\vartheta}}{\gamma_\eps^\vartheta}+\frac{\alpha_\eps^{2\vartheta}}{\beta_\eps^\vartheta}\Big)\leq C_T\Big(\frac{\alpha_\eps^2}{\eta_\eps\gamma_\eps}+\frac{\alpha_\eps^\vartheta}{\gamma_\eps^\vartheta}\Big).
	\end{align*}
	Combining the above computations, we arrive at
	\begin{align*}
	\cR_2^2(\eps)\leq C_T\Big( \frac{\alpha_\eps^\vartheta}{\gamma_\eps^\vartheta}+\frac{\alpha_\eps\gamma_\eps}{\beta_\eps}\Big).
	\end{align*}
	
	\noindent{\bf Case $\ell=3$.}  (Regime 2.1 in (\ref{regime3})).  In this case, we have
	$$
	\p_s\tilde u^2_3+(\bar\sL_2+\bar \sL^2_5+\bar\sL_3^0+\bar\sL_3^2+\bar\sL_4^2)\tilde u^2_3=0.
	$$
	Thus we write
	\begin{align*}
	\sum_{i=1}^5\sS_i(\eps)&\leq \mE\left(\int_0^T\!\big(\sL_1-\bar\sL_1\big)\tilde u^2_3(s,Y_s^\eps,Z_s^{2,\eps})\dif s\right)\\
	&\quad+ \mE\left(\int_0^T\!\big(\sL_5^{2,\eps}-\bar \sL^2_5\big)\tilde u^2_3(s,Y_s^\eps,Z_s^{2,\eps})\dif s\right)\\
	&\quad+\mE\left(\int_0^T\sL_3\tilde\Phi^2_3(s,X_s^\eps,Y_s^\eps,Z_s^{2,\eps})-\overline{c\cdot\nabla_x\Phi}(s,Y_s^\eps)\cdot\nabla_y\tilde u^2_3(s,Y_s^\eps,Z_s^{2,\eps})\dif s\right)\\
	&\quad+\mE\left(\int_0^T\sL_3\tilde\Upsilon^2_3(s,X_s^\eps,Y_s^\eps,Z_s^{2,\eps})-\bar\sL_3^0\tilde u^2_3(s,Y_s^\eps,Z_s^{2,\eps})\dif s\right)\\
	&\quad+\mE\left(\int_0^T\sL_3\tilde\Psi^2_3(s,X_s^\eps,Y_s^\eps,Z_s^{2,\eps})-\bar\sL_3^2\tilde u^2_3(s,Y_s^\eps,Z_s^{2,\eps})\dif s\right)\\
	&\quad+\mE\left(\int_0^T\sL_4^2\hat\Phi^2_3(s,X_s^\eps,Y_s^\eps,Z_s^{2,\eps})-\bar\sL_4^2\tilde u^2_3(s,Y_s^\eps,Z_s^{2,\eps})\dif s\right).
	\end{align*}
	By combining the above two cases we deduce that
	\begin{align*}
	\cR^2_3(\eps)\leq C_T\Big(\alpha_\eps^{\vartheta}+\frac{\alpha_\eps^2}{\beta_\eps} +\frac{\alpha_\eps^2}{\eta_\eps\gamma_\eps}+\frac{\alpha_\eps^{\vartheta}}{\gamma_\eps^\vartheta}+\frac{\alpha_\eps^{2\vartheta}}{\beta_\eps^\vartheta}\Big) \leq C_T\frac{\alpha_\eps^\vartheta}{\gamma_\eps^\vartheta},
	\end{align*}
	and the whole proof is finished.	
\end{proof}

\end{document}